\documentclass[letterpaper,12pt]{article}
\usepackage{mathtools} 
\usepackage{latexsym,amssymb,amsmath,amscd,amsfonts}
\usepackage{graphicx} 
\usepackage{url}
\usepackage{color}
\usepackage{hyperref}
\hypersetup{colorlinks=true}

\newcommand{\bdot}{\boldsymbol{\cdot}}
\newcommand{\ovl}{\overline}

\newcommand{\real}{\mathbb{R}}

\newcommand{\lbilin}{\langle\langle}
\newcommand{\rbilin}{\rangle\rangle}

\newcommand\mystack[2]{\genfrac{}{}{0pt}{}{#1}{#2}}

\newcommand{\half}{{\frac{1}{2}}}
\newcommand{\thalf}{{\frac{3}{2}}}

\newcommand{\dx}{\triangle x}
\newcommand{\dy}{\triangle y}
\newcommand{\dz}{\triangle z}
\newcommand{\dt}{\triangle t}

\newcommand{\grad}{{\vec{\nabla}}}
\newcommand{\divg}{{\vec{\nabla} \bdot }}
\newcommand{\curl}{{\vec{\nabla} \times }}
\newcommand{\lap}{{\Delta}}
\newcommand{\lapv}{{\mathbf \Delta}}

\newcommand{\Nodes}{\mathcal N}
\newcommand{\Edges}{\mathcal E}
\newcommand{\Faces}{\mathcal F}
\newcommand{\Cells}{\mathcal C}
\newcommand{\StarNodes}{{{\mathcal N}^\star}}
\newcommand{\StarEdges}{{{\mathcal E}^\star}}
\newcommand{\StarFaces}{{{\mathcal F}^\star}}
\newcommand{\StarCells}{{{\mathcal C}^\star}}

\newcommand{\SpaceN}{{S_\Nodes}}
\newcommand{\SpaceE}{{V_\Edges}}
\newcommand{\SpaceF}{{V_\Faces}}
\newcommand{\SpaceC}{{S_\Cells}}
\newcommand{\StarSpaceN}{{S_{\Nodes^\star}}}
\newcommand{\StarSpaceE}{{V_{\Edges^\star}}}
\newcommand{\StarSpaceF}{{V_{\Faces^\star}}}
\newcommand{\StarSpaceC}{{S_{\Cells^\star}}}

\newcommand{\GRAD}{{\mathcal G}}
\newcommand{\CURL}{{\mathcal R}}
\newcommand{\DIVG}{{\mathcal D}}

\newcommand{\GRADstar}{{{\mathcal G}^\star}}
\newcommand{\CURLstar}{{{\mathcal R}^\star }}
\newcommand{\DIVGstar}{{{\mathcal D}^\star }}

\newcommand{\Dx}{{\Delta x}}
\newcommand{\Dy}{{\Delta y}}
\newcommand{\Dz}{{\Delta z}}

  \newcommand{\norm}[1]{\left|\left|#1\right|\right|}
\usepackage[mathscr]{euscript}
\numberwithin{figure}{section}
\numberwithin{table}{section}

\begin{document}
\title{Explicit Time Mimetic Discretiztions \thanks{ }}
\author{
    Stanly L. Steinberg \\
    Department of Mathematics and Statistics \\
    University of New Mexico, 
    Albuquerque NM 87131-1141 USA
    }

\maketitle
\newpage \clearpage
\tableofcontents
\listoftables
\listoffigures
\newpage
\begin{abstract}

This paper is part of a program to combine a staggered time and
staggered spatial discretization of continuum mechanics problems
so that any property of the continuum that is proved using vector
calculus can be proven in an analogous way for the discretized
system. It is required that the discretizations be second order
accurate and have a conserved quantity that approximates the
energy for the system and guarantees stability for a reasonable
constraint on the time step. It is also require that the discretization
is time explicit so as to avoid the solution of large system of
possibly nonlinear algebraic equations.
The well known Yee grid discretization of Maxwell's equations
is the same as discretization described here and is an early
example of using a staggered space and time grid .

Motivation of the discussion begins by studying the staggered time
or leapfrog discretization of the harmonic oscillator and use this
to introduce the modification of the energy that is conserved.
Next systems of linear equations are used to motivate the definition
of the modified energy for more complex systems of ordinary 
differential equations and then apply these ideas to the scalar wave
equation in one spatial dimension. Discretizing the three dimensional
scalar wave and Maxwell's equations shows the power of the mimetic method.
Because the spatial discretization is mimetic, the divergence of the
electric and magnetic fields are constant when there are no sources.
Using the mimetic properties the proof of this trivial and is essentially
the same as in the continuum.

Updated versions of this paper will appear at arXiv.

\end{abstract}


\newpage \clearpage
\setcounter{equation}{0}
\section{Introduction}

The first 7 section have been rewritten. Sections 8, 9, and 10 are
to be done soon.

The main goal is to show that mimetic finite difference spatial
discretizations can be combined with an explicit finite difference time
discretizations to discretize wave equations to produce stable
second-order accurate simulations and to create a general method for
extending mimetic discretizations to model inhomogeneous and anisotropic
materials. The mimetic spatial discretizations require the use of two grids
where the corners of cells in one grid are the same points as the centers
of the cells in the other grid.

The time discretizations require that the wave equation be written as a first
order system and then the leapfrog discretization is used for the system each
equation.  An early example of this type of time discretization was given by
Yee \cite{Yee1966}.  Most importantly, it is shown how to derive a discrete
quantity that is conserved by the solution of the disctrized wave equation
and also provides a second order accurate in time approximation to the energy
of the waves and is positive for a sufficiently small time step.

Mimetic spatial discretizations have been used extensively to create
simulation programs for problems in continuum mechanics, see
\cite{LipnikovMS14} and the volume \cite{Koren2014} in which this work
appeared. They have also been used to model inhomogeneous and anisotropic
materials in two dimensions \cite{HymanMorel02,HymanSS97}.  More recently
mimetic methods have been used in Magnetohydrodynamics \cite{Krauss18,Krauss16}
For a comparison of mimetic finite difference, finite volume and finite
element discretizations see \cite{BochevHyman06}.

The spatial discretization studied here are extensions of those described in
\cite{RobidouxSteinberg2011} where it is shown that mimetic discretizations
of the gradient, curl and divergence satisfy all of the important properties
of the continuum operators. For example the discrete divergence of the
discrete curl operator is identically zero and the adjoint of the discrete
gradient is the minus the discrete divergence.  To extend this work to
anisotropic materials it is critical that in three dimensions the anisotropic
properties of many important materials are describe by a $3\times3$ symmetric
positive definite matrix \cite{Nye04}, see Chapter 4 Sections 1 and 4 for the
permittivity tensor and Chapter 11 Section 5 for heat flow.  Additionally,
\cite{LipnikovMS14,BrezziLS05,HymanMSS02,HymanShashkov01} discuss the
incorporation of material properties into mimetic discretizations using
symmetric positive definite matrices.

Sections \ref{Harmonic Oscillator}, \ref{ODEs} and \ref{1D Wave} show how to
create time conserved quantities for leapfrog time discretizations for
successively more complex problems. Section \ref{Harmonic Oscillator}
discusses the harmonic oscillator, section \ref{ODEs} discusses constant
coefficient linear systems of ordinary differential equations while
Section \ref{1D Wave} introduces a staggered in time and space discretization
for the one dimensional wave equation. It is important that the conserved
quantities depend on the time step but converge second order to a conserved
quantity for the continuum problem. Simulations show that the quantity is
conserved to at least one part in $10^{14}$.

In Section \ref{Harmonic Oscillator} begins with the harmonic oscillator for
which it is well known that the implicit Crank-Nicolson discretization
conserves a discrete energy, see Appendix \ref{C-N}.  The oscillator
equations are written as a first order system and discretized using a
leapfrog scheme. Then the continuum energy or Hamiltonian for oscillator
are approximated to make a conserved quantity.  Simulations show that this
quantity is conserved to one part in $10^{15}$, see
{\tt StaggeredOscillator.m}. Keeping the conserved quantity positive gives a
time step constraint for stability that is far less restrictive than the
constraint for accuracy.

In Section \ref{ODEs} the discussion of the harmonic oscillator is extended
to finite and infinite systems of linear ordinary differential equations
that are wave equations. Such equations have conserved continuum quantities
and their leapfrog discretizations are shown to also have an approximate
conserved quantities analogous to that of the harmonic oscillator.  These
systems are set up to be analogous to the discrete systems obtained from
the mimetic discretizations of continuum wave equations.  The equations
studied are are important motivation for extending the harmonic oscillator
discretization to partial differential equations.  Again simulations show
that the approximate quantities are conserved to about 2 parts in $10^{15}$,
see {\tt SystemsODE.m}.

In Section \ref{1D Wave}, the scalar wave equation in one space dimension
is written as a system of two first order equations and discretized using
a grid staggered in space and time. A natural continuum conserved quantity
which implies that the classical energy is conserved is introduced.
The first order system is discretized using a staggered space time
discretization like those used in Yee method \cite{Yee1966}.  The discrete
conserved quantities are extended to cover this 1D case.  Simulation show
that the new quantities are conserved to one part in
$10^{14}$, see {\tt OneDWave.m}. This provides insight into how to combine
mimetic time and spatial discretizations. 

Next Sections \ref{3D Wave Equations}, \ref{Second Order D0s}, and
\ref{Wave Equations With Variable Materials} discuss three dimensional continuum
second order differential operators for anisotropic and inhomogeneous
materials that are then used to define general scalar and vector wave
equations. The main difference between the discussion here and that in
\cite{RobidouxSteinberg2011} is the introduction of variable material
properties.

Section \ref{3D Wave Equations} reviews some continuum wave equations in 3
dimensions when the material properties are constant.  The main issue is
understanding the role of the spatial dimension for distance $d$ plays in the
partial differential equations. The scalar and vector wave equations, the
elastic wave equation and Maxwell's equations are introduced and a conserved
quantity is given for each equation.

In Section \ref{Second Order D0s} continuum second order differential
operators for anisotropic and inhomogeneous materials are introduced.  The
main idea is to use the notion of a double exact sequence and diagram
chasing to define a large class of second order spatial differential operators
that can be used to define wave equations.  This idea is motivated by the
exact sequences used in differential geometry.  A knowledge of differential
forms is not required for understanding this material but can be helpful
\cite{Arapura99}.  A discrete double exact sequence is critical for the
discussion of discretizations using staggered grids.  Importantly, the
analogous discrete double exact sequence cannot be reduced to a single
sequence.  The paper \cite{Palha20141394}, Figure 9, uses a double exact
sequence that is called a De Rham complex.

Additionally, weighted inner products for scalar and vector functions are
introduced and used to define adjoint operators and to show that the second
order operators are either positive or negative.  The main difference between
the discussion here and that in \cite{RobidouxSteinberg2011} is the
introduction of variable material properties.  This discussion depends
heavily on the spatial units of the dependent variables, the differential
operators, and the material properties.

In Section \ref{Wave Equations With Variable Materials} the second order
differential operators defined in the previous section are combined with a
second time derivative to define several types of wave equations with variable
material properties.  The second order equations are written as a first
order system that has properties similar to the systems studied earlier.
This then gives an automatic definition of a conserved quantity. The
conserved quantity can be used to easily show that energy is conserved.
At the end of the section Maxwell's equations and the general elastic wave
equations examples are studied.

The material below to be revised soon!

In Section \ref{Mimetic Discretizations} primal and dual staggered grids in
3D are introduced. These grids are the same as those introduce by Yee 
\cite{Yee1966} in 1966 to discretize Maxwell's equations. Consequently there
are two types of scalar fields and two types of discrete vector fields.
The differential operators divergence, gradient and curl are discretized as
in \cite{RobidouxSteinberg2011}.  Because two grids are used there are two of
each of the first order discrete operators $\divg$, $\curl$ and $\grad$.
Additionally it is shown how to discretize the material properties.  This
section continues by defining discrete inner products and adjoint operators
critical for understanding important properties of the discrete operators.
Note that the paper \cite{MohamedHiraniSataney16} also used a dual grid
differential form method to discretize the Navier-Stokes Equations.
For an introduction to the relationship of vector calculus to
differential forms see the notes \cite{Arapura99}.

In Section \ref{Discretize Wave} the scalar wave equation and Maxwell's
equations are discretized for constant scalar material properties and
constant time the identity for matrix material properties. This is actually
a easy task. There are programs \cite{RobidouxSteinberg2011} available
for computing the action of the divergence, curl on discrete vector fields
and the gradient on discrete scalar fields. As there are two types of fields
so there are two types of each discrete operator. Additionally a with a bit
of algebra a conserved quantity similar to the ones in the previous sections
can be computed. Simulation show that the conserved quantity is constant to
one part in $10^{15}$, see {\tt ScalarWave.m}.  Also the curl of the velocity
field is constant to one part in $10^{13}$. So everything works in the case
of the scalar wave equation and Maxwell's equations with constant material
properties.

In Section \ref{E and M section} we show how to create a conserved quantity
for Maxwell's equations that is constant to one in $10^{15}$, and
also show that the divergence of the electric and magnetic fields
are constant to one part in $10^{13}$, see {\tt Maxwell.m}.

Should section 10 be an appendix?

\subsection{Notes}

For the latest, see the minisimposium at a recent SIAM meeting \cite{SIAM}.
For more information on steady state problems see \cite{LipnakovBVM14}
For an idea of the difficulties encountered in discretizing Maxwell's
equations see \cite{ChristliebTang14}.
Others have used approximate quantities to study time discretizations 
\cite{HairerLW05,NGI,EngleSD05,GansS00,Quispel08}.
It appears that most of the energy preserving methods are implicit, but by
introducing additional variables, explicit methods that conserve a modified
energy are discussed in \cite{Tao2016}.

One complexity of mimetic spatial discretizations is caused by having
primal and dual grids.  This is leads to there being a primal gradient,
curl, and divergence and dual gradient, curl and divergence. The dual
operators are labeled with a star $\star$. This complexity was already
present in the paper by Yee \cite{Yee1966} which has evolved into the FDTD
discretization method \cite{wikiFDTD}.

There are several minor problems caused by writing wave equations as a
second order differential equations or as a system of two first order
equations.  For example second order equations are not exactly equivalent
to first order system.  Additionally, for the discrete equations there are
problems in converting the initial data for the second order equations to
data for the first order equations and vice versa.  Additionally, because
the equations studied are linear, if  they conserve some quantity, they will
conserve infinitely many quantities and thus there are choices in what
conserved quantity to study. For the first order system there is a natural
{\em primitive} conserved quantity.

This paper was inspired by the papers \cite{WanBihloNave2015} and
\cite{Yee1966}.  We note that in \cite{SchuhmannWeiland2001}
(see equation (45)) the same stability constraint was found as the
one in this paper for conserving the classical energy by modifying
the discretization of Maxwell's equations.  In \cite{GaoZang13} an
implicit (ADI) method is developed that has a modified energy that
is similar to the one used here but the added term is positive while
the added term here is negative.  For a finite element approach that
produce many of the same results that as in this paper see
\cite{TaylorFournier10,BrezziBM14}.

The paper \cite{Sanderse13} gives an overview of energy conserving methods
for Navier-Stokes equations and develops some implicit Runge-Kutta methods
for doing this. The thesis \cite{Capuano15} addresses energy conservation for
turbulent flows.  For a differential forms approach to discretization
see \cite{PerotZusi14,Teixeira13} and additionally for multisympletic time
integration approach to Maxwell's equations see \cite{StToDeMa2015}.
For two dimensional problems see
\cite{ChenLL08,Perot2011,PalhaGerritsma16,Salmon07,VeigaLV15,MorinishiLCM98,AdjointHymanShashkov97, NaturalHymanShashkov97}.
The papers \cite{Tonti14, WanBihloNave2015} take a novel approach to finding
discrete models.  For a finite-element approach to vector wave equations
see Section 2.3.2 of \cite{arnoldelastic2010}.

For isotropic and homogeneous materials, simulations show that the three
dimensional scalar wave equation and Maxwell's equations without sources
the approximate energy is constant to less than one part in $10^{15}$.
Additionally, for the scalar wave equation the curl of the velocity is
constant to less than one part in $10^{13}$ and the divergence of the
electric and magnetic fields are constant to less than one part in
$10^{13}$, see {\tt ScalarWave.m} and {\tt Maxwell.m}.

For higher order mimetic methods, see \cite{2015SanchezETALAlgorithms}.

\newpage \clearpage
\setcounter{equation}{0}
\section{The Harmonic Oscillator \label{Harmonic Oscillator}}

The goal is to use the discretization of the harmonic oscillator to motivate
time discretizations of three dimensional wave equations that conserve a
discrete approximate energy.  First the second order continuum oscillator
equation and its energy are described. Then the second order equation is written
as a first order system and a conserved quantity for the system that is
essentially the energy is described.

Next the central difference approximation of the second order wave is equation
is described and a discrete conserved quantity is derived that is an
approximation of the energy. Next an staggered in time discretization of
the first order system is described and a conserved approximate energy
is derived. The two conserved quantities are essentially the same.
These quantities give a constraint on the time step for stability that
is far less restrictive than the constraint of the time step to obtain
an accurate solution.
The first section in \cite{PalhaGerritsma15} has a complementary discussion
of energy conservation for the harmonic oscillator.

Appendix \ref{C-N} reviews the implicit Crank-Nicholson discretization
of the oscillator which conserves a natural discretization of the
continuum energy. The paper \cite{WanBihloNave2015} uses a natural
discretization of the classical energy to derive a discretization of the
oscillator equation that is equivalent to the Crank-Nicholson discretization.

\subsection{The Harmonic Oscillator and Conserved Quantities}

The linear harmonic oscillator equation is given by 
\[
u'' + \omega^2 u = 0 \,,
\]
where $u = u(t)$ is a smooth function of time $t$ and $u' = du/dt$,
$u'' = d^2u/dt^2$ and $\omega>0$ is a real constant.  The total energy is a
multiple of the average of the kinetic and potential energies which is 
\[
E = \frac{(u')^2 + (\omega \, u)^2}{2} \,.
\]
This is conserved quantity because
\[
E' =
u'' \, u' + \omega^2 u \, u' = 
\left(u'' + \omega^2 u \right) \, u' = 0 \,.
\]

The oscillator equation can be written as a first order system by
introducing $v = v(t)$ and requiring
\begin{equation}\label{First Order System}
u' = \omega \, v \,,\quad v' = -\omega \, u \,.
\end{equation}
The minus sign can be put in either equation. For the system, set
\begin{equation}\label{Conserved Quantity Simple}
C = \frac{1}{2} \left( u^2 + v^2 \right)\,.
\end{equation}
This quantity is conserved because
\[
C' = u\,u' + v\,v' = u \, \omega \, v - v \, \omega \, u = 0 \,,
\]
Use \ref{First Order System} to remove $v$ to get
\begin{equation}
C  = \frac{E}{\omega^2} \,.
\end{equation}

Because the second order equation is linear with constant coefficients
the time derivatives of $u$ satisfy the same equation a $u$ and because
the system is linear with constant coefficients, the time derivatives
of $u$ and $v$ also satisfy the system thus creating an infinity of
conserved quantities.

The condition that $\omega > 0$ and not that $\omega \geq 0 $ is important
because for $\omega = 0$ the second order equation with $u(0)= 0$ and
$u'(0) = 1$ has the solution $u(t) = t$ for which the energy is unbounded.
However, for $\omega = 0$ the system only has constant a solution
which have bounded energy. So the second order equation and the
system are not consistent for $\omega = 0$. This will have some 
impact on general system of linear constant coefficient ODEs and PDEs.

\subsection{Discretizing the Second Order Oscillator Equation}

If $\Delta t > 0$ then a standard explicit discretization of the second order
oscillator equation using the discrete times
$t_n = n \, \Delta t$, $0 \leq n < \infty$ is where $n$ is and integer and
\[
\frac{u^{n+1} -2 \, u^n + u^{n-1}}{\Delta t ^2} + \omega^2 \, u^n = 0 ,
\quad n \geq 1 \,.
\]
Given the two initial conditions $u(0)$ and $u'(0)$ set
\begin{align}
u^0 & = u(0) \,  \nonumber \\
u^1 & = u(\Delta t) = u(0) + \Delta t \, u'(0).
\nonumber
\end{align}
The discrete equation is then 
\[
u^{n+1} = (2 -(\omega \Delta t)^2) u^n - u^{n-1} , \quad n \geq 1 \,.
\]

\begin{figure}
\begin{center}
   \includegraphics[width = 4.00in,trim = 50 0 90 350,clip]{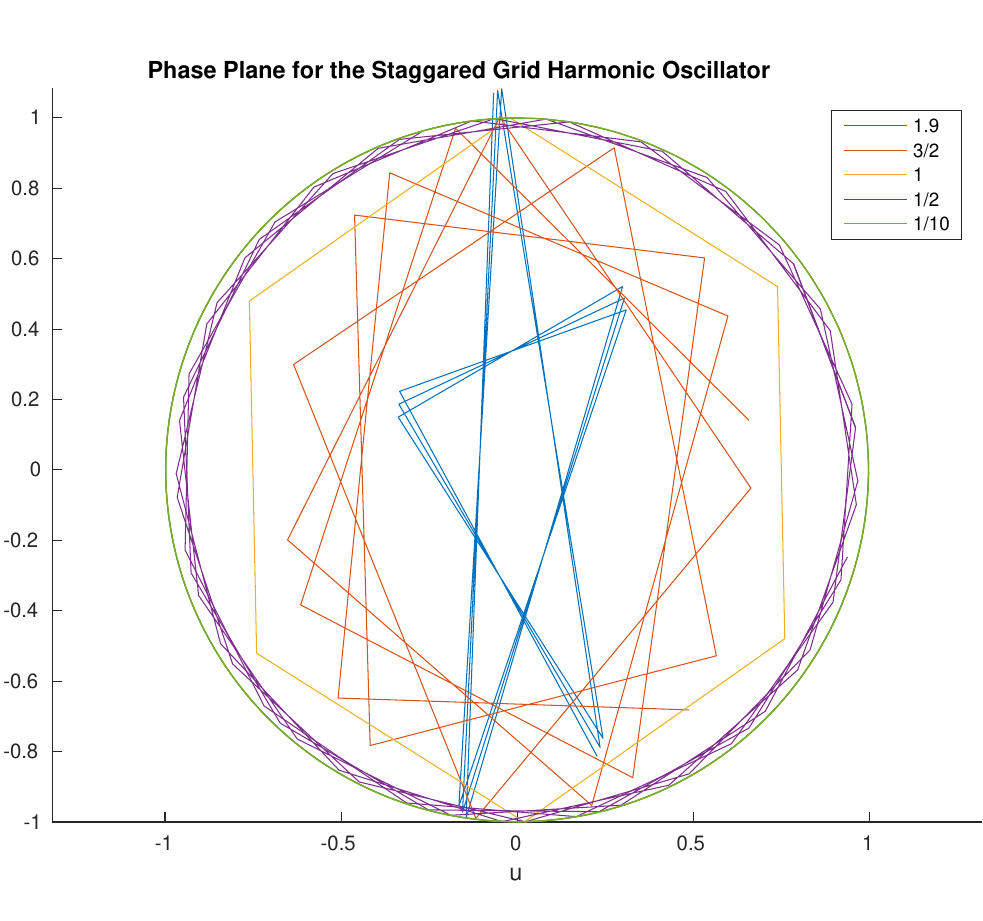}
\caption{Phase plane plots for the second order harmonic oscillator model
with $\omega = 1$ and $\Delta t = 3/2,\, 1 ,\, 1/2,\, 1/10$}.
\label{SecondOrderOsc}
\end{center}
\end{figure}

A natural proposal for a second-order accurate discrete conserved
quantity is 
\[
C^n =
(u^n)^2 + \left( \frac{u^{n+1}-u^{n-1}}{2 \, \omega \, \Delta t} \right)^2 \,.
\]
A little algebra shows that $C^n$ is not conserved. However this computation
shows that
\[
C^n = \left(1-\left(\frac{\omega \, \Delta t}{2}\right)^2 \right) (u^n)^2 +
\left( \frac{u^{n+1}-u^{n-1}}{2 \, \omega \, \Delta t} \right)^2 
\]
is conserved. Consequently for $0 < \Delta t /\omega < 2$ the discretization
is stable.  It is important that this constraint is less restrictive than
requiring an accurate solution.
Thus there seems to be no advantage to using
a discretization that is stable for all $\Delta t$?

The initial condition 
$u^1 = u(0) + \Delta t \, u'(0)$
is only first order accurate but a Taylor series expansion can
be used to make the order of accuracy higher:
\begin{align}
u^1 & = u(0)
+ \Delta t \,  u'(0)
+ \frac{1}{2} \, \Delta t^2 \, u''(0)
+ \frac{1}{6} \, \Delta t^3 \, u'''(0)
+ \cdots \nonumber \\
    & = u(0) 
+ \Delta t \,  u'(0)
- \frac{1}{2} \, \Delta t^2 \, \omega^2 \, u(0)
- \frac{1}{6} \, \Delta t^3 \, \omega^2 \, u'(0)
+  \cdots \label{HigherOrder} \\
\nonumber
\end{align}

The program {\tt Oscillator2ndOrder.m} produced Figure
\ref{SecondOrderOsc}, confirming that the algorithm is stable for
$0 < \Delta t /\omega < 2$ and that $C^n$ is constant to less than one
part in $10^{15}$.

\subsection{Staggering the Time Discretization}

A time staggered grid is used to discretize the first order system
\eqref{First Order System} which is given by a primal grid
$t^n = n \, \Delta t$ and a dual grid $t^{n+1/2} = (n+1/2) \, \Delta t$,
where $\Delta t > 0$ and $0 \leq n < \infty$ is an integer.
The staggered or leapfrog discretization of the harmonic oscillator is
then given by
\begin{equation}\label{Main Difference Equations}
\frac{u^{n+1}-u^n}{\Delta t} = \omega \, v^{n+1/2} , \quad  n \geq 0
\,,\quad
\frac{v^{n+1/2}-v^{n-1/2}}{\Delta t} = -\omega \, u^{n} , \quad n \geq 1
\,.
\end{equation}
The minus sign can be put in either equation, but it is important
to have an $\omega$ in both equations.  As before, the initial conditions
$u(0)$ and $u'(0)$ are given and then $u^0 = u(0)$ and
\begin{align}
v^{1/2} & = v\left(\frac{\Delta t}{2}\right) \nonumber \\
& = v(0) +
\frac{\Delta t }{2} \, v'(0) +
\frac{1}{2} \, \left(\frac{\Delta t }{2} \right)^2 \, v''(0) +
\frac{1}{6} \, \left(\frac{\Delta t }{2} \right)^3 \, v'''(0) + \cdots 
	\nonumber \\
& = v(0)+
\frac{\Delta t }{2} \, \omega \, u(0) +
\frac{1}{2} \, \left(\frac{\Delta t }{2} \right)^2 \, \omega^2 \, v(0) +
\frac{1}{6} \, \left(\frac{\Delta t }{2} \right)^3 \, \omega^3 \, u(0) +
\cdots \label{Higher Order IC} \\
\nonumber
\end{align}
The update algorithm starts with $u^0$ and $v^{1/2}$ and then for
$n \geq 0$
\begin{align*}
u^{n+1} = u^n + \Delta t \,  \omega \, v^{n+1/2} \,,\quad
v^{n+3/2} = v^{n+1/2} - \Delta t \, \omega \, u^{n+1} \,.
\end{align*}
Note that the second equation depends on the update in the first equation,
so the order of evaluation is critical.

This staggered grid discretization gives two standard single grid
discretization of the second order oscillator equation:
\[
\frac{u^{n+2} - 2 \, u^{n+1} +u^n }{\Delta t^2} + \omega^2 u^{n+1} = 0
\,;\quad
\frac{v^{n+3/2} - 2 \, v^{n+1/2} +v^{n-1/2} }{\Delta t^2} +
	\omega^2 v^{n+1/2} = 0 \,.
\]
So the solution of the fractional step methods is identical to the 
solution of the second order equations.

Again a simple proposed conserved quantity for
\eqref{Main Difference Equations} is
\begin{equation}\label{Simple Conserved}
C^n = \frac{1}{2}
\left( (u^n)^2 +\left( \frac{v^{n+1/2}+v^{n-1/2)}}{2}\right)^2 \right) \,.
\end{equation}
A little algebra gives
\[
C^{n+1} - C^n = \frac{\omega^2 \, \Delta t^2 }{4} \,
    \left((u^{n+1})^2 - (u^n)^2 \right) \,.
\]
So $C^n$ is {\em not} conserved. However, set
\[
\alpha = \frac{\omega \, \Delta t}{2} \,,
\]
and then the following two quantities are conserved:
\begin{equation}\label{Scaler Conserved n}
C^n = \frac{1}{2}
\left(
\left(1 - \alpha ^2 \right)
(u^n)^2 +\left( \frac{v^{n+1/2}+v^{n-1/2)}}{2}
\right)^2 \right) \,;
\end{equation}
\begin{equation}\label{Scaler Conserved n+1/2}
C^{n+1/2} = \frac{1}{2} \left(
\left(\frac{u^{n+1}+u^{n}}{2}\right)^2 + 
\left(1 - \alpha^2 \right) (v^{n+1/2})^2 \right) \,.
\end{equation}
The important properties for the staggered scheme are that it is explicit,
second order accurate and stable for $ \alpha = \omega \, \Delta t / 2 < 1$. 
By modifying the discretization, a similar result was obtained in
\cite{SchuhmannWeiland2001}, Equation 45, for the Yee time discretization
of Maxwell's equations.

The code {\tt StaggeredOscillator.m} confirms that the two energies
are constant to one part in $10^{15}$. The phase plane plots for the
staggered grid and the second order equation are identical. The code
also estimates that for $\omega = 1$, $\Delta t < 2$ is required for
stability, but for such a large $\Delta t$ the numerical solution
is very inaccurate, as made clear in Figure \ref{SecondOrderOsc}.
So the stability constraint on the time step is far less stringent
than the accuracy constraint.

\subsection{Summary}

If conserved quantities for the harmonic oscillator are allowed to
depend on $\Delta t$ then it is possible to derive conserved quantities
that converge quadratically to the energy of the continuum differential
equation.  The restriction on $\Delta t$ to keep the conserved quantity
positive is less stringent than the restriction for reasonably accurate
solutions.  All discretizations considered are confirmed to be second order
accurate in using {\tt StaggeredOscillator.m}. For generalizing these
results to more complex wave equation it is important that there are
no division by $\omega$ in the numerical algorithms.

\newpage \clearpage
\setcounter{equation}{0}
\section{Systems of Ordinary Differential Equations \label{ODEs}}

The next task is to consider a special class of systems of linear ordinary
differential equations that are wave equations. The discrete conservation
laws are easy to find by following the harmonic oscillator example.

\subsection{Continuous Time} \label{Continuous Time}

Let $X$ and $Y$ be linear spaces (finite or infinite dimensional).  It is
important that it is not assumed that $X$ and $Y$ have the same dimension.
If $f$ and $g$ are in $X$ then their inner product is $\langle f , g \rangle$
and the norm of $f$ is given by $\norm{f}^2 = \langle f , f \rangle$,
with the same notation for $Y$.  Let $A$ be a linear operator
mapping $X$ to $Y$ with adjoint $A^*$,  then
\[
X \overset{A}{\rightarrow} Y \,,\quad
Y \overset{A^*}{\rightarrow} X \,,
\]
and if $f \in X$ and $g \in Y$ then
\[
\langle A \, f , g \rangle = \langle f, A^* \, g \rangle \,.
\]

Next, if $f = f(t) \in X$ and $g = g(t) \in Y$ then a generalization of the
harmonic oscillator system is given by
\begin{equation} \label{Wave System}
f' = A \, g \,,\quad g' = - A^* \, f  \,,
\end{equation}
or in matrix form 
\[
\left[ \begin{matrix}
f'  \\
g'
\end{matrix} \right]
= 
\left[ \begin{matrix}
0 & A \\
- A^*& 0 
\end{matrix} \right]
\left[ \begin{matrix}
f \\
g
\end{matrix} \right] \,.
\]
Because the matrix
\[
\left[ \begin{matrix}
0 & A \\
- A^*& 0 
\end{matrix} \right]
\]
is skew adjoint, it must have purely imaginary spectra and the solutions of
this system must be made up of waves and constant solutions.  All solutions
are bounded in $t$.
Both $f$ and $g$ are solutions of second order linear wave equations:
\[
f'' + A^* \, A f = 0 \,;\quad g'' + A \, A^* g \ = 0.
\]
There are three natural initial conditions:
for the system specify $f(0)$ and $g(0)$;
for the second order equation in $f$ specify, $f(0)$ and $f'(0)$;
and for the second order equation in $g$ specify, $g(0)$ and $g'(0)$.

If the dimensions of $X$ and $Y$ are the same so that it make sense to
assume that $A$ is invertible then the system \ref{Wave System} will
have properties similar to the harmonic oscillator system
\eqref{First Order System} when $\omega > 0$. The most interesting case
is when the dimensions of the spaces are different which provides insight
in to the discretization of the scalar and vector wave equation and also
the Maxwell equations.  Also the case when $A$ is self adjoint, $A^* = A$,
provides insight into the discretization of Maxwell's equations.

If $f \in X$ and $g \in Y$ then
\begin{align*}
\langle  A^* \, A f , f \rangle = & 
   \langle  A f , A f \rangle = 
   \langle  f , A^* \, A f \rangle \, \\
\langle  A \, A^* g , g \rangle = &
   \langle  A^* g , A^* g \rangle = 
   \langle  g , A \, A^* g \rangle  \,.
\end{align*}
Consequently both $A\,A^*$ and $A^* \, A$ are self-adjoint positive operators,
but they may not be positive definite.  Also if $h \neq 0$ and $A \, h = 0$
then $ g(t) = t\, h$ is an unbounded solution of the second second order
equation while if $A^*h = 0$ then $f(t) = t \, h$ is an unbounded solution
of the first second order equation. For this $f(t)$ the system becomes
$ h = A \, g(t) \,,\quad g'(t) = 0  \,.  $ So $g(t) = k$ a
constant and then
$ \langle h , h \rangle = \langle h , A \, k \rangle = \langle A^* h ,  k \rangle =  0 $,
that is $h=0$ and then $f(t) = 0$ and $g(t) = k$ and $A \, k = 0$.
So the unbounded solution of the second order equation is not a solution
of the system, an advantage of using the system.
If $f$ and $g$ are vectors of the same length and $A$
is invertible then the system and second order equations are consistent.

There is also a problem with the initial conditions for the system
and the second order equations. If $A$ is an $n$ by $m$ matrix, then
$A^*$ is $m$ by $n$ matrix and consequently $A \, A^*$ is an $n$ by
$n$ matrix and $ A^* \, A$ is an $m$ by $m$ matrix.  So the first of
the second order equation needs $2\,n$ initial conditions, and the second
of the second-order equations needs $2\,m$ initial conditions. The
system needs $n+m$ initial conditions.
However, for example, if one knows $f(t)$ then $g(t)$ can be found using
simple integration and the initial condition for $g(t)$ and conversely for
knowing $g(t)$.
If $n=m$ then the number of initial conditions are the same for all three variants of the ordinary differential equations.
The $n \neq m$ is far more analogous to the situation for the scalar
and vector wave and Maxwell's equations than the $n=m$ case. 

\subsection{Continuous Time Conserved Quantities}

An important point here is that there is a conserved quantity that is
not analogous to energy but implies that the energy is conserved.
The fundamental conserved quantity is
\[
C(t) = \frac{1}{2} \left( \norm{ f(t) }^2 + \norm {g(t) }^2  \right) \,,
\]
which is analogous to \eqref{Conserved Quantity Simple}.
Because
\[
C'(t) = \langle f'(t), f(t) \rangle + \langle g'(t), g(t) \rangle  =
\langle A \, g(t), f(t) \rangle + \langle - A^* f(t), g(t) \rangle  = 0 \,,
\]
this quantity is conserved.

For the second order equations, an analog of the total energy that is the
sum of the kinetic plus the potential energy is given by
\begin{align}
\label{Continuum Total Energy}
E(t)
& = \frac{1}{2} \left( \norm{ f'(t) }^2 + \norm {A^* \, f(t) }^2  \right)
\nonumber \\
& = \frac{1}{2} \left( \norm{ A \, g(t) }^2 + \norm {g'(t) }^2  \right) 
\nonumber \\
& = \frac{1}{2} \left( \norm{ f'(t) }^2 + \norm {g'(t) }^2  \right) \,,
\end{align}
and is conserved because if $f, g$ are a solutions of the system then so are
$f'$ and $ g'$.  Note that the linearly growing solution has constant energy
but $C(t)$ is unbounded.  The $C(t$) type conserved quantities will used
from now on.

\subsection{Staggered Time Discretization\label{Staggered Time Discretization}}

A second order centered leapfrog discretization for the first order system is  
\[
\frac{f^{n+1} - f^n}{\Delta t} = A \, g^{n+1/2}  \,,\quad
\frac{g^{n+1/2} - g^{n-1/2}}{\Delta t} = -A^* \, f^{n} \,.
\]
Assuming that $f^0$ and $g^{1/2}$ are given then for $n \geq 0$ the
leapfrog time stepping scheme is 
\[
f^{n+1} = f^n + \dt \, A \, g^{n+1/2} \,,\quad
g^{n+3/2} = g^{n+1/2} - \dt \, A^* \, f^{n+1} \,.
\]
Again that the order of evaluation is important.

The initial conditions for the discretized system require $f^0 = f(0)$
and $g^{1/2} = g(\Delta t/2)$.  If $g(0)$ and $g'(0)$ are given then
\[
g^{1/2} \approx g^0 + \frac{\Delta t}{2} \, g'(0)  \,.
\]
If $f(0)$ and $g(0)$ are given then
\[
g^{1/2} \approx
g^0 + \frac{\Delta t}{2} \, g'(0) =
g^0 - \frac{\Delta t}{2} \, A^* \, f^0 \,.
\]
If $f(0)$ and $f'(0)$ are given, interchange $g$ and $f$ in the
discretization. The estimates can be more made more accurate as was 
done in section \ref{Higher Order IC}.

Both $f$ and $g$ satisfy a second order difference equation:
\begin{align} \label{Second Order}
\frac{f^{n+1} - 2 f^n + f^{n-1}}{\Delta t^2}
    & = - A \, A^* f^{n} \,;  \nonumber \\
\frac{g^{n+3/2} - 2 g^{n+1/2} + g^{n-1/2}}{\Delta t^2}
    & = - A^* \, A \, g^{n+1/2} \,. \\
\nonumber
\end{align}
Additionally a second order average is needed for computing
conserved quantities:
\begin{align} \label{Second Order Average}
\frac{f^{n+1} + 2 \, f^{n} + f^{n-1}}{4}
& = f^{n}  - \frac{\Delta t^2}{4} A \, A^* f^{n} \,; \nonumber \\
\frac{g^{n+3/2} + 2 \, g^{n 1/2} + g^{n-1/2}}{4} 
& = 
g^{n+1/2}  - \frac{\Delta t^2}{4} A^* \, A \, g^{n+1/2} \,. \\
\nonumber
\end{align}

When comparing this discretization to the simple oscillator
discretization it is important that $\omega > 0$, while here
the operators $A$ and $A^*$ may not be invertible which is typically
the case when studying spatially dependent partial differential wave equations.

\subsection{Discrete Time Conserved Quantities}

To show that $A$ not being invertible is not serious problem a detailed
derivation of the conservation laws that are analogs of
\eqref{Scaler Conserved n+1/2} and \eqref{Scaler Conserved n} are given.  Let
\begin{align*}
C_1^{n+1/2} & = \norm{\frac{f^{n+1} + f^{n}}{2}}^2 \,, \\
C_2^{n+1/2} & = \norm{ g^{n+1/2}}^2 \,, \\
C_3^{n+1/2} & = \Delta t^2 \norm{ A \, g^{n+1/2} }^2  \,.
\end{align*}
As before compute:
\begin{align*}
C_1^{n+1/2}- C_1^{n-1/2}
& = \langle \frac{f^{n+1} + 2 \, f^{n}  + f^{n-1}}{4} , f^{n+1} - f^{n-1} \rangle \\
& = \langle f^{n}  - \frac{\Delta t^2}{4} A \, A^* f^{n} , f^{n+1} - f^{n-1} \rangle \,;
\end{align*}
\begin{align*}
C_2^{n+1/2}- C_1^{n-1/2}
& = \langle g^{n+1/2} + g^{n-1/2} , g^{n+1/2} - g^{n-1/2} \rangle \\
& = \langle g^{n+1/2} + g^{n-1/2} , - \Delta t \, A^* f^n \rangle \\
& = - \Delta t \, \langle A \, g^{n+1/2} + A \, g^{n-1/2} , f^n \rangle \\
& = -\Delta t \, \langle \frac{f^{n+1}-f^{n-1}}{\Delta t} , f^n \rangle \\
& = \langle -f^n , f^{n+1}-f^{n-1} \rangle \,;
\end{align*}
\begin{align*}
C_3^{n+1/2}- C_1^{n-1/2} & = \Delta t^2
\langle A \, g^{n+1/2} - A \, g^{n-1/2} \,,\, A \, g^{n+1/2} + A \, g^{n-1/2} \rangle \\
& = \Delta t^ 2 \langle A \left( g^{n+1/2} - g^{n+1/2} \right) \,,\, \frac{f^{n+1}-f^{n-1}}{\Delta t} \rangle \\
& = \Delta t^ 2 \langle - \Delta t \, A \, A^* f^{n} \,,\,
	\frac{f^{n+1}-f^{n-1}}{\Delta t} \rangle \\
& = - \Delta t^ 2 \langle A \, A^* f^{n} \,,\, f^{n+1}-f^{n-1} \rangle \,.
\end{align*}
Consequently $C = C_1 + C_2 - C_3/4$ is a conserved quantity, that is
\begin{equation}
C^{n+1/2} = \norm{\frac{f^{n+1} + f^{n}}{2}}^2 +
\left( 1 - \frac{\Delta t^2}{4} \norm{A}^2 \right)
\norm{ g^{n+1/2}}^2  \,.
\label{Conserved A half}
\end{equation}
is positive and constant for $\Delta t$ sufficiently small.

Next let
\begin{align*}
C_1^{n} & = \norm{\frac{g^{n+1/2} + g^{n-1/2}}{2}}^2 \,, \\
C_2^{n} & = \norm{f^{n}}^2 \,, \\
C_3^{n} & = \Delta t^2 \norm{ A^* \, f^{n} }^2 \,. 
\end{align*}

\begin{align*}
C_1^{n+1}- C_1^{n}
& = \langle
\frac{g^{n+3/2} + 2 \, g^{n 1/2} + g^{n-1/2}}{4} \,,\,
   g^{n+3/2} -g^{n-1/2} \rangle \\
& = \langle
g^{n+1/2}  - \frac{\Delta t^2}{4} A^* \, A \, g^{n+1/2} \,,\,
   g^{n+3/2} -g^{n-1/2} \rangle \,.
\end{align*}
\begin{align*}
C_2^{n+1}- C_2^{n} 
& = \langle f^{n+1}-f^{n} \,,\, f^{n+1}+f^{n} \rangle \\
& = \langle \Delta t \,  A \, g^{n+1/2} \,,\, f^{n+1}+f^{n} \rangle \\
& = \Delta t \langle g^{n+1/2} \,,\, A^* f^{n+1}+A^* f^{n} \rangle \\
& = \Delta t \langle g^{n+1/2} \,,\,
	-\frac{g^{n+3/2}- g^{n-1/2}}{\Delta t} \rangle \\
& = \langle - g^{n+1/2} \,,\,
	g^{n+3/2}- g^{n-1/2} \rangle  \,.
\end{align*}
\begin{align*}
C_3^{n+1}- C_3^{n}
& = \Delta t^2 
\langle
   A^* f^{n+1}- A^* f^n  \,,\,  A^* f^{n+1}+ A^* f^n
\rangle \\
& = \Delta t^2 \langle
   A^* f^{n+1}- A^* f^n  \,,\, \frac{g^{n+3/2}-g^{n-1/2}}{\Delta t}
\rangle \\
& = \Delta t^2 \langle
   \Delta t \, A^* A \, g^{n+1/2} \,,\, \frac{g^{n+3/2}-g^{n-1/2}}{\Delta t}
\rangle \\
& = \Delta t^2 \langle
   A^* A \, g^{n+1/2} \,,\, g^{n+3/2}-g^{n-1/2} \rangle \,.
\end{align*}
Consequently $C^n = C_1^n + C_2^n - C_3^n/4$ is a conserved quantity,
that is
\begin{equation}
\norm{C^n} = \left(1 - \frac{\Delta t^2}{4} \norm{A^*}^2\right) \norm{f^n}^2
	+ \norm{\frac{g^{n+1/2} + g^{n-1/2}}{2}}^2 \,,
\label{Conserved A* full}
\end{equation}
is positive and constant for $\Delta t$ sufficiently small.

The program {\tt SystemsODEs.m} tests these conservation laws for $A$ a
$2 \times 3$ random matrix showing that the energies are constant with
an error less that than one part in $10^{14}$.

\newpage \clearpage
\setcounter{equation}{0}
\section{Discretizing the 1D  Wave Equation\label{1D Wave}}

\begin{figure}
\begin{center}
   \includegraphics[width=5.00in,trim = 0 320 0 50, clip ]{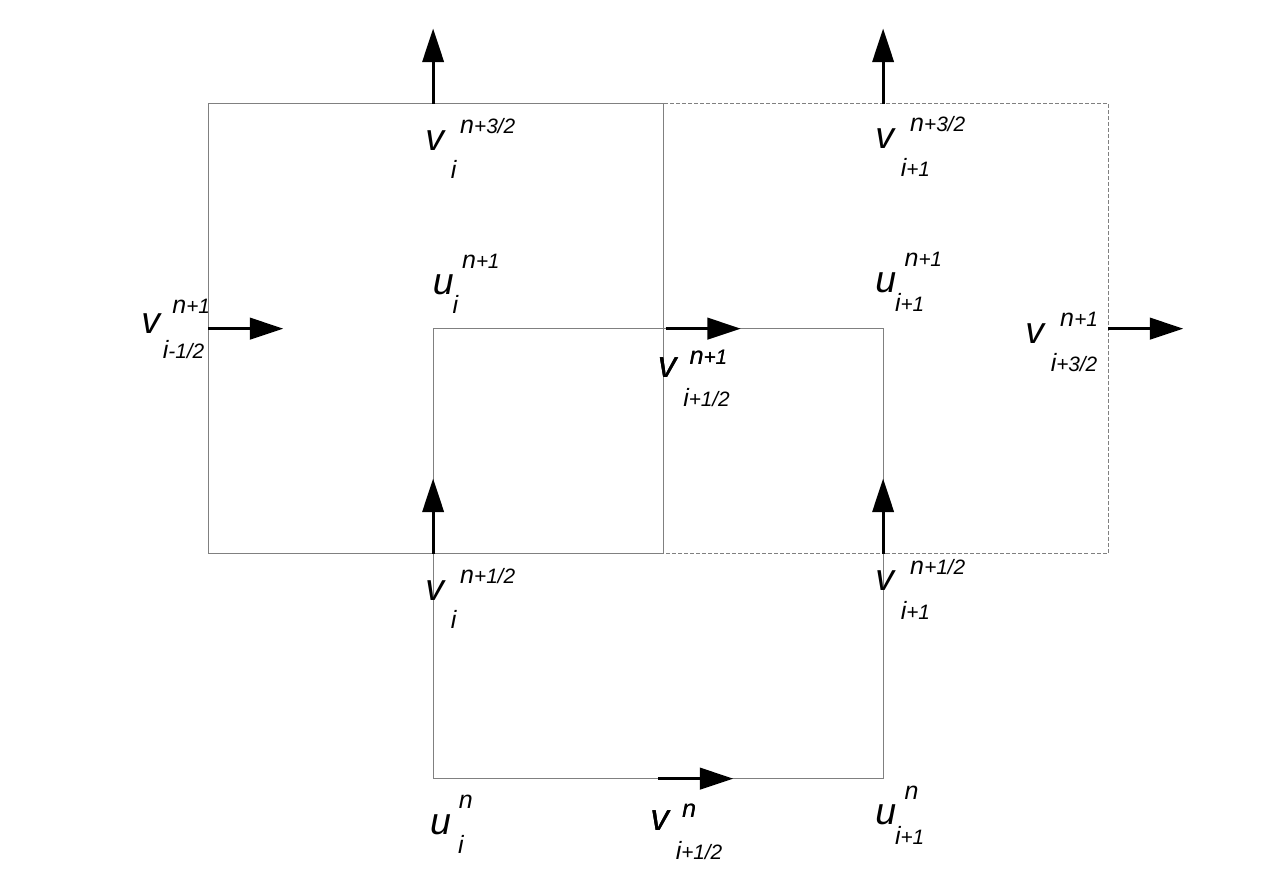}
\caption{Space-Time Staggered Grid}
\label{TimeSpaceGrid}
\end{center}
\end{figure}

The 1D scalar wave equation will be discretized by writing the equation
as a system of two first order equations and then using a staggered time
and spatial discretization. The time discretization is the same as
leapfrog discretization used before while the spatial
discretization is the mimetic discretization specialized to one dimension.
Most importantly, a conserved quantity $C$ is introduced that is not
the classical energy, but the conservation of $C$ implies the conservation
of the energy $E$. This will play an important role in 3D discretizations.

Let $u = u(t,x)$ be a smooth real valued function of the real variables
$x$ and $t$ such that $u(t,\pm \infty) = 0$.
Then let
$u_t = \partial u / \partial t$,
$u_x = \partial u / \partial x$,
$u_{tt} = \partial^2 u / \partial t^2$,
and
$u_{xx} = \partial^2 u / \partial x^2$.
The 1D wave equation is 
\[
u_{tt} = c^2 \, u_{xx} \,,
\]
where $c>0$.
The initial conditions for this equation are $u(0,x)$ and $u_t(0,x)$.

This equation can also be written as a system
\begin{equation} \label{Wave-1D-System}
u_t = c \, v_x \,,\quad v_t = c \, u_x \,,
\end{equation}
where again $v$ is smooth and $v(t,\pm \infty) = 0$. The
initial conditions are $u(0,x)$ and $v(0,x)$.
As before $v$ also satisfies a second order wave equation
\[
v_{tt} = c^2 \, v_{xx} \,.
\]

The vector spaces $X$ and $Y$ from Section \ref{ODEs} are replaced by
$L^2$, the functions defined on the real line and are square integrable.
The inner product of two functions $f = f(x)$ and $g = g(x)$ is
\[
\langle f,g \rangle = \int_{-\infty}^{\infty} f(x) \, g(x) \, dx \,.
\]
If $f(\pm \infty) = 0$ and $g(\pm \infty) = 0$ then integration by parts
gives $\langle f', g \rangle = \langle f, -g' \rangle$, so if
$A = \partial/\partial x$ then $-A^* = A$.  So the wave equation has
the has the same structure as the equations in the previous sections.

The usual energy $E = E(t)$ for the wave equation is the kinetic plus the
potential energies,
\[
E =
\frac{1}{2} \int_{-\infty}^{\infty}
\left( u_{t}^2 + c^2 \, u_{x}^2 \right) d x \,.
\]
Use integration by parts to see that
\begin{align*}
E_t
= & \int_{-\infty}^{\infty}
   \left( u_t \, u_{tt} + c^2 \, u_x  \, u_{tx} \right) \, d x \,. \\
= &\int_{-\infty}^{\infty}
\left( u_t \, u_{tt} - c^2 \, u_{xx} u_t \right) \, d x  = 0 \,,
\end{align*}
that is, the energy $E(t)$ is conserved.
As indicated in the previous sections a preferred conserved quantity
is $C = C(t)$ where
\begin{equation} \label{Conserved 1D Wave}
C =
\frac{1}{2} \int_{-\infty}^{\infty} \left( u^2 + v^2 \right) d x \,,
\end{equation}
because 
\begin{align*}
C_t & =  \int_{-\infty}^{\infty}
	\left( u \, u_t+ v \, v_t \right) d x \,, \\
   & =  \int_{-\infty}^{\infty}
	\left( u \, c \, v_x + v \, c \, u_x \right) d x \,, \\
   & =  c \, \int_{-\infty}^{\infty}
	\left( u \, v_x + v \, u_x \right) d x \,, \\
   & =  c \, \int_{-\infty}^{\infty}
	\left( u \, v \right)_x d x \,, \\
   & = 0 \,.	
\end{align*}
Again note that if $u,v$ are solutions of the system \eqref{Wave-1D-System}
then so are $u_t$, $v_t$ and then \eqref{Conserved 1D Wave} implies that
a conserved quantity is given by the energy
\[
\frac{1}{2} \int_{-\infty}^{\infty} \left( u_t^2 + v_t^2 \right) d x = 
\frac{1}{2} \int_{-\infty}^{\infty} \left( u_t^2 + c^2 u_x^2 \right) d x 
= E \,.
\]
So if $C$ is conserved then so is $E$.

\subsection{A Staggered Discretization of the Wave Equation}

Let $\dt > 0$ and $\dx > 0$ be given and then the primal and
dual grids points are given by
\begin{align*}
\left(t^n,x_i\right) = & \left(n\,\dt, i \, \dx \right) \,, \\
\left(t^{n+1/2},x_{i+1/2}\right) = &
	\left((n+1/2)\,\dt, (i+1/2) \, \dx \right) \,,
\end{align*}
where $ -\infty < n < \infty$ and $ -\infty < i < \infty$.
The discretization of $u(x,t)$ on the primal grid and $v(x,t)$ on the dual
grid are $u^n_i$ and $v^{n+1/2}_{i+1/2}$. Then the system \eqref{Wave-1D-System}
is discretized as
\begin{align}\label{Discretized System}
\frac{u^{n+1}_{i} - u^{n}_{i}}{\dt}
	& = c \, \frac{v^{n+1/2}_{i+1/2} - v^{n+1/2}_{i-1/2}}{\dx} \,,
	\nonumber \\
\frac{v^{n+1/2}_{i+1/2} - v^{n-1/2}_{i+1/2}}{\dt} 
	& = c \, \frac{u^n_{i+1}-u^n_i}{\dx} \,.
\end{align}
Assume that $u^0$ and $v^{\half}$ are given then the leapfrog
time stepping scheme for $n>0$ is
\[
u^{n+1}_{i} =
u^{n}_{i} +
c \, \frac{\dt}{\dx} \left(v^{n+1/2}_{i+1/2} - v^{n+1/2}_{i-1/2}\right) \,,\quad
v^{n+3/2}_{i+1/2} =
v^{n+1/2} + c \, \frac{\dt}{\dx} \, \left(u^{n+1}_{i+1}-u^{n+1}_i\right) \,.
\]

This implies that both $u$ and $v$ satisfy a discretization of the
second order wave equation:
\begin{align*}
\frac{u^{n+1}_{i} - 2\,u^{n}_{i} +u^{n-1}_{i}}{\dt^2} 
& = \frac{1}{\dt} \left( \frac{u^{n+1}_{i} - \,u^{n}_{i}}{\dt} -  
   \frac{u^{n}_{i} - u^{n-1}_{i}}{\dt} \right) \\
& = \frac{c}{\dt} \left(
	\frac{ v^{n+1/2}_{i+1/2} - v^{n+1/2}_{i-1/2} }{\dx} -
	\frac{ v^{n-1/2}_{i+1/2} - v^{n-1/2}_{i-1/2} }{\dx} \right) \\
& = \frac{c}{\dx} \left(
	\frac{ v^{n+1/2}_{i+1/2} - v^{n-1/2}_{i+1/2} }{\dt} -
	\frac{ v^{n+1/2}_{i-1/2} - v^{n-1/2}_{i-1/2} }{\dt} \right) \\
& = \frac{c^2}{\dx} \left( \frac{u^n_{i+1}- u^n_{i}}{\dx}
	- \frac{u^n_{i}- u^n_{i-1}}{\dx} \right) \\
& =  c^2 \,\frac{u^n_{i+1}- 2 \, u^n_{i} +u^n_{i-1} }{\dx^2} \,.
\end{align*}
A similar calculation shows that
\[
\frac{v^{n+1/2}_{i-1/2}-2\,v^{n-1/2}_{i-1/2}+v^{n-3/2}_{i-1/2}}{\dt^2}
=
c^2 \, \frac{v^{n-1/2}_{i+1/2}-2\,v^{n-1/2}_{i-1/2}+v^{n-1/2}_{i-3/2}}{\dx^2} \,.
\]

The inner product of two grid functions
$a = (\cdots, a_{-1}, a_0 , a_1 , \cdots)$
and
$b = (\cdots, b_{-1}, b_0 , b_1 , \cdots)$:
is given by
\[
\langle a,b \rangle = \sum_{i=-\infty}^\infty a_i \, b_i \, , \quad
\norm{a}^2 = \langle a,a \rangle \,.
\]
Similarly, if 
$c = (\cdots, c_{-1/2}, c_{1/2} , c_{3/2} , \cdots)$
and
$d = (\cdots, d_{-1/2}, d_{1/2} , d_{3/2} , \cdots)$
then
\[
\langle c,d \rangle = \sum_{i=-\infty}^\infty c_{i+1/2} \, d_{i+1/2} \,, \quad
\norm{c}^2 = \langle c,c \rangle \,.
\]
The discrete analogs of the integration by parts formula will be needed
so let 
\[
\delta(a)_{i+1/2} = a_{i+1}-a_i \,,\quad \delta(c)_{i} = c_{i+1/2}-c_{i-1/2} \,.
\]
Then the summation by parts formula is given by
\begin{align*}
\langle \delta(a) , c \rangle
& = \sum_{i=-\infty}^\infty \left( a_{i+1} - a_i \right)  c_{i+1/2}  \\
& =  \sum_{i=-\infty}^\infty a_{i+1} c_{i+1/2} 
   - \sum_{i=-\infty}^\infty a_i c_{i+1/2}  \\
& =  \sum_{i=-\infty}^\infty a_{i} c_{i-1/2} 
   - \sum_{i=-\infty}^\infty a_i c_{i+1/2}  \\
& = - \sum_{i=-\infty}^\infty a_i  \left( c_{i+1/2} - c_{i-1/2} \right)  \\
& = - \langle a , \delta (c) \rangle \,.
\end{align*}

The difference equations \eqref{Discretized System} can now be written
\begin{align}
\label{Vector Discretized System}
\frac{u^{n+1}_{i} - u^{n}_{i}}{\dt}
	& = c \, \frac{\delta(v^{n+1/2})_{i}}{\dx} \,,
	\nonumber \\
\frac{v^{n+1/2}_{i+1/2} - v^{n-1/2}_{i+1/2}}{\dt} 
	& = c \, \frac{\delta(u^n)_{i+1/2}}{\dx} \,.
\end{align}
To find a conserved quantity define:
\begin{align*}
C1(n) & = \norm{u^n}^2 \,; \\
C2(n) & = \norm{\frac{v^{n+1/2}+v^{n-1/2}}{2}}^2 \,; \\
C3(n) & = \norm{\delta(u^n)}^2 \,.
\end{align*}
Now
\begin{align*}
C1(n+1) - C1(n) & = \langle u^{n+1}-u^n , u^{n+1}+u^n \rangle \\
& = 
c \, \frac{\dt}{\dx}
	\langle \delta(v^{n+1/2}) ,
	u^{n+1}+u^n \rangle \quad
	\text{see} \eqref{Vector Discretized System} \\
& =
	c \, \frac{\dt}{\dx} \langle v^{n+1/2} , 
	- \delta(u^{n+1}+u^n) \rangle  \\
& =
	- c \, \frac{\dt}{\dx} \langle v^{n+1/2} , 
	 \delta(u^{n+1})+\delta(u^n) \rangle  \\
& =
	-\langle v^{n+1/2} , v^{n+3/2} - v^{n-1/2}  \rangle  \quad
	\text{see} \eqref{Vector Discretized System} \,. 
\end{align*}
\begin{align*}
C2(n+1) -C2(n)
& = \frac{1}{4} \langle
v^{n+3/2}+v^{n+1/2}+v^{n+1/2}+v^{n-1/2} ,
v^{n+3/2}+v^{n+1/2}-v^{n+1/2}-v^{n-1/2} \rangle \\
& = \frac{1}{4}
   \langle v^{n+3/2}+2\,v^{n+1/2}+v^{n-1/2} , v^{n+3/2}-v^{n-1/2} \rangle  \\
& = 
\langle v^{n+1/2} , v^{n+3/2}-v^{n-1/2} \rangle +
\frac{1}{4}
\langle v^{n+3/2}-2\,v^{n+1/2}+v^{n-1/2} , v^{n+3/2}-v^{n-1/2} \rangle \\
& = \langle v^{n+1/2} , v^{n+3/2}-v^{n-1/2} \rangle +
c \, \frac{\dt}{4 \, \dx}
\langle \delta (u^{n+1}) - \delta (u^n) , v^{n+3/2}-v^{n-1/2} \rangle  \\
& = \langle v^{n+1/2} , v^{n+3/2}-v^{n-1/2} \rangle +
\left( c \, \frac{\dt}{2 \, \dx}\right)^2
\langle \delta ( \delta (v^{n+1/2})) , v^{n+3/2}-v^{n-1/2} \rangle  \\
\end{align*}
\begin{align*}
C3(n+1)- C3(n) 
& = \langle
\delta ( u^{n+1}) - \delta ( u^n) , \delta ( u^{n+1} ) + \delta ( u^n )
\rangle \\
& = \langle
\delta ( u^{n+1} -  u^n) , \delta ( u^{n+1} )  +  \delta( u^n )
\rangle \\
& = c \, \frac{\dx}{\dt} \langle
\delta ( u^{n+1} ) - \delta ( u^n ) , v^{n+3/2} - v^{n-1/2} \rangle \\
& =
\langle \delta ( \delta ( v^{n+1/2} ) ) , v^{n+3/2} - v^{n-1/2} \rangle \\
\end{align*}
Consequently the quantity
\[
C(n) = 
\norm{u^n}^2 
-\left(\frac{c \, \dt}{2 \, \dx}\right)^2
\norm{\delta  u^n }^2 +
\norm{\frac{v^{n+1/2}+v^{n-1/2}}{2}}^2  \,,
\]
is conserved.
A similar argument shows that 
\[
C(n+1/2) =
\norm{ v^{n+1/2}}^2
-\left(\frac{c \, \dt}{2 \, \dx}\right)^2 \norm{ \delta \, v^{n+1/2} }^2
+ \norm{\frac{u^{n+1} + u^{n}}{2}}^2 \,.
\]
is a conserved quantity. The program {\tt OneDWave.m} shows that the
conserved quantities are constant to within an error of less than 10e-14.

The first conserved quantity will be positive provided that
\[
\frac{c \, \dt}{2 \, \dx}
\frac{\norm{\delta  u^n}}{\norm{u^n}} < 1 \,.
\]
But 
\[
\frac{\norm{\delta  u^n}}{\norm{u^n}} \leq \norm{\delta} \,,
\]
so the conserved quantity will be positive if
\[
c \, \frac{\dt}{ \dx} < \frac{2}{\norm{\delta}} \,.
\]
Because $\norm{\delta} = 2$ (see {\tt Normdelta.m}) this is the
Courant-Friedrichs-Lewy (CFL) condition for stability.

The program {\tt Wave1D.m} shows that conservation errors are less than
one part in $10^{15}$ and that the convergence rate is 2.

\newpage \clearpage
\setcounter{equation}{0}

\section{3D Wave Equations\label{3D Wave Equations}}

\begin{table}[ht]
\begin{center}
\begin{tabular}{|l|c|l|}
\hline
quantity & units & name \\
\hline
$ \vec{x}$  & $d$     & spatial position   \\
$\vec{u}$   & $d$     & displacement       \\
$\rho > 0$  & $1/d^3$ & density            \\
$\sigma$    & $1/d^2$ & stress             \\
$e$         & $1$     & strain             \\
$C$         & $1/d^2$ & elastic properties \\
$\lambda>0$ & $1/d$   & Lam\'e parameter   \\
$\mu > 0$   & $1/d$   & Lam\'e parameter   \\
$K$         & $1/d$   & bulk modulus       \\
$Pa$        & $1/d$   & Pascal             \\
\hline
\end{tabular}
\caption{Quantities and their spatial units.}
\label{Elastic Units}
\end{center}
\end{table}

The main interest is in three dimensional wave equations of which there 
are several variants. Here the material properties are assumed to be
constant. What is important is the spatial dimension $d$ of the variables
and operators, see table \ref{Elastic Units}. The differential operators
divergence $\divg$, curl $\curl$ and gradient $\grad$ all have spatial
dimension $1/d$.  Many wave equations are derived from Newton's laws and
thus have the form
\begin{equation}
\rho \, \frac{d^2 W}{d t^2} = \mathscr{A} \, W \,,
\label{General Wave Equation}
\end{equation}
where $W = W(t,x,y,z)$ is a scalar or vector function, $\rho$ is the density
of the material that the wave is traveling in and $\mathscr{A}$ is $\pm$ a
constant times a second order differential operator. Because such equations
are linear in $W$ the dimensions of $W$ are not important. So the spatial
dimension of $\mathscr{A}$ must be the same as $\rho$. Consequently the
dimensionless form of this wave equation is
\[
\frac{d^2 W}{d t^2} = \frac{1}{ \rho } \mathscr{A} \, W \,.
\]
The differential operators in the wave equations will be compositions of
two of the operators divergence $\divg$, curl $\curl$ and gradient $\grad$
and thus will have spatial dimension $1/d^2$.

A critical point about wave equations is that the operator $\mathscr{A}$
must be negative definite that is all eigen values of this operator are
real and strictly less than zero. This will guarantee that the solutions
are oscillatory. For constant material properties functions that are linear
in time will be solutions of the second order equation. These solutions
do not go to zero at large distance from the origin and so are theoretically
not allowed but can still cause problems in simulations. 

Inner products and norms of scalar functions are needed to describe
the conservation laws for wave equations. So assume that the scalar and
vector functions are smooth and converge rapidly to zero for large distances
from the origin. If $f$ and $g$ are such scalar function then their inner
product and norm are
\[
\langle f, g \rangle = \int_{\real^3} f(x,y,z) \, g(x,y,z)
\, dx \,  dy \, dz \,,\quad \norm {f}^2 = \langle f , f \rangle \,,
\]
while if $\vec{u}$ and $\vec{v}$ are smooth vector functions then
their inner product and norm are given by
\[
\langle \vec{u}, \vec{v} \rangle =
\int_{\real^3} \vec{u}(x,y,z) \bdot \vec{v}(x,y,z) 
\, dx \,  dy \, dz \,,\quad
\norm {\vec{u}}^2 = \langle \vec{u} , \vec{u} \rangle \,.
\]

\subsection{The Scalar Wave Equation}
For example the scalar wave equation is given by 
\begin{equation}\label{3D Scalar Wave Equation}
\frac{d^2 f}{d t^2} =  c^2 \, \divg \, \grad \,  f \,,
\end{equation}
where $f = f(x,y,z,t)$ is a dimensionless scalar and $c$ is the constant sound
speed with spatial dimension $d$. Consequently this scalar wave equation is
dimensionless. Note that the sound speed is given by 
\[
c^2 = \frac{K}{\rho} \,,
\]
where $K$ is the bulk modulus of the material the sound is traveling in.

To apply the stagged time discretization to second order wave equations
they must be converted to a first order system. For the scalar wave
equation introduce $\vec{v} = \grad f$ to get
\begin{equation}
\frac{d f}{d t} = c \, \divg \vec{v} \,,\quad 
\frac{d \vec{v}}{d t} = c \, \grad f  \,.
\end{equation}
Putting a $c$ is each of the first oder equations makes them dimensionless.
Note that $\vec{v}$ satisfies a simple vector wave equation
\[
\frac{d^2 \vec{v}}{d t^2} =  c^2 \, \grad \, \divg \,  \vec{v} \,.
\]
Consequently solving the system will produce a solution to both the simple
scalar and vector wave equations.

For the first order system the quantity 
\[
C = \frac{\norm{f}^2 + \norm{\vec{v}}^2}{2}
\]
is conserved because
\begin{align}
\frac{d C}{d t} & = \langle f, \frac{df}{d t} \rangle
                +   \langle \vec{v}, \frac{d\vec{v}}{d t} \rangle \nonumber \\
                & = \langle f, c \, \divg \vec{v} \rangle
                +   \langle \vec{v}, c \, \grad f \rangle \nonumber \\
                & = c \, \langle f, \divg \vec{v} \rangle
                -   c \, \langle \divg \vec{v}, f \rangle \nonumber \\
		& = 0 
\end{align}

Note that if $f$ and $\vec{v}$ are solutions of the system then so are $df/dt$
and $d \vec{v}/dt$. So the conservation $C$ implies the conservation $E$. In
fact and time or spatial derivative of the solution of the wave equation is
also a solution so there are infinitely many conserved quantities.
Consequently
\[
E = 
\norm{\frac{d f}{d t}}^2 +
\norm{\frac{d \vec{v}}{d t}}^2 =
\norm{\frac{d f}{d t}}^2 + c^2 \, \norm{\grad f}^2 \,.
\]
is conserved and is the energy for the scalar wave equation. 

\subsection{The Elastic Wave Equation}
In the case that $W$ is a vector $\vec{v}$ there are only two second order
differential operators that can be made from the gradient $\grad$, the
divergence $\divg$ and the curl $\curl$ which are $\curl \curl \vec{v}$ and
$\grad \divg \vec{v}$.  So it is no surprise that the elastic wave equation
\cite{Igel06} is made up of these operators:
\begin{equation}
\label{Elastic Wave Equation}
\rho \, \frac{d^2 \vec v}{d t^2} =
(\lambda + 2 \, \mu) \grad \divg \vec{v}
- \mu \curl \curl \vec{v} \,,
\end{equation}
where $\mu$ and $\lambda$ are constant scalars with spatial dimension $1/d$.
Note that $\mu = 0$ produces
\begin{equation}
\label{Simple Vector Wave Equation}
\frac{d^2 \vec v}{d t^2} =
\frac{\lambda}{\rho} \grad \divg \vec{v} \,,
\end{equation}
which is the dimensionless simple vector wave equation.

If 
\[
a = \sqrt{\frac{(\lambda + 2 \, \mu)}{\rho}} \,,\quad
b = \sqrt{\frac{\mu}{\rho}}
\]
then, like for $c$ in the scalar wave equation, $a$ and $b$ have units $d$,
so the elastic wave equation \ref{Elastic Wave Equation} can be written in
dimensionless form as:
\begin{equation} \label{dimensionless elastic equation}
\frac{d^2 \vec v}{d t^2} =
a^2  \grad \divg \vec{v} - b^2 \curl \curl \vec{v} \,.
\end{equation}
This equation can be converted to a system of three dimensionless first order
equations:
\begin{align*}
\frac{d g}{d t} & = a \, \divg \vec{v} \\
\frac{d \vec{u}}{d t} & = b \, \curl \vec{v} \\
\frac{d \vec{v}}{d t} & = a \, \grad g - b \,\curl \vec{u} \,.
\end{align*}
The quantity
\[
C = \frac{\norm{\vec{v}}^2 + \norm{\vec{u}}^2 + \norm{g}^2}{2} \\,
\]
is conserved because
\begin{align}
\frac{d C}{d t}
& =\langle \frac{d \vec{v}}{d t}, \vec{v} \rangle 
+  \langle \frac{d \vec{u}}{d t}, \vec{u} \rangle 
+  \langle \frac{d \vec{g}}{d t}, \vec{g} \rangle \nonumber \\
& =\langle a \, \grad g - b \,\curl \vec{u}, \vec{v} \rangle 
+  \langle b \, \curl \vec{v} , \vec{u} \rangle 
+  \langle a \, \divg \vec{v}, \vec{g} \rangle \nonumber \\
& = a \, \langle \grad g , \vec{v} \rangle 
  - b \, \langle \curl \vec{u}, \vec{v} \rangle 
  + b \, \langle  \curl \vec{v} , \vec{u} \rangle 
  + a \, \langle  \divg \vec{v}, \vec{g} \rangle \nonumber \\
& = - a \, \langle g , \divg \vec{v} \rangle 
  - b \, \langle  \vec{u}, \curl \vec{v} \rangle 
  + b \, \langle  \curl \vec{v} , \vec{u} \rangle 
  + a \, \langle  \divg \vec{v}, \vec{g} \rangle \nonumber \\
& = 0 \nonumber
\end{align}

\subsection{Maxwell Equations\label{Section Maxwell Equations}}

\begin{table}[ht]
\begin{center}
\begin{tabular}{|l|c|l|l|}
\hline
quantity & units & name \\
\hline
$\vec{E}$ & $1/d$   & electric field \\
$\epsilon$& $1/d$   & permeability tensor \\
$\vec{D}$ & $1/d^2$ & electric displacement \\
\hline
$\vec{H}$ & $1/d$   & magnetic field \\
$\mu$     & $1/d$   & permittivity tensor \\
$\vec{B}$ & $1/d^2$ & magnetic flux \\
\hline
$\curl$   & $1/d$   & curl operator  \\
$\curl$   & $1/d$   & curl operator \\
\hline
$\vec{J}$ & $1/d^2$     & current \\
\hline
\end{tabular}
\caption{Quantities and their units in the Maxwell equations \cite{Crain17}.}
\label{Maxwell Units}
\end{center}
\end{table}

The Maxwell Equations
\begin{equation}\label{Maxwell Equations}
\frac{d \vec{B}}{d t} + \curl \vec{E} = 0 \,,\quad
\frac{d \vec{D}}{d t} - \curl \vec{H} = {\vec J} \,.
\end{equation}
\[
\vec{B} = \mu \, \vec{H} \,,\quad \vec{D} = \epsilon \, \vec{E} \,.
\]
provide an example that was studied by Yee \cite{Yee1966} with essentially the
same ideas that are used in this paper.  Here $\vec{B}$,
$\vec{E}$, $\vec{D}$ and $\vec{H}$ are vector functions of $(x,y,z,t)$
while $\mu$ and $\epsilon$ are symmetric positive definite matrices that
depend only on the spatial variables.  The meaning of variables and their
distance units are given in Table \ref{Maxwell Units}. The Maxwell equations
are a bit different in that they start as a first order system. 

Eliminate $\vec{B}$ and $\vec{D}$ from the equation to get
\begin{equation} \label{Maxwell System}
\frac{d \vec{E}}{d t} =   \epsilon^{-1} \curl \vec{H} \,,\quad
\frac{d \vec{H}}{d t} = - \mu^{-1} \curl \vec{E}  \,,
\end{equation}
which is a dimensionless system.  This system can be written as either of
two dimensionless second order equations:
\begin{equation}
\frac{d^2 \vec{E}}{d t^2} = 
- \epsilon^{-1} \curl \mu^{-1} \curl \vec{E} \,,\quad
\frac{d^2 \vec{H}}{d t^2} = 
- \mu^{-1} \curl \epsilon^{-1} \curl \vec{H}\,.
\end{equation}

Because $\epsilon$ and $\mu$ are matrices and not numbers the conserved
quantity must be changed a bit:
\begin{equation}
C = \frac{
\langle \epsilon \vec{E} , \vec{E} \rangle +
      \langle \mu \vec{H} , \vec{H} \rangle}{2}\\
\end{equation}
Because $\epsilon$ and $\mu$ are symmetric matrices
\begin{align}
\frac{d C}{d t}
& = \langle \epsilon \frac{d \vec{E}}{d t} , \vec{E} \rangle +
 \langle \mu \frac{d \vec{H}}{d t} , \vec{H} \rangle \nonumber\\
& = \langle \curl \vec{H} , \vec{E} \rangle -
      \langle \curl \vec{E} , \vec{H} \rangle \nonumber \\
& = \langle \vec{H} , \curl \vec{E} \rangle -
      \langle \curl \vec{E} , \vec{H} \rangle \nonumber \\
& = 0
\end{align}

The vector identity
\[
\grad \bdot ( \vec{E} \times \vec{H} ) =
(\curl \vec{E}) \bdot \vec{H} - \vec{E} \bdot (\curl \vec{H}) \,.
\]
can also be used to see that the energy is constant. The time derivative
of $C$ can be written using integrals as
\begin{align*}
\frac{d C}{d t}
& = \int_{\real^3}
\left( \epsilon \, \frac{d \vec{E}}{d t} \bdot \vec{E} +
\mu \, \frac{\vec{H}}{d t} \bdot \vec{H} \right)
\, dx \, dy \, dz \,, \\
   & = \int_{\real^3} \left(
	 \curl \vec{H} \bdot \vec{E}
	-\curl \vec{E} \bdot \vec{H}
	  \right) \, dx \, dy \, dz \,, \\ 
   & = - \int_{\real^3} \divg ( \vec{H} \times \vec{E} ) \, dx \, dy \, dz \,, \\
   & = \int_{\real^3} \divg ( \vec{E} \times \vec{H} ) \, dx \, dy \, dz \,, \\
   & = \int_{\real^3} \divg \vec{S} \, dx \, dy \, dz \,, \\
   & = 0 \,.
\end{align*}
The last integral is zero because it was assumed that $\vec{E}$ and $\vec{H}$
are zero far from the origin. Also $\vec{S} = \vec{E} \times \vec{H} $ is
called the Poynting vector which has spatial units $1/d^2$.
The integrand is the standard energy density confirming $C$ is spatially
dimensionless.

\newpage \clearpage
\setcounter{equation}{0}

\section{Variable Coefficient Differential Operators
\label{Second Order D0s}}

The first task is to describe all second order operators with variable
coefficients that can be generated using diagram chasing in the double exact
sequence shown in Figure \ref{Exact-Sequences}.  Next inner products are
introduced on the $H$ spaces and the adjoints of all of the operators in the
double exact sequence are computed. This is then used to show that the second
order operators are self-adjoint and either positive or negative. The second
order operators and the adjoint operators can be found quickly by diagram
chasing as will be described. In the continuum setting there are two
equivalent derivations of each operator.  In the discrete setting these will
be different.

\subsection{Exact Sequences}

\begin{figure}
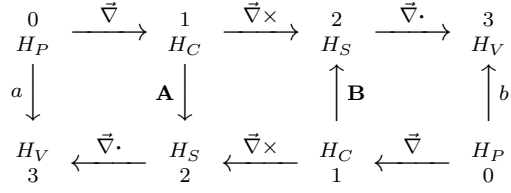

\begin{equation*}
  \begin{CD}
		   {\mystack{0}{H_P}}
        @>\grad >> {\mystack{1}{H_C}}
        @>\curl >> {\mystack{2}{H_S}}
        @>\divg >> {\mystack{3}{H_V}}
\\
@V{a}VV @V{{\bf A}}VV @AA{{\bf B}}A @AA{b}A @. \\
                   {\mystack{H_V}{3}}
        @<\divg << {\mystack{H_S}{2}}
        @<\curl << {\mystack{H_C}{1}}
        @<\grad << {\mystack{H_P}{0}}
\\
  \end{CD}
\end{equation*}
\caption{Continuum Double Exact Sequence Diagram}
\label{Exact-Sequences}
\end{figure}

In the double exact sequence diagram \ref{Exact-Sequences} the bottom row is
the same as the top row written in opposite order. For reasons which will
become clear when the operators are discretized,
$P$ stands for points,
$C$ stands for curves,
$S$ stands for surfaces,
$V$ stands for volumes.
In this diagram $H_P$ and $H_V$ are linear spaces of smooth scalar functions
depending on the spatial variables $(x,y,z)$ and $H_C$ and $H_S$ are linear
spaces of smooth vector functions that depend also on $(x,y,z)$.
All of the functions converge rapidly to zero as $x^2+y^2+z^2$ becomes large.
The first order differential operators are the gradient $\grad$, curl or
rotation $\curl$, and divergence $\divg$.
The scalar functions $a$ and $b$ and also the matrix valued functions
${\bf A}$ and ${\bf B}$ are also smooth function of the spatial variables
that are used to describe material properties.
The functions $a$ and $b$ are bounded above and below by positive constants.
The matrix functions are symmetric positive definite and the eigenvalues
of the matrices are bound above and below by positive constants.

\begin{table}[ht]
\begin{center}
\begin{tabular}{|l|l|}
\hline
$\text{if }f \in H_P$        $\text{ then } \grad f \in H_C$ &
$\text{if }f \in H_P$        $\text{ then } a \, f \in H_V$ \\
$\text{if }\vec{v} \in H_C$  $\text{ then } \curl \vec{v} \in H_S$ &
$\text{if }\vec{v} \in H_C$  $\text{ then } {\bf A} \, \vec{v} \in H_S$ \\
$\text{if }\vec{w} \in H_S$  $\text{ then } \divg \vec{w} \in H_V$ &
$\text{if }\vec{w} \in H_C$  $\text{ then } {\bf B} \, \vec{w} \in H_S$ \\
& $\text{if }g \in H_P$        $\text{ then } b \, g \in H_V$  \\
\hline
\end{tabular}
\caption{First order operators on the left, material property operators on
the right.}
\label{Basic Operators}
\end{center}
\end{table}

The differential operators and material property functions give mappings
between the spaces in the double exact sequence as described in Table
\ref{Basic Operators}. Note that the differential operators are not
invertible, but that the conditions on the material properties function
imply that they give invertible mappings.
The horizontal arrows represent the action of the differential operators
while the vertical arrows represent multiplication by scalar functions $a$ and
$b$ and by $3 \times 3$ matrices ${\bf A}$ and ${\bf B}$ that are known as
{\em star} operators in differential geometry.  The directions of the
vertical arrows in the double exact sequence can be chosen to be
either up or down as the operators are invertible.

The integers in the double exact sequence give the spatial dimension
of the function in the spaces, that is if
$f \in H_P$ the $f$ has no spatial dimension while if $g \in H_V$ then
$g$ has spatial dimension $1/d^3$.
Also if $\vec{v} \in H_C$ then $\vec{v}$ has spatial dimension $1/d$
and if $\vec{w} \in H_S$ then $\vec{w}$ has spatial dimension $1/d^2$.
The differential operators all have dimension $1/d$. Moreover 
$a$ and $b$ have spatial dimensions $1/d^3$ while $\bf A$ and $\bf B$
have dimensions $1/d$. The directions of the vertical arrows were chosen
so that $a$, $b$, $\bf A$ and $\bf B$ had dimensions $1/d^k$ for $k>0$.

Importantly, all of the operators given Table \ref{Basic Operators} are
dimensionally consistent and all of the spaces in the double exact sequence

The most important properties of exactness are that $\divg \, \curl = 0$ and
$\curl \, \grad = 0$. Discretization that violate either of these two conditions
are not mimetic!  Additionally, exactness requires the existence of
local scalar and vector potentials.  That is, if $\vec{v} \in H_C$ and
$\curl \vec{v} = 0$ then there an $f \in H_P$ so that $\grad f = \vec{v}$
and if $\vec{w} \in H_S$ and $\divg \vec{w} = 0$ then there is an
$\vec{v} \in H_C$ so that $ \vec{w} = \curl \vec{v}$.

\subsection{Diagram Chasing and Second Order Operators}

\begin{table}[ht]
\begin{center}
\begin{tabular}{|r|r||r|r|}
\hline
\multicolumn{2}{|c|}{Upper Row} &
\multicolumn{2}{|c|}{Bottom Row} \\
\hline
\multicolumn{4}{|c|}{First Box} \\
\hline
  $f \in H_P$                       
& $\vec{v} \in H_C$
& $\vec{w} \in H_S$
& $g \in H_V$ \\
  $\grad f \in H_C$                 
& ${\bf A} \vec{v} \in H_S$
& $\divg \vec{w} \in H_V$
& $a^{-1} g \in H_P$ \\
  ${\bf A}  \grad f \in H_S$              
& $\divg {\bf A} \vec{v} \in H_V$
& $a^{-1} \divg \vec{w} \in H_P$
& $\grad a^{-1} g \in H_C$ \\
  $\divg {\bf A}  \grad f \in H_V$        
& $a^{-1} \divg {\bf A} \vec{v} \in H_P$
& $\grad a^{-1} \divg \vec{w} \in H_C$
& ${\bf A} \grad a^{-1} g \in H_S$ \\
  $a^{-1} \divg {\bf A}  \grad f \in H_P$ 
& $\grad a^{-1} \divg {\bf A} \vec{v} \in H_C$
& ${\bf A} \grad a^{-1} \divg \vec{w} \in H_S$
& $\divg {\bf A} \grad a^{-1} g \in H_V$ \\
\hline 
\multicolumn{4}{|c|}{Second Box} \\
\hline
  $\vec{v} \in H_C$                       
& $\vec{w} \in H_S$
& $\vec{v} \in H_C$
& $\vec{w} \in H_S$ \\
  $\curl \vec{v} \in H_S$                       
& ${\bf B}^{-1} \vec{w} \in H_C$
& $\curl \vec{v} \in H_S$
& ${\bf A}^{-1} \vec{w} \in H_C$ \\
  ${\bf B}^{-1} \curl \vec{v} \in H_C$                       
& $\curl {\bf B}^{-1} \vec{w} \in H_S$
& ${\bf A}^{-1} \curl \vec{v} \in H_C$
& $\curl {\bf A}^{-1} \vec{w} \in H_S$ \\
  $\curl {\bf B}^{-1} \curl \vec{v} \in H_S$                       
& ${\bf A}^{-1} \curl {\bf B}^{-1} \vec{w} \in H_C$
& $\curl {\bf A}^{-1} \curl \vec{v} \in H_S$
& ${\bf B}^{-1} \curl {\bf A}^{-1} \vec{w} \in H_C$ \\
  ${\bf A}^{-1} \curl {\bf B}^{-1} \curl \vec{v} \in H_C$
& $\curl {\bf A}^{-1} \curl {\bf B}^{-1} \vec{w} \in H_S$
& ${\bf B}^{-1} \curl {\bf A}^{-1} \curl \vec{v} \in H_C$
& $\curl {\bf B}^{-1} \curl {\bf A}^{-1} \vec{w} \in H_S$ \\
\hline 
\multicolumn{4}{|c|}{Third Box} \\
\hline
  $\vec{w} \in H_S$
& $g \in H_V$
& $f \in H_P$
& $\vec{v} \in H_C$ \\
  $\divg \vec{w} \in H_V$
& $b^{-1} g \in H_P$
& $\grad f \in H_C$
& ${\bf B}^{-1} \vec{v} \in H_S$ \\
  $b^{-1} \divg \vec{w} \in H_P$
& $\grad {b^{-1}} g \in H_C$ 
& ${\bf B} \grad f \in H_S$
& $\divg {\bf B}^{-1} \vec{v} \in H_V$ \\
  $\grad b^{-1} \divg \vec{w} \in H_C$
& ${\bf B} \grad b^{-1} g \in H_S$
& $\divg {\bf B} \grad f \in H_V$
& $b^{-1} \divg {\bf B}^{-1} \vec{v} \in H_P$ \\
  ${\bf B} \grad b^{-1} \divg \vec{w} \in H_S$
& $\divg {\bf B} \grad b^{-1} g \in H_V$
& $b^{-1} \divg {\bf B} \grad f \in H_P$
& $\grad b^{-1} \divg {\bf B}^{-1} \vec{v} \in H_C$ \\
\hline 
\end{tabular}
\caption{Fundamental second order differential operators. The
left two columns start with spaces in the top row in Figure
\ref{Exact-Sequences}
while the right two columns start with spaces in the bottom row.}
\label{General Fundamental Operators}
\end{center}
\end{table}

Table \ref{General Fundamental Operators} gives all of the possible second
order operators given by diagram chasing. As an example of diagram chasing,
consider
$f \in H_P$ so that
$\grad f \in H_C$ and then
${\bf A} \grad f \in H_S$ so that
$\divg {\bf A} \grad f \in H_V $ and finally
$a^{-1} \divg {\bf A} \grad f \in H_S $.
The gives the upper left entry in Table \ref{General Fundamental Operators}.
The remaining operators are created similarly.  For diagram chasing it is
important that the mappings $a$, $b$, ${\bf A}$ and ${\bf B}$ are invertible
while $\grad$, $\curl$ and $\divg$ are not invertible. Consequently to create
a second order operator, only going clockwise around a square in 
Figure \ref{Exact-Sequences} is allowed. However it is possible to start and
any corner, so this gives twelve operators, four corners times three squares.
In the continuum some of these operators are essentially the same, for
example $a^{-1} \divg {\bf A}  \grad$ and $b^{-1} \divg B \grad$. The
assumption that ${\bf B}= {\bf A}$ and $b = a$ reduces the number of operators
to six. Finally, if $a = b = 1$ and ${\bf A} = {\bf B} = I$ the identity matrix
then the operators simplify to those in Table
\ref{Simplified General Operators} that is there are only three distinct
second order operators:
\[
\lap f = \divg \grad f \,; \quad \quad
\curl \curl \vec{v} \,; \quad \quad
\grad \divg  \vec{w} \,.
\]

\begin{table}[ht]
\begin{center}
\begin{tabular}{|r|r|r|r|}
\hline
\multicolumn{4}{|c|}{First Box} \\
\hline
  $f \in H_P$                       
& $\vec{v} \in H_C$
& $\vec{w} \in H_S$
& $g \in H_V$ \\
  $ \divg  \grad f \in H_P$ 
& $\grad \divg \vec{v} \in H_C$
& $\grad \divg \vec{w} \in H_S$
& $\divg \grad  g \in H_V$ \\
\hline 
\multicolumn{4}{|c|}{Second Box} \\
\hline
  $\vec{v} \in H_C$                       
& $\vec{w} \in H_S$
& $\vec{v} \in H_C$
& $\vec{w} \in H_S$ \\
  $\curl \curl \vec{v} \in H_C$                       
& $\curl \curl \vec{w} \in H_S$
& $\curl \curl \vec{v} \in H_C$
& $\curl \curl \vec{w} \in H_S$ \\
\hline 
\multicolumn{4}{|c|}{Third Box} \\
\hline
  $\vec{w} \in H_S$
& $g \in H_V$
& $f \in H_P$
& $\vec{v} \in H_C$ \\
  $ \grad \divg \vec{w} \in H_S$
& $\divg \grad g \in H_V$
& $\divg \grad f \in H_P$
& $\grad \divg \vec{v} \in H_C$ \\
\hline 
\end{tabular}
\caption{The operators can be simplified by assuming that
$a = b = 1$ and that ${bf A} = {\bf B} ={\bf I}$, the identity
matrix.
}
\label{Simplified General Operators}
\end{center}
\end{table}

\subsection{Additional Second Order Operators}

Note that in Table \ref{General Fundamental Operators}
there are two operators defined on $H_P$, four operators defined on $H_C$,
four operators defined on $H_S$ and two operators defined on $H_V$. If
linear operators are defined on the same space then linear combinations
of these operators are again linear operators. The two operators defined
on $H_P$ and the two defined on $H_V$ are essentially the same so linear
combinations are not interesting. For any operator in the {\em Top Row}
boxes, there is a similar operator in the {\em Bottom Row} that can
be obtained by interchanging $a$ with $b$ and ${\bf A}$ with ${\bf B}$.
For $\vec{v} \in H_C$ and for $\vec{w} \in H_S$ define
\begin{align}
{\bf VL}_1 (\vec{v}) & = 
\grad a^{-1} \, \divg \, {\bf A} \, \vec{v} -
{\bf A}^{-1} \, \curl \, {\bf B}^{-1} \, \curl \, \vec{v} \, \nonumber \\
{\bf VL}_2 (\vec{w}) & = 
{\bf B} \, \grad b^{-1} \, \divg \, \vec{w} -
\curl \, {\bf A}^{-1} \, \curl \, {\bf B}^{-1} \, \vec{w} \, \nonumber \\
{\bf VL}_3 (\vec{v}) & =
\grad \, b^{-1} \, \divg \, {\bf B} \, \vec{v} - 
{\bf B}^{-1} \, \curl \, {\bf A}^{-1} \, \curl \, \vec{v} \label{Vector Laplacians} \\
{\bf VL}_4 (\vec{w}) & =
{\bf A} \, \grad \, a^{-1} \, \divg \, \vec{w}-
\curl \, {\bf B}^{-1} \, \curl \, {\bf A}^{-1} \, \vec{w} \nonumber \\
\nonumber
\end{align}

Under the simplifying assumptions that $a = b = 1$ and
${\bf A = \bf B = \bf I}$ these operators become
\[
\lapv \vec{v} = \grad \divg \vec{v} - \curl \curl \vec{v} \,,\quad
\lapv \vec{w} = \grad \divg \vec{w} - \curl \curl \vec{w} \,,
\]
which in Cartesian coordinates gives (see {\tt CurlCurl.nb})
the vector Laplacian
\[
\lapv \left(v_1,v_2,v_3\right) = \left(\lap v_1,\lap v_2,\lap v_3\right) \,,\quad
\lapv \left(w_1,w_2,w_3\right) = \left(\lap w_1,\lap w_2,\lap w_3\right) \,.
\]
Operators like these appear in the elastic wave equation
\ref{Elastic Wave Equation}.

\subsection{Inner Products}

Applying the mimetic ideas to physical problems requires the use of inner
products on the spaces $H_P$, $H_C$, $H_S$ and $H_V$. It is important that
the inner products do not have a spatial dimension.  Two bilinear forms
will help simplify the notations. Let
$f = f(x,y,z) \in H_P$,\,
$g = g(x,y,z) \in H_V$,\,
$\vec{v} = \vec{v(x,y,z)} \in H_C$
and
$\vec{v} = \vec{w(x,y,z)} \in H_S$ and then define 
\begin{align}
\lbilin f , g \rbilin & = 
\int_{\real^3} f(x,y,z) \, g(x,y,z) \, dx \, dy \, dz \,, \nonumber \\
\lbilin \vec{v} , \vec{w} \rbilin & = 
\int_{\real^3} \vec{v}(x,y,z) \bdot \vec{w}(x,y,z) \, dx \, dy \, dz \,.
\label{bilinear form}
\end{align}
These bilinear forms are dimensionless because $dx$, $dy$ and $dy$ have
dimension $d$ while $f$ has dimension $0$, $g$ has dimension $1/d^3$,
$\vec{v}$ has dimension $1/d$ and $\vec{w}$ has $1/d^2$.

The inner product on the function spaces must use a weight function to
be dimensionless: 
\begin{itemize}
\item
for $f_1 \,, f_2 \in H_P$ set
$\langle f_1 , f_2 \rangle_P = \lbilin a \, f_1 , f_2 \rbilin$
\item
for $\vec{v}_1 \,, \vec{v}_2 \in H_C$ set
$ \langle \vec{v}_1 , \vec{v}_2 \rangle_C = 
\lbilin {\bf A} \vec{v}_1 , \vec{v}_2 \rbilin$
\item
for $\vec{w}_1 \,, \vec{w}_2 \in H_S$ set
$\langle \vec{w}_1 , \vec{w}_2 \rangle_S = 
\lbilin {\bf A}^{-1} \, \vec{w}_1 , \vec{w}_2 \rbilin$ 
\item
for $g_1 \,, g_2 \in H_V$ set
$\langle g_1 , g_2 \rangle_V = 
\lbilin a^{-1} \, g_1 , g_2 \rbilin$
\end{itemize}
As usual
$\norm{f}_P^2 = \langle f , f \rangle_P$,
$\norm{\vec{v}}_C^2 = \langle \vec{v} , \vec{v} \rangle_C$,
$\norm{\vec{w}}_S^2 = \langle \vec{w} , \vec{w} \rangle_S$,
$\norm{g}_V^2 = \langle g , g \rangle_V$.
Additional inner products can be made by replacing $a$ by $b$ and
${\bf A}$ by ${\bf B}$. To be inner products it is important that
$a$ and $b$ are positive and that ${\bf A}$ and ${\bf B}$ are
symmetric and positive definite matrices.

\subsection{Adjoint Operators \label{Continuum Adjoint Operators}}

Adjoints are commonly defined for operators mapping a space into itself
but most of the operators used here are mapping between two different
spaces, so the adjoint is defined as in Section \ref{Continuous Time}.
The discussion in that section shows that if $X$, $Y$, and $Z$ are 
linear spaces and $O_1$ and $O_2$ are linear operator such that:
\[
X \overset{O_1}{\rightarrow} Y \overset{O_2}{\rightarrow} Z
\quad
\text{then }
Z \overset{O_2^*}{\rightarrow}
Y \overset{O_1^*}{\rightarrow} X \,
\]
and 
\[
(O_1 O_2)^* = O_2^* \, O_1^* \,,\quad  (O_1^*)^* = O_1 \,.
\]
Because diagram chasing gives operators as compositions of other
operators, this will be used many times.

The adjoints of the operators in Table \ref{Basic Operators} are
\begin{align}
\grad^* & = - a^{-1} \, \divg \, {\bf A} &
\grad^* & = - b^{-1} \, \divg \, {\bf B}
\,,\nonumber\\
\curl^* & = + {\bf A}^{-1} \, \curl \, {\bf B}^{-1} &
\curl^* & = + {\bf B}^{-1} \, \curl \, {\bf A}^{-1}
\,,\nonumber\\
\divg^* & = - {\bf B} \grad  b^{-1} &
\divg^* & = - {\bf A} \grad  a^{-1} 
\,,\nonumber\\
{\bf A}^* & = {\bf A}^{-1} &
{\bf B}^* & = {\bf B}^{-1}
\,, \label{Adjoints} \\
a^* & = a^{-1} &
b^* & = b^{-1}
\nonumber\,,
\end{align}
where the column on the left contains differential operators from the top
row in Figure \ref{Exact-Sequences} and the column on the right contains
differential operators from the bottom row in Figure \ref{Exact-Sequences}.

Now the proofs of the formulas in  \ref{Adjoints}  are straight forward.
For the gradient let $f \in H_P$ and $\vec{v} \in H_C$ so that
\begin{align*}
\langle \grad f , \vec{v} \rangle_C 
 & = \lbilin {\bf A} \, \grad f , \vec{v} \rbilin \\
 & = \lbilin \grad f , {\bf A} \, \vec{v} \rbilin \\
 & = - \lbilin f , \divg \, {\bf A} \vec{v} \rbilin \\
 & = - \lbilin a \, f , a^{-1} \, \divg \, {\bf A} \vec{v} \rbilin \\
 & =   \langle f , - a^{-1} \, \divg \, {\bf A} \vec{v} \rangle_P
\end{align*}
For the curl let $\vec{v} \in  H_C$ and $\vec{w} \in H_S$ so that
\begin{align*}
\langle \curl \vec{v} , \vec{w} \rangle_S 
 & = \lbilin {\bf A}^{-1} \, \curl \vec{v} , \vec{w} \rbilin \\
 & = \lbilin \curl \vec{v} , {\bf A}^{-1} \, \vec{w} \rbilin \\
 & = \lbilin \vec{v} , \curl \, {\bf A}^{-1} \vec{w} \rbilin \\
 & = \lbilin {\bf B} \, \vec{v} ,
	{\bf B}^{-1} \, \curl \, {\bf A}^{-1} \vec{v} \rbilin \\
 & = \langle \vec{v} , {\bf B}^{-1} \, \curl \, {\bf A}^{-1} \vec{w} \rangle_C
\end{align*}
For the divergence let $\vec{w} \in H_S$ and $ g \in H_V$ so that
\begin{align*}
\langle \divg \vec{w} , g \rangle_V
 & = \lbilin a^{-1} \divg \vec{w} , g \rbilin  \\
 & = \lbilin \divg \vec{w} , a^{-1} g \rbilin  \\
 & = - \lbilin \vec{w} , \grad ( a^{-1} g ) \rbilin  \\
 & = - \lbilin {\bf A}^{-1}  \vec{w} , {\bf A} \grad ( a^{-1} g ) \rbilin  \\
 & =   \langle \vec{w} , - {\bf A} \grad ( a^{-1} g ) \rangle_S  \\
\end{align*}
For the operator ${\bf A}$ let $\vec{v} \in H_C$ and $\vec{w} \in H_S$ so that
\begin{align*}
\langle {\bf A} \vec{v} , \vec{w} \rangle_S  
 & = \lbilin {\bf A}^{-1} {\bf A} \vec{v} , \vec{w} \rbilin \\
 & = \lbilin {\bf A} \vec{v} , {\bf A}^{-1} \vec{v} \rbilin \\
 & = \langle \vec{v} , {\bf A}^{-1} \vec{v} \rangle_C
\end{align*}
For the operator $a$ let $ f \in H_P$ and $g \in H_V$ so that
\begin{align*}
\langle a \, f , g \rangle_V 
 & = \lbilin a^{-1} a \, f , g \rbilin  \\
 & = \lbilin a \, f , a^{-1} g \rbilin  \\
 & = \langle f , a^{-1} g \rangle_P  \\
\end{align*}
Similar arguments give the adjoint operators for operators containing
$b$ and ${\bf B}$.
To keep the notation easy to read it has not been specified whether to use
$a$ or $b$ and whether to use ${\bf A}$ or ${\bf B}$ in the inner products
when computing adjoints. This is clear from the context.

It is now straight forward to compute the adjoints of the second order
operators in Table \ref{General Fundamental Operators}:
\begin{align}
\left( a^{-1} \divg {\bf A} \grad \right)^* 
& = \grad^* \, {\bf A}^* \, \divg^* \, {a^{-1}}^* \nonumber\\
& = a^{-1}\, \divg\, {\bf A}\, {\bf A}^{-1} \,
	{\bf A} \, \grad \, a^{-1} \, a \label{Example 1}\\
& = a^{-1} \, \divg \, {\bf A} \, \grad \,. \nonumber
\end{align}
\begin{align}
\left({\bf A}^{-1} \, \curl \, {\bf B}^{-1} \, \curl \, \right)^*
& = \curl^* \, {{\bf B}^{-1}}^* \, \curl^* \, {{\bf A}^{-1}}^* \nonumber\\
& = {\bf A}^{-1} \, \curl \, {\bf B}^{-1} \, {\bf B} \,
	{\bf B}^{-1} \, \curl \, {\bf A}^{-1} \, {\bf A} \label{Example 2} \\
& = {\bf A}^{-1} \, \curl \, {\bf B}^{-1} \, \curl \, . \nonumber
\end{align}
\begin{align}
\left( {\bf B} \, \grad b^{-1} \, \divg \right)^*
& = \divg^* \, {b^{-1}}^* \, \grad^* \, {\bf B}^* \nonumber \\
& = {\bf B} \, \grad \, b^{-1} \, b \, b^{-1} \,
	\divg \, {\bf B} \, {\bf B}^{-1} \label{Example 3} \\
& = {\bf B} \, \grad \, b^{-1} \, \divg  \nonumber
\end{align}
So these three operators are self-adjoint and similar arguments show that all
operators in Table \ref{General Fundamental Operators} are self-adjoint.

\subsection{Positive and Negative Second Order Operators}

Arguments like those in the previous sections can be use to show that
the second order operators in \ref{General Fundamental Operators} are
either positive or negative, those that contain two curl operators are
positive while those that contain a gradient and divergence are negative.

Let $f \in H_P$ and then consider
\begin{align*}
\langle a^{-1} \, \divg \, {\bf A} \, \grad \, f, f \rangle_P
& = \langle \divg \, {\bf A} \, \grad \, f, a \, f \rangle_V
\\
& = - \langle {\bf A} \, \grad \, f, {\bf A} \, \grad \, a^{-1} \, a \,f \rangle_S \\
& = - \langle {\bf A} \, \grad \, f, {\bf A} \, \grad \, f \rangle_S \\
& \leq 0
\end{align*}
Let $\vec{v} \in H_C$ and then consider
\begin{align*}
\langle {\bf A}^{-1} \, \curl \, {\bf B}^{-1} \, \curl \, \vec{v} , \vec{v} \rangle_C
& = \langle  \curl \, {\bf B}^{-1} \, \curl \, \vec{v} , {\bf A} \, \vec{v} \rangle_S \\
& = \langle {\bf B}^{-1} \, \curl \, \vec{v} , {\bf B}^{-1} \, \curl \, {\bf A}^{-1} {\bf A} \vec{v} \rangle_C \\
& = \langle {\bf B}^{-1} \, \curl \, \vec{v} , {\bf B}^{-1} \, \curl \, \vec{v} \rangle_C \\
& \geq 0
\end{align*}
Let $\vec{w} \in H_S$ and then consider
\begin{align*}
\langle {\bf B} \, \grad b^{-1} \, \divg \, \vec{w} , \vec{w} \rangle_S
& = \langle \grad b^{-1} \, \divg \, \vec{w} , {\bf B}^{-1} \, \vec{w} \rangle_S \\
& = 
\langle b^{-1} \, \divg \, \vec{w} , \grad^* \, {\bf B}^{-1} \, \vec{w} \rangle_P \\
& = - \langle b^{-1} \, \divg \, \vec{w} , b^{-1} \, \divg \, {\bf B} \, {\bf B}^{-1} \, \vec{w} \rangle_P \\
& = - \langle b^{-1} \, \divg \, \vec{w} , b^{-1} \, \divg \, \vec{w} \rangle_P \\
& \leq 0
\end{align*}
These results capture the important features of the homogeneous and
isotropic second order differential operators.

\newpage \clearpage
\setcounter{equation}{0}
\section{Wave Equations With Variable
Material Properties\label{Wave Equations With Variable Materials}}
This section describes all of the second order wave equations that can
be generated using the operators that were generated using diagram chasing
in Section \ref{Second Order D0s}. These wave equations have the form 
\[
\frac{d^2 W}{d t^2} = \mathscr{A} \, W \,,
\]
where $W = W(t,x,y,z)$ is a scalar or vector function and where $\mathscr{A}$
is $\pm$ an operator from Table \ref{General Fundamental Operators} with the
$+$ or $-$ chosen so that $\mathscr{A}$ is negative operator.  This provides
12 possible wave equations. Another four equations are obtained using the
operators in \ref{Vector Laplacians}.  Many of these equations are
equivalent but will not be equivalent in the discrete setting.
All of these equations are dimensionless as the map a space onto it
self so the spatial dimensions of the functions play no role.

First six types of wave equations are introduced and then reduced to a easily
recognized form for uniform material properties given by using the simplifying
assumptions that
$a = 1/c^3$, $b = 1/c^3$, ${\bf A} = {\bf I}/c$ and ${\bf B} = {\bf I}/c$.
These assumptions keep the spatial dimensions correct.

Next first order systems are created using operators generate by going half way
around the squares in the double exact sequence diagram \ref{Exact-Sequences}.
Note that these first order systems are not dimensionless because they map
between space containing function with different spatial units. The first
order systems are used to derived to create conservation laws for the first
systems and consequently for the second order equations.

\subsection{Scalar and Vector Wave Equations}

To create a scalar wave equations choose $f = f(x,y,z,t) \in H_P$
and then define
\begin{equation}\label{Second Order 1}
\frac{d^2 f}{d t^2} = a^{-1} \divg {\bf A} \grad  f \,.
\end{equation}
Setting $f = a^{-1} \, g$ with $g \in H_V$ gives
\begin{equation}\label{Second Order 2}
\frac{d^2 g}{d t^2} = \divg {\bf A} \grad a^{-1} g \,.
\end{equation}
Under the simplifying assumptions these become the standard scalar wave 
equation
\[
\frac{d^2 f}{d t^2} = c^2 \, \divg \, \grad  f \,,\quad
\frac{d^2 g}{d t^2} = c^2 \, \divg \, \grad  g \,.
\]

To create a vector wave equation choose $\vec{v} \in H_C$ to get
\begin{equation}\label{Second Order 3}
\frac{d^2 \vec{v}}{d t^2} =
- {\bf A}^{-1} \, \curl \, {\bf B}^{-1} \, \curl \, \vec{v} \,.
\end{equation}
Setting $\vec{v} = {\bf A}^{-1} \vec{w}$ with $\vec{w} \in H_S$ gives
\begin{equation}\label{Second Order 4}
\frac{d^2 \vec{w}}{d t^2} =
- \curl \, {\bf B}^{-1} \, \curl \, {\bf A}^{-1} \, \vec{w} \,.
\end{equation}
Under the simplifying assumptions these become
\[
\frac{d^2 \vec{v}}{d t^2} = - c^2 \, \curl \, \curl \, \vec{v} \,,\quad
\frac{d^2 \vec{w}}{d t^2} = - c^2 \, \curl \, \curl \, \vec{w} \,,
\]
which are Maxwell's equations \ref{Maxwell Equations} in uniform materials.

More wave equations can be generated by choosing $\vec{w} \in H_S$:
\begin{equation}\label{Second Order 5}
\frac{d^2 \vec{w}}{d t^2} =
{\bf B} \, \grad b^{-1} \, \divg \, \vec{w} \,.
\end{equation}
Setting $\vec{w} = {\bf B} \, \vec{v}$ with $\vec{v} \in H_C$ gives
\begin{equation}\label{Second Order 6}
\frac{d^2 \vec{v}}{d t^2} =
\grad \, b^{-1} \, \divg \, {\bf B} \, \vec{v} \,.
\end{equation}
Under the simplifying assumptions these become
\[
\frac{d^2 \vec{w}}{d t^2} = c^2 \, \grad \, \divg \, \vec{w} \,,\quad
\frac{d^2 \vec{v}}{d t^2} = c^2 \, \grad \, \divg \, \vec{v} \,.
\]

In total twelve equations can be created offering useful flexibility in
modeling physical problems.

Additional second order wave equations can be made from the two term
second order operators ${\bf VL}_1$, ${\bf VL}_2$, ${\bf VL}_3$, and
${\bf VL}_4$ in \ref{Vector Laplacians}. For example
\begin{equation}
\frac{d^2 \vec{v}}{d t^2} =
\grad a^{-1} \divg {\bf A} \vec{v} -
{\bf A}^{-1} \curl {\bf B}^{-1} \curl \vec{v}  \,.
\label{Curl Curl}
\end{equation}
To simplify this equation assume $a$ and $b$ are constants and 
${\bf A} = a^3 \, {\bf I}$
and
${\bf B} = 1/(a^3 \, b^2) {\bf I}$ to get
\[
\frac{d^2 \vec{v}}{d t^2} =
a^2 \, \grad \, \divg \, \vec{v} -  b^2 \, \curl \, \curl \vec{v}  \,,
\]
which is the elastic wave equation \ref{dimensionless elastic equation}.
All four equations created this way reduce to the elastic wave equation
under simplifying assumptions.

\subsection{First Order Systems and Conserved Quantities}

There is a natural way to use diagram chasing to write the second order
wave equations as a system and then use this to define a conserved quantity.
Note that the first order equations are not dimensionless as in the constant
coefficient case because the functions in the spaces the
Continuum Double Exact Sequence Diagram \ref{Exact-Sequences} are not
dimensionless. The main idea is two choose two function that in diagonally
opposite corners of one of the squares in \ref{Exact-Sequences} and do a diagram
chase.  Consequently there are lots of first order systems!

For example for equation \eqref{Second Order 1}, because $f \in H_P$,
choose $\vec{w} \in H_S$ and then set
\begin{equation}
\frac{d \vec{w}}{d t}  = {\bf A} \grad  f 
\,,\quad
\frac{d f}{d t} = a^{-1} \divg \, \vec{w}
\label{System 1}
\end{equation}
Note that $\vec{w}$ also satisfies the second order vector wave equation 
\eqref{Second Order 5} with ${\bf B}$ and $b$ replaced by
${\bf A}$ and $a$.

Table \ref{Adjoints} implies that
\[
(a^{-1} \divg )^* = \divg^* (a^{-1})^* = -A \grad a^{-1} a = -A \grad \,,
\]
so this system has the form of the equation discussed in \ref{Wave System}
and consequently should have a conserved quantity given by 
\[
C = \frac{ \norm{f}_P^2 + \norm{\vec{w}}_S^2 }{2}
\]
This can be checked explicitly:
\begin{align*}
\frac{d C}{d t}
& =
\langle f , \frac{d f}{d t} \rangle_P +
\langle \vec{w} , \frac{d \vec{w}}{d t} \rangle_S \\
& = \langle f , a^{-1} \divg \, \vec{w} \rangle_P +
\langle \vec{w} , {\bf A} \grad  f \rangle_S \\
& = \langle (a^{-1} \divg)^* f , \, \vec{w} \rangle_S +
\langle \vec{w} , {\bf A} \grad  f \rangle_S \\
& = - \langle {\bf A} \, \grad \, f , \vec{w} \rangle_S +
\langle \vec{w} , {\bf A} \grad  f \rangle_S \\
& = 0
\end{align*}
As discussed in Section \ref{Harmonic Oscillator} if $f$ and $\vec{w}$
are solutions of \eqref{System 1} then so are $df/dt$ and $d \vec{w}/dt$
and consequently the classical energy
\[
E  = \frac{\norm{\frac{d f}{d t}}_P^2
+ \norm{\frac{d \vec{w}}{d t}}_S^2}{2} 
   = \frac{\norm{\frac{d f}{d t}}_P^2
+ \norm{{\bf A} \, \grad f }_S^2}{2}
\]
is conserved.

For equation \eqref{Second Order 2}, because $g \in H_V$,
choose $\vec{v} \in H_C$ and then set
\begin{equation}
\frac{d \vec{v}}{d t}  = \grad a^{-1} \, g \,,\quad
\frac{d g}{d t}  = \divg {\bf A} \, \vec{v}\,.
\label{System 2}
\end{equation}
Again, with a change of notation this is \eqref{Second Order 5}.
A conserved quantity is given by
\[
C = \frac{ \norm{g}_V^2 + \norm{\vec{v}}_C^2 }{2}
\]
because 
\begin{align*}
\frac{d C}{d t} & =
\langle g , \frac{d g}{d t} \rangle_V +
\langle \vec{v} , \frac{d \vec{v}}{d t} \rangle_C \\
& =
\langle g , \divg {\bf A} \, \vec{v} \rangle_V +
\langle \vec{v} , \grad a^{-1} \, g \rangle_C \\
& =
\langle - {\bf A}^{-1} \,  {\bf A} \,  \grad \,  a^{-1} \,  g , \vec{v} \rangle_C +
\langle \vec{v} , \grad a^{-1} \, g \rangle_C \\
& =
- \langle \grad \,  a^{-1} \,  g , \vec{v} \rangle_C +
\langle \vec{v} , \grad a^{-1} \, g \rangle_C \\
& = 0
\end{align*}
As before
\[
E = \frac{ \norm{\frac{d g}{d t}}_V^2 +
\norm{ \frac{d\vec{w}}{d t} }_S^2 }{2} 
  = \frac{ \norm{\frac{d g}{d t}}_V^2 +
\norm{ \grad \, a^{-1} g }_S^2 }{2} 
\]
is conserved.

For equation \eqref{Second Order 3}, because $\vec{v} \in H_C$,
choose $\vec{u} \in H_C$ and then set
\begin{equation*}
\frac{d \vec{u}}{d t} = {\bf B}^{-1} \, \curl \, \vec{v} \,,\quad 
\frac{d \vec{v}}{d t} = - {\bf A}^{-1} \, \curl \, \vec{u} \,.
\end{equation*}
A conserved quantity is given by
\[
C = \frac{ \norm{\vec{v}}_C^2 + \norm{\vec{u}}_C^2 }{2}
\]
because
\begin{align*}
\frac{d C}{d t} & =
\langle \vec{v} , \frac{d \vec{v}}{d t} \rangle_C +
\langle \vec{u} , \frac{d \vec{u}}{d t} \rangle_C \\
& =
- \langle \vec{v} , {\bf A}^{-1} \, \curl \, \vec{u} \rangle_C +
\langle \vec{u} , {\bf B}^{-1} \, \curl \, \vec{v} \rangle_C \\
& = -
\langle {\bf B}^{-1} \curl {\bf A}^{-1} \, {\bf A} \vec{v} , \vec{u} \rangle_C +
\langle \vec{u} , {\bf B}^{-1} \, \curl \, \vec{v} \rangle_C \\
& = -
\langle {\bf B}^{-1} \curl
\vec{v} , \vec{u} \rangle_C +
\langle \vec{u} , {\bf B}^{-1} \, \curl \, \vec{v} \rangle_C \\
& = 0
\end{align*}
Additionally,
\[
E = \frac{\norm{\frac{d \vec{u}}{d t}}_C^2 +
          \norm{\frac{d \vec{v}}{d t}}_C^2 }{2}
  = \frac{ \norm{ {\bf B}^{-1} \, \curl \, \vec{v} }_C^2 +
	\norm{ {\bf A}^{-1} \, \curl \, \vec{u} }_C^2 }{2}
  = \frac{ \norm{\curl \, \vec{v} }_S^2 +
	\norm{\curl \, \vec{u} }_S^2 }{2} \,,
\]
is a conserved quantity.
In this case the first order system is essentially Maxwell's equations and
$C$ is essentially the physical energy.

For equation \eqref{Second Order 4}, because $\vec{v} \in H_C$,
choose $\vec{u} \in H_C$ and then set
\begin{equation*}
\frac{d \vec{u}}{d t} =
	- \curl \, {\bf A}^{-1} \, \vec{w} \,,\quad 
\frac{d \vec{v}}{d t} =
	\curl \, {\bf B}^{-1} \, \vec{u} \,.
\end{equation*}
So $\vec{v}$ and $\vec{u}$ satisfy the previous first order system with
$\vec{v}$ and $\vec{w}$ interchanged and thus has the same conserved quantity.

For equation \eqref{Second Order 5}, because $\vec{w} \in H_S$,
choose $f \in H_P$ and then set
\begin{equation*}
\frac{d f}{d t}        = b^{-1} \, \divg \, \vec{w} \,,\quad 
\frac{d \vec{w}}{d t}  =  {\bf B} \, \grad f \,.
\end{equation*}
Again $f$ satisfies \eqref{Second Order 1} with $b$ and ${\bf B}$ replaced
by $a$ and ${\bf A}$, so
\[
C = \frac{ \norm{\vec{w}}_C^2 + \norm{f}_P^2 }{2}
\]
is conserved.

For equation \eqref{Second Order 6}, because $\vec{v} \in H_C$,
choose $g \in H_V$ and then set
\begin{equation*}
\frac{d g}{d t}       = \divg \, {\bf B} \, \vec{v} \,,\quad
\frac{d \vec{v}}{d t} = \grad \, b^{-1} \, g \,.
\end{equation*}
Again $g$ satisfies \eqref{Second Order 2} with $b$ and ${\bf B}$
replaced by $a$ and ${\bf A}$ so
\[
C = \frac{ \norm{\vec{v}}_C^2 + \norm{g}_V^2 }{2}
\]
is conserved.

Also first order systems can be made from the second order equations made
from the two term second order operators
${\bf VL}_1$, ${\bf VL}_2$, ${\bf VL}_3$ and ${\bf VL}_4$ given in
\ref{Vector Laplacians}. However this requires three first
order equations.  For example, for ${\bf VL}_1$, because
$\vec{v} \in H_C$ let $g \in H_V$ and $\vec{u} \in H_C$ and then set 
\begin{align*} 
\frac{d g}{d t} & =  \divg {\bf A} \vec{v} \\
\frac{d \vec{u}}{d t} & = {\bf B}^{-1} \curl \vec{v} \\
\frac{d \vec{v}}{d t} & =
	\grad a^{-1} g - {\bf A}^{-1} \curl \vec{u} \,.
\end{align*}
For a conserved quantity set
\[
C = \frac{\norm{\vec{v}}_C^2 + \norm{\vec{u}}_C^2 + \norm{g}_V^2}{2} \\,
\]
so that 
\begin{align*}
\frac{d C}{d t} & = 
  \langle \vec{v} \, \frac{d \vec{v}}{d t} \rangle_C
+ \langle \vec{u} \, \frac{d \vec{u}}{d t} \rangle_C
+ \langle g \, \frac{d g}{d t} \rangle_V \\ 
& =
  \langle \vec{v} , \grad a^{-1} g - {\bf A}^{-1} \curl \vec{u} \rangle_C
+ \langle \vec{u} , {\bf B}^{-1} \curl \vec{v} \rangle_C
+ \langle g , \divg {\bf A} \vec{v} \rangle_V \\ 
& =
\langle \vec{v} , \grad a^{-1} g \rangle -
\langle \vec{v} , {\bf A}^{-1} \curl \vec{u} \rangle_C
+ \langle
{\bf A}^{-1} \curl {\bf B}^{-1} {\bf B} \vec{u} , \vec{v} \rangle_C
- \langle
{\bf A}^{-1} {\bf A} \grad a^{-1} g , \vec{v} \rangle_V \\ 
& =
\langle \vec{v} , \grad a^{-1} g \rangle -
\langle \vec{v} , {\bf A}^{-1} \curl \vec{u} \rangle_C
+ \langle {\bf A}^{-1} \curl \vec{u} , \vec{v} \rangle_C
- \langle
\grad a^{-1} g , \vec{v} \rangle_V \\ 
& = 0 \,.
\end{align*}
As before this implies that
\begin{align*}
E & = \frac{
		\norm{\frac{d \vec{v}}{d t}}_C^2 +
		\norm{\frac{d \vec{u}}{d t}}_C^2 +
		\norm{\frac{d g}{d t}}_V^2}{2} \\
  & = \frac{
		\norm{\frac{d \vec{v}}{ d t}}_C^2 +
		\norm{ {\bf B}^{-1} \curl \vec{v}}_C^2 +
		\norm{\divg {\bf A} \vec{v}}_V^2}{2} 
\end{align*}
is conserved.

\subsection{Examples}

Maxwell equations for electrodynamics fit into the diagram chasing paradigm
easily and will be discussed first. The general elastic wave equations are
do not fit into the diagram chasing but many special cases of elastic or
acoustic wave equations do fit well .
The notation will be changed that used in the applications.

\subsubsection{Maxwell's Equation}

\begin{figure}[ht]
\begin{equation*}
  \begin{CD}
        {\mystack{1}{H_{CU}}}
        @>\curl >> {\mystack{2}{H_{SU}}}
\\
@V{\epsilon }VV @AA{\mu}A @. \\
         {\mystack{H_{SD}}{2}}
        @<\curl << {\mystack{H_{CD}}{1}}
\\
  \end{CD}
\end{equation*}
\caption{Maxwell Exact Sequences}
\label{Maxwell-Sequences}
\end{figure}

Representing Maxwell's equations \ref{Maxwell Equations}
using diagram chasing uses the center
square in Figure \ref{Exact-Sequences} which is reproduced in Figure
\ref{Maxwell-Sequences} using notation appropriate to Maxwell's equations,
that is, by setting ${\bf A} = \epsilon$ and ${\bf B} = \mu$.
Because in Figure \ref{Maxwell-Sequences} the upper left space and the lower
right space have the same labels as well as the upper right space has the
same label as as the lower left space they are relabeled with
$U$ and $D$ standing for up and down, that is, as $H_{CU}$, $H_{CD}$, $H_{SU}$
and $H_{SD}$.  Maxwell's equations can be represented using
$\vec{E} \in H_{CU}$,
$\vec{H} \in H_{CD}$,
$\vec{B} \in H_{SU}$ and
$\vec{D} \in H_{SD}$.
Here it is assumed that $\vec{J} = 0 $, but if this is not the case then
$\vec{J} \in H_{SD}$.

To derive first order system of Maxwell's equations by diagram chasing,
start with $\vec{E} \in H_{CU}$ and then define $\vec{H} \in H_{CD}$ by
\[
\frac{d \vec{H}}{d t} = - \mu^{-1} \curl \vec{E}  \,,
\]
and then introduce 
\[
\frac{d \vec{E}}{d t} = \epsilon^{-1} \curl \vec{H} \,.
\]

The diagram chasing implies that a conserved quantity is given by
\begin{equation}
C = \frac{\langle \vec{E} , \vec{E} \rangle_{CU} +
      \langle \vec{H} , \vec{H} \rangle_{CD} }{2}\\
\end{equation}
It is easy to check that this quantity is conserved and is the
energy for the Maxwell equations, see Section \ref{Section Maxwell Equations}.

\subsubsection{General Elastic Wave Equations}

This section is based on the discussion in \cite{Etgen87} and a summary of the
notation is given in Table \ref{Elastic Units}.  When there are no external
forces the general elastic wave equation in a material with spatially variable
properties is given by Newton's law applied to the displacements of the
material:
\begin{equation}
\rho \frac{d^2 u_i}{d t^2} =
\sum_{j=1}^3 \frac{ d{\sigma_{i,j}}}{ d {x_j}}
	\,,\quad 1 \leq i \leq 3
\label{General Elastic Wave Equation}
\end{equation} 
where $t$ is time, $\vec{u} = \vec{u}(\vec{x},t)$, are the displacements
of the material, $\rho = \rho(\vec{x})$ is the density of the material,
and $\sigma = \sigma(\vec{x})$ is the symmetric stress tensor:
\[
\sigma_{i,j} = \sum_{k,l = 1}^3 C_{k,l,i,j} \, e_{k,l} \,.
\]
This equation has the same form as the general wave equation in
\ref{General Wave Equation}.
The strains are 
\[
e_{i,j} = \frac{1}{2} 
\left(
\frac{\partial u_i}{\partial x_j} + \frac{\partial u_j}{\partial x_i}
\right) 
\]
which are dimensionless.
The material properties other than density are given by the
$C_{k,l,i,j}$ where $C = C(\vec{x})$ and where $C$ has the
symmetries $C_{k,l,i,j} = C_{l,k,i,j} = C_{k,l,j,i}$.
Consequently $C$ has only 21 independent entries \cite{Etgen87}.

An important message from this is that the right had side of
\ref{General Elastic Wave Equation}
is the sum of many terms that have a form 
\[
\frac{\partial}{\partial x_2} \, C_{1211} (\vec {x}) \,
\frac{\partial u_1}{\partial x_1} \,,
\]
that is, the elastic parameters $C_{k,l,i.j}$ appear between two first
derivatives. If $C$ does depend on $\vec{x}$ then any other form of such
terms must contain a derivative of $C$.

To fit the general elastic wave equation \ref{General Elastic Wave Equation}
into the mimetic framework there are two obvious choices: divide the equation
by $\rho$ or set $ \vec{w} = \rho \vec{v}$ with $\vec{w} \in H_S$ because
the displacements have dimension $d$. The problem with the first choice
is that none of the vectors in \ref{Vector Laplacians} have a scalar 
or vector multiplier. For the second choice both ${\bf VL}_2$ and
${\bf VL}_4$ could work. However these operator only contain 9 parameters
which are is far less than the 21 parameters in the elastic wave 
equation \ref{General Elastic Wave Equation}.

Using ${\bf VL}_4$ gives the most general mimetic wave equation as
\begin{equation}
\frac{\partial^2 \vec{w}}{\partial t^2} = 
{\bf A} \, \grad \, a^{-1} \, \divg \, \vec{w} -
\curl \, {\bf B}^{-1} \, \curl \, {\bf A}^{-1} \, \vec{w} \,.
\label{Mimetic Elastic Wave Equation}
\end{equation}
However, this is not in the form of the elastic wave equation unless
${\bf A}= I$ the identity matrix. So in fact this equation only has
7 parameters.

For an alternative related approach that works see \cite{arnoldfalk06}.

\newpage \clearpage
\setcounter{equation}{0}
\section{Mimetic Discretizations \label{Mimetic Discretizations}}

This discussion and notation will follow that in \cite{RobidouxSteinberg2011}.
However, that work was set up to rigorously prove that the discrete operators
in mimetic discretizations have the same properties as the continuum operators
used in vector calculus. Here the focus will be on applying mimetic methods to
physical problems by adding a time variable and its discretization and focus
on how to use physical spatial units to correctly discretize physical problems.

\subsection{Primal and Dual Grids}

\begin{figure}
\begin{center}
\begin{tabular}{c}
\includegraphics[width=4.00in,trim = 0 50 0 150,clip]{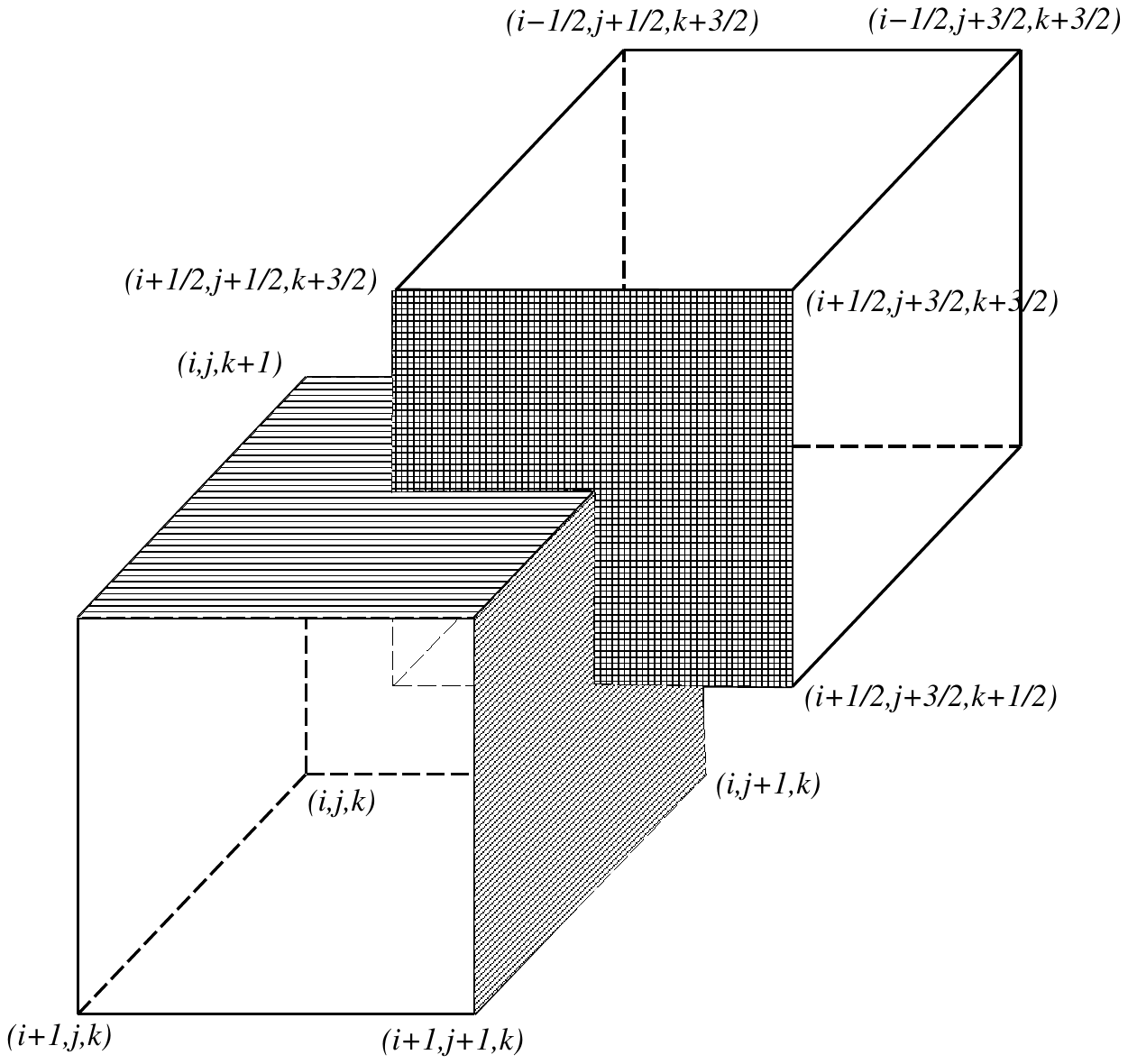}
\end{tabular}
\end{center}
\caption{The Primal and Dual Grids Taken From \cite{RobidouxSteinberg2011}
\label{Figure-Dual-Grid}}
\end{figure}

\begin{table}
\begin{center}
{\renewcommand{\arraystretch}{1.4}
\begin{tabular}{|c|c|c|}
\hline
primal   & & dual     \\
\hline
nodes   &$\left( i\,\Dx,j\,\Dy,k\,\Dz \right)$   & cells  \\
\hline
        &$\left( \left( i+\half \right) \,\Dx, j \,\Dy, k \,\Dz \right)$          &          \\
edges   &$\left( i \,\Dx, \left( j+\half \right) \,\Dy, k \,\Dz \right)$          & faces    \\
        &$\left( i \,\Dx, j \,\Dy, \left( k+\half \right) \,\Dz \right)$          &          \\
\hline
        &$\left( i \Dx, \left( j+\half \right) \Dy, \left( k+\half \right) \Dz \right) $    &          \\
faces   &$\left( \left( i+\half \right) \Dx, j \Dy, \left( k+\half \right) \Dz \right) $    & edges    \\
        &$\left( \left( i+\half \right) \Dx, \left( j+\half \right) \Dy , k \Dz \right)$    &          \\
\hline
cells   &$\left( \left( i+\half \right) \Dx, \left( j+\half \right) \Dy, \left( k+\half \right) \right) \Dz$ &  nodes \\
\hline
primal  & & dual                                  \\
\hline
\end{tabular}
}
\caption{Notation for the indices of the nodes and the center points of the
edges, faces and cells in the primal and dual grids where
$-\infty < i,j,k < \infty$.
\label{Dual Primal Grids}}
\end{center}
\end{table}

\begin{table}
\begin{center}
{\renewcommand{\arraystretch}{1.4}
\begin{tabular}{|c|c||c|c|}
\hline
units   & primal                         & dual                                 &units  \\
\hline
$1$     &$s_{i,j,k}$                     &$d^\star_{i,j,k}$                     &$1/d^3$\\
\hline
        &$tx_{i+\half, j, k}$             &$nx^\star_{i+\half, j, k}$             &       \\
$1/d  $ &$ty_{i, j+\half, k}$             &$ny^\star_{i, j+\half, k}$             &$1/d^2$\\
        &$tz_{i, j, k+\half}$             &$nz^\star_{i, j, k+\half}$             &       \\
\hline
        &$nx_{i, j+\half, k+\half}$       &$tx^\star_{i, j+\half, k+\half}$       &       \\
$1/d^2$ &$nz_{i+\half, j, k+\half}$       &$ty^\star_{i+\half, j, k+\half}$       &$1/d$  \\
        &$nz_{i+\half, j+\half, k}$       &$tz^\star_{i+\half, j+\half, k}$       &       \\
\hline
$1/d^3$ &$d_{i+\half, j+\half, k+\half}$ &$s^\star_{i+\half, j+\half, k+\half}$ &$1$    \\
\hline
units   & primal                         & dual                                 &units   \\
\hline
\end{tabular}
}
\caption{Notation for the primal and dual, scalar and vector fields where
$-\infty < i,j,k < \infty$.
\label{scalar and vector fields}}
\end{center}
\end{table}

Mimetic discretizations use primal and dual spatial grids as shown in
Figure \ref{Figure-Dual-Grid} and the notation for the nodes, edges, faces
and cells of the grid are given in Table \ref{Dual Primal Grids} while the
notation for scalar and vector fields are given in Table \ref{scalar and
vector fields}.  It is important that the components of vector fields are
not located at the same points in the grid.
All scalar and vector fields are defined on all of space are assumed
to converge to zero far from the origin. The discussion for boundary
value problems is more complex and will be started in Section
\ref{Implementaion 2D}.

There two types of scalar fields and also two type of vector fields on
both the primal and dual grids.  On the primal grid there are scalar
fields $s$ that do not have a spatial dimension, vector fields $\vec t$
(for tangent) that have spatial dimension $1/d$, vector fields $\vec n$
(for normal) that have units $1/d^2$, and scalar fields with spatial
dimension $1/d^3$ (as in densities) while the dual grid has the same
types of fields labeled with a superscript star as in $s^\star$.
Note that at each point in the grid there is a value from both
the primal and dual fields.  These fields are different because their
spatial dimensions are not the same.

Historically, this type of discretization appeared in the Yee grid for
Maxwell equations \cite{Yee1966} which will be discussed in Section
\ref{E and M section}.

\subsection{The Discrete Double Exact Sequences}

This section describes the discrete double exact sequences shown in
Figure \ref{Discrete-Exact-Sequences}.  This begins with a description
of the discrete difference operators gradient, curl and divergence on
the primal and dual grids. Next the star or multiplication operators that
describe material properties are discretized.

\begin{figure}
\begin{equation*}
  \begin{CD}
                   \SpaceN
        @>\GRAD >> \SpaceE
        @>\CURL >> \SpaceF
        @>\DIVG >> \SpaceC
\\
@V{a}VV @V{\bf A}VV @A{\bf B}AA @A{b}AA @. \\
                       \StarSpaceC
        @<\DIVGstar << \StarSpaceF
        @<\CURLstar << \StarSpaceE
        @<\GRADstar << \StarSpaceN
\\
  \end{CD}
\end{equation*}
\caption{Discrete Exact Sequences \label{Discrete-Exact-Sequences} }
\end{figure}

\subsubsection{Difference Operators}

The discrete gradient $\GRAD$, curl or rotation $\CURL$ and divergence $\DIVG$
are difference operators on a scalar or vector fields.  The formulas for the
dual grid are obtained by making the changes $i \rightarrow i+1/2$,
$j \rightarrow j+1/2$ and $k \rightarrow k+1/2$.

\noindent {\bf The Gradient:}
If $s \in \SpaceN $ is a discrete scalar field, then its gradient
$\GRAD s = (\GRAD sx, \GRAD sy, \GRAD sz) \in \SpaceE$ is an edge
vector field:
\begin{align}
\GRAD sx_{i+\half, j, k} \equiv &
	\frac{ s_{i+1,j,k}- s_{i,j,k} }{\dx} \,; \nonumber \\
\GRAD sy_{i, j+\half, k} \equiv &
	\frac{s_{i,j+1,k}-s_{i,j,k}}{\dy} \,; \\
\GRAD sz_{i, j, k+\half} \equiv &
	\frac{s_{i,j,k+1}-s_{i,j,k}}{\dz} \,. \nonumber
\end{align}

\noindent {\bf The Curl:}
If $\vec{t} = (tx, ty, tz) \in \SpaceE$ is a discrete edge vector field,
then its curl $\CURL {\vec t} \in \SpaceF$ is a discrete face vector field:
\begin{align}
\CURL tx_{i, j+\half, k+\half} & \equiv
	  \frac{tz_{i, j+1, k+\half}- tz_{i, j, k+\half}}{\dy}
	- \frac{ty_{i, j+\half, k+1}- ty_{i, j+\half, k}}{\dz}
		\,; \nonumber \\
\CURL ty_{i+\half, j, k+\half} & \equiv
	  \frac{tx_{i+\half, j, k+1}- tx_{i+\half, j, k}}{\dz}
	- \frac{tz_{i+1, j, k+\half}- tz_{i, j, k+\half}}{\dx}
		\,; \\
\CURL tz_{i+\half, j+\half, k} & \equiv
	  \frac{ty_{i+1, j+\half, k}- ty_{i, j+\half, k}}{\dx}
	- \frac{tx_{i+\half, j+1, k}- tx_{i+\half, j, k}}{\dy}
		\,. \nonumber 
\end{align}

\noindent {\bf The Divergence:}
If ${\vec n} = (nx, ny, nz ) \in \SpaceF$ is a discrete face vector field,
then its divergence $\DIVG {\vec n} \in \SpaceC $ is a cell scalar field:
\begin{align}
\DIVG {\vec n}_{i+\half, j+\half, k+\half} & \equiv  
	      \frac{nx_{i+1,j+\half,k+\half}- nx_{i,j+\half,k+\half}}{\dx}
			\nonumber \\
	& +  \frac{ny_{i+\half,j+1,k+\half}- ny_{i+\half,j,k+\half}}{\dy}
			\\
	& +  \frac{nz_{i+\half,j+\half,k+1}- nz_{i+\half,j+\half,k}}{\dz}
			\,. \nonumber 
\end{align}

\noindent {\bf The Star Gradient:}
If $s^\star \in \StarSpaceN $ is a discrete star scalar field then
its star gradient $\GRADstar \, s^\star \in \StarSpaceE$ is a star edge
vector field:
\begin{align}
\GRADstar s^\star x_{i, j+\half, k+\half} & \equiv
\frac{s^\star_{i+\half,j+\half,k+\half}-s^\star_{i-\half,j+\half,k+\half} }{\Delta x}; \nonumber \\
\GRADstar s^\star y_{i+\half, j, k+\half} & \equiv
\frac{s^\star_{i+\half,j+\half,k+\half}-s^\star_{i+\half,j-\half,k+\half} }{\Delta x}; \\
\GRADstar s^\star z_{i+\half, j+\half, k} & \equiv
\frac{s^\star_{i+\half,j+\half,k+\half}-s^\star_{i+\half,j+\half,k-\half} }{\Delta x}; \nonumber 
\end{align}

\noindent {\bf The Star Curl:}
If $\vec{t}^\star = (tx^\star, ty^\star, tz^\star ) \in \StarSpaceC$
is a discrete star edge vector field then its curl
$\CURLstar \, \vec{t}^\star \in \StarSpaceF$ is a discrete star face
vector field:
\begin{align}
\CURLstar t^\star x_{i+\half, j, k} & \equiv
          \frac{tz^\star_{i+\half, j+\half, k}-
        tz^\star_{i+\half, j-\half, k}}{\dy}
        - \frac{ty^\star_{i+\half, j, k+\half}-
        ty^\star_{i+\half, j, k-\half}}{\dz}
                \,; \nonumber \\
\CURLstar t^\star y_{i, j+\half, k} & \equiv
          \frac{tx^\star_{i, j+\half, k+\half}-
        tz^\star_{i, j+\half, k-\half}}{\dz}
        - \frac{tz^\star_{i+\half, j+\half, k}-
        tz^\star_{i-\half, j+\half, k}}{\dx}
                \,; \\
\CURLstar t^\star z_{i, j, k+\half} & \equiv
          \frac{ty^\star_{i+\half, j, k+\half}-
        ty^\star_{i-\half, j, k+\half}}{\dx}
        - \frac{tx^\star_{i, j+\half, k+\half}-
        tx^\star_{i, j-\half, k+\half}}{\dy}
                \,. \nonumber
\end{align}

\noindent {\bf The Star Divergence:}
If $\vec{n}^\star = (nx^\star, ny^\star, nz^\star) \in \StarSpaceF$
is a discrete star face vector field then it divergence
$\DIVGstar \vec{n}^\star \in \StarSpaceC$ is a discrete star cell field.
In terms of components
\begin{align}
\DIVGstar {\vec n^\star}_{i, j, k}
                & \equiv
              \frac{nx^\star_{i+\half,j,k}-
        n^\star x_{i-\half,j,k}}{\dx} \nonumber \\
        & +  \frac{ny^\star_{i,j+\half,k}-
        ny^\star_{i,j-\half,k}}{\dy} \\
        & +  \frac{nz^\star_{i,j,k+\half}-
        nz^\star_{i,j,k-\half}}{\dz} \,.  \nonumber
\end{align}

The second order accuracy of the difference operators is confirmed in
{\tt TestAccuracy3.m}.

\subsubsection{Mimetic Properties of Difference Operators}

If $c$ is a constant scalar field then a direct computation 
\cite{RobidouxSteinberg2011} shows that:
\begin{equation}
\GRAD c \equiv 0 \,,\quad
\CURL \GRAD \equiv 0 \,,\quad
\DIVG \CURL \equiv 0 \,,\quad
\GRADstar c \equiv 0 \,,\quad
\CURLstar \GRADstar \equiv 0 \,,\quad
\DIVGstar \CURLstar \equiv 0 \,.
\end{equation}
These relationships are confirmed in {\tt TestZero3.m}. These properties
are summarized by saying that the discretization is exact or that the
sequences in Figure \ref{Discrete-Exact-Sequences} are exact, see
\cite{RobidouxSteinberg2011} for a precise definition of exact and
a proof that the diagram is exact.

\subsubsection{Discrete Star or Multiplication Operators}

The star operators are multiplication operators that model the material
properties and are given by two scalar functions
$a = a(x,y,z)$ and $b=b(x,y,z)$
and two $3 \times 3$ matrix functions
${\bf A} = {\bf A}(x,y,z)$ and ${\bf B} = {\bf B}(x,y,z)$
that are symmetric and positive. The spatial dimensions of $a$ and $b$
must be $1/s^3$ while for ${\bf A}$ and ${\bf B}$ must be $1/d$.

If $a_{i,j,k} = a(i\,\Dx,j\,\Dy,k\,\Dz )$ and if $s \in \SpaceN$ and
$ d^\star = a \, s \in \StarSpaceC$ then 
\[
 d^\star_{i,j,k} = a_{i,j,k} \, s_{i,j,k} \,.
\]
If $b_{i,j,k} = b(i\,\Dx,j\,\Dy,k\,\Dz )$ and if $s^\star \in \StarSpaceN$ and
$ d = b \, s^\star \in \StarSpaceC$ then 
\[
 d_{i,j,k} = b_{i,j,k} \, s^\star_{i,j,k} \,.
\]
The assumption that $a$ and $b$ are not zero implies that these star operators
are invertible.

Let 
\begin{equation}
{\bf A} = 
\left[
\begin{matrix}
Axx & Axy & Axy \\
Ayx & Ayy & Ayz \\
Azx & Azy & Azz
\end{matrix}
\right] \,,
\end{equation}
where ${\bf A}$ is symmetric, positive definite, and the entries in ${\bf A}$
are functions of $(x,y,z)$.  The discretized ${\bf A}$ will not be symmetric
but will be nearly symmetric.

For ${\vec t} \in \SpaceE$, computing
${\vec n}^\star = {\bf A} \vec{t} \in \StarSpaceN$
requires averaging of the off diagonal terms in ${\bf A}$ to maintain
second order accuracy.  First set
\begin{align}
nx^\star_{i+\half,j,k} = &
   Axx_{i+\half,j,k} \,      tx_{i+\half,j,k} +
   Axy_{i+\half,j,k} \, \ovl{ty}_{i+\half,j,k} +
   Axz_{i+\half,j,k} \, \ovl{tz}_{i+\half,j,k} \nonumber \\
ny^\star_{i,j+\half,k} = & 
   Ayx_{i,j+\half,k} \, \ovl{tx}_{i,j+\half,k} +
   Ayy_{i,j+\half,k} \,      ty_{i,j+\half,k} +
   Ayz_{i,j+\half,k} \, \ovl{tz}_{i,j+\half,k} \label{A times t} \\
nz^\star_{i,j,k+\half} = & 
   Ayx_{i,j,k+\half} \, \ovl{tx}_{i,j,k+\half} +
   Ayy_{i,j,k+\half} \, \ovl{ty}_{i,j,k+\half} +
   Ayz_{i,j,k+\half} \,      tz_{i,j,k+\half} \,.  \nonumber
\end{align}
The average values are given by
\begin{align*}
\ovl{ty}_{i+\half,j,k} & = 
    \frac{ ty_{i,  j-\half,k} + ty_{i,  j+\half,k} + 
           ty_{i+1,j-\half,k} + ty_{i+1,j+\half,k}}{4} \,, \\
\ovl{tz}_{i+\half,j,k} & =
    \frac{ tz_{i,  j,k-\half} + tz_{i,  j,k+\half} +
           tz_{i+1,j,k-\half} + tz_{i+1,j,k+\half}}{4} \,,\\
\ovl{tx}_{i,j+\half,k} & = 
    \frac{tx_{i,  j-\half,k} +  tx_{i,  j+\half,k} + 
          tx_{i+1,j-\half,k} +  tx_{i+1,j+\half,k}}{4} \,, \\
\ovl{tz}_{i,j+\half,k} & = 
    \frac{tz_{i,  j,k-\half} +  tz_{i,  j,k+\half} + 
          tz_{i,j+1,k-\half} +  tz_{i,j+1,k+\half}}{4} \,, \\
\ovl{tx}_{i,j,k+\half} & =
    \frac{tx_{i-\half,j,k} + tx_{i+\half,j,k} +
          tx_{i-\half,j,k+1} + tx_{i+\half,j,k+1}}{4} \,, \\
\ovl{ty}_{i,j,k+\half} & =
    \frac{ty_{i,j-\half,k} + ty_{i,j+\half,k} +
          ty_{i,j-\half,k+1} + ty_{i,j+\half,k+1}}{4} \,. \\
\end{align*}

Note that if ${\bf A}$ is diagonal then multiplication by ${\bf A}$ is
simply multiplication by the diagonal entries of ${\bf A}$ and no averaging
is required.  There are similar formulas for multiplication by
${\bf B}$, ${\bf A}^{-1}$ and ${\bf B}^{-1}$.
Many of the second order operators in Table \ref{General Fundamental Operators}
have a multiplication by the inverse of ${\bf A}$ and/or ${\bf B}$.
It is assumed that the matrix operators are given by formulas, and formulas
for the inverse operators can be found and so that the formulas above
can be used to multiply by the inverse matrices.

\subsection{Discrete Inner Products \label{Discrete Inner Products}}

To study conserved quantities an inner product is needed for each of the eight
linear spaces in the dual exact sequences
Figure \eqref{Discrete-Exact-Sequences}. The inner products will be defined
in terms of four bilinear forms as in \ref{bilinear form}.
As in the continuum, an important property of the inner products is that
they need to be symmetric, positive definite and importantly dimensionless.

Four bilinear forms will be needed.
Set $\Delta V = \Delta x \, \Delta y \, \Delta z$.

\noindent
If $s \in \SpaceN$ and $d^\star \in \StarSpaceC$ then
\[
\lbilin s , d \rbilin  = \sum s_{i,j,k} \, d^\star_{i,j,k} \,
\Delta V \,.
\]
\noindent
If $\vec{t} \in \SpaceE$ and $\vec{n}^\star \in \StarSpaceF$ then
\[
\lbilin \vec{t} , \vec{n}^\star \rbilin  = 
\left(
\sum tx_{i+\half,j,k} \, nx^\star_{i+\half,j,k} +
\sum ty_{i,j+\half,k} \, ny^\star_{i,j+\half,k} +
\sum tz_{i,j,k+\half} \, nz^\star_{i,j,k+\half}
\right)
\Delta V \,.
\]

\noindent
If $\vec{n} \in \SpaceE$ and $\vec{t}^\star \in \StarSpaceF$ then
\[
\lbilin \vec{n} , \vec{t}^\star \rbilin  = 
\left(
\sum nx_{i,      j+\half,k+\half} \, tx^\star_{i,      j+\half,k+\half} +
\sum ny_{i+\half,j,      k+\half} \, ty^\star_{i+\half,j,      k+\half} +
\sum nz_{i+\half,j+\half,k}  \,      tz^\star_{i+\half,j+\half,k}
\right)
\Delta V \,.
\]
\noindent
If $g \in \SpaceC$ and $ f^\star \in \StarSpaceN$ then
\[
\lbilin g , f^\star \rbilin = 
\sum g_{i+\half,j+\half,k+\half} \, f^\star_{i+\half,j+\half,k+\half}
\Delta V \,.
\]

The eight inner products are given by are given by the bilinear forms.

\noindent
If $s1, s2 \in \SpaceN$ then
$ \langle s1, s2 \rangle_\Nodes = \lbilin a \, s1 , s2 \rbilin $.

\noindent
If $s1^\star, s2^\star \in \StarSpaceN$ then
\[
\langle s1^\star, s2^\star \rangle_\StarNodes = \sum_{i,j,k}
 b_{i+\half,j+\half,k+\half} \,
s1^\star_{i+\half,j+\half,k+\half} \,
s2^\star_{i+\half,j+\half,k+\half} \dx \dy \dz \,.
\]
If $\vec{t1}, \vec{t2} \in \SpaceE$ then
\begin{align*}
\langle \vec{t1}, \vec{t2} \rangle_\Edges = \sum_{i,j,k} (
&
({\bf A}\,\vec{t1})_{i+\half,j,k} \, \vec{t2}_{i+\half,j,k} + \\
&({\bf A}\,\vec{t1})_{i,j+\half,k} \, \vec{t2}_{i,j+\half,k} + \\
&({\bf A}\,\vec{t1})_{i,j,k+\half} \, \vec{t2}_{i,j,k+\half}
) \dx \dy \dz
\end{align*}
If $\vec{t1}^\star, \vec{t2}^\star \in \StarSpaceE$ then
\begin{align*}
\langle \vec{t1}^\star, \vec{t2}^\star \rangle_\StarEdges = 
\sum_{i,j,k} ( 
& ({\bf B} \, \vec{t1}^\star)_{i, j+\half, k+\half} \, \vec{t2}^\star_{i, j+\half, k+\half} + \\
& ({\bf B} \, \vec{t1}^\star)_{i,+\half j, k+\half} \, \vec{t2}^\star_{i+\half,j,k+\half} + \\
& ({\bf B} \, \vec{t1}^\star)_{i+\half, j+\half, k} \, \vec{t2}^\star_{i+\half, j+\half, k} ) \dx \dy \dz
\end{align*}
If $\vec{n1}, \vec{n2} \in \SpaceF$ then
\begin{align*}
\langle \vec{n1}, \vec{n2} \rangle_\Faces = 
\sum_{i,j,k} ( 
& ({\bf B}^{-1} \, n1)_{i, j+\half, k+\half} \, n2{i, j+\half, k+\half} + \\
& ({\bf B}^{-1} \, n1)_{i,+\half j, k+\half} \, n2{i+\half,j,k+\half} + \\
& ({\bf B}^{-1} \, n1)_{i+\half, j+\half, k} \, n2{i+\half, j+\half, k} ) \dx \dy \dz
\end{align*}
If $\vec{n1}^\star, \vec{n2}^\star \in \StarSpaceF$ then
\begin{align*}
\langle \vec{n1}^\star, \vec{n2}^\star \rangle_\StarFaces = 
\sum_{i,j,k} (
&({\bf A}^{-1}\, n1)^\star_{i+\half, j, k}\, n2^\star_{i+\half, j, k} + \\
&({\bf A}^{-1}\, n1)^\star_{i, j+\half, k}\, n2^\star_{i, j+\half, k} + \\
&({\bf A}^{-1}\, n1)^\star_{i, j, k+\half}\, n2^\star_{i, j, k+\half}
) \dx \dy \dz
\end{align*}
If $d1, d2 \in \SpaceC$ then
\[
\langle d1, d2 \rangle_\Nodes = \sum_{i,j,k}
 b^{-1}_{i+\half,j+\half,k+\half} \,
d1_{i+\half,j+\half,k+\half} \,
d2_{i+\half,j+\half,k+\half} \dx \dy \dz \,.
\]
If $d1^\star, d2^\star \in \StarSpaceC$ then
\[
\langle d1^\star, d2^\star \rangle_\StarCells = 
 \sum_{i,j,k} a^{-1}_{i,j,k} \, d1^\star_{i,j,k} \, d2^\star_{i,j,k} \dx \dy \dz \,.
\]

\subsection{Adjoint Operators}

The derivation of the adjoints of the discrete operators in the discrete
exact sequences shown in Figure \ref{Discrete-Exact-Sequences} are 
similar to the derivation for the continuum the adjoint operators
defined in Section \ref{Continuum Adjoint Operators}. Again note that the
discrete operators are not mapping of  a space into itself. The adjoints
can easily be derived by diagram chasing using \ref{Discrete-Exact-Sequences}.

\begin{align}
\GRAD^*& = - a^{-1} \, \DIVGstar \, {\bf A} &
\GRADstar^* & = - b^{-1} \, \DIVG \, {\bf B}
\,,\nonumber\\
\CURL^* & = + {\bf A}^{-1} \, \CURLstar \, {\bf B}^{-1} &
\CURLstar^* & = + {\bf B}^{-1} \, \CURL \, {\bf A}^{-1} 
\,,\nonumber\\
\DIVG^* & = - {\bf B} \GRADstar  b^{-1} &
\DIVGstar^* & = - {\bf A} \GRAD  a^{-1}
\,,\nonumber\\
{\bf A}^* & = {\bf A}^{-1} &
{\bf B}^* & = {\bf B}^{-1}
\,, \label{DiscreteAdjoints} \\
a^* & = a^{-1} &
b^* & = b^{-1}
\nonumber\,.
\end{align}

The proofs of the adjoint formulas rely of summation by parts which is
illustrate by computing
\begin{align*}
\sum_{i,j,k}
\left(s_{i+1,j,k}-s_{i,j,k}\right) \, n^\star_{i+\half, j, k}
& = \sum_{i,j,k} s_{i+1,j,k} \, n^\star_{i+\half, j, k}
   -\sum_{i,j,k} s_{i,j,k} \, n^\star_{i+\half, j, k} \\
& = \sum_{i,j,k} s_{i,j,k} \, n^\star_{i-\half, j, k}
    -\sum_{i,j,k} s_{i,j,k} \, n^\star_{i+\half, j, k} \\
& = \sum_{i,j,k}
s_{i,j,k} \, \left( n^\star_{i-\half, j, k} - n^\star_{i+\half, j, k} \right) \\
& = - \sum_{i,j,k}
s_{i,j,k} \, \left( n^\star_{i+\half, j, k} - n^\star_{i-\half, j, k} \right) 
\end{align*} 
The remaining proofs are straight forward.
Others have used summation by parts to obtain mimetic like discretizations
\cite{WangKreiss2017,FernandezHZ14,NordstromL13}.

When working with systems of first order wave equation and second order
wave equation additional adjoints will be needed:
\[
({\bf A} \GRAD)^*
= \GRAD^* {\bf A}
= a^{-1} \DIVG^* {\bf A} {\bf A}^-1
= a^{-1} \DIVG^* \,.
\]

\subsection{Positive and Negative Discrete Operators}

\noindent If $s \in \in \SpaceN$ and $\vec{n}^\star \in \StarSpaceF$ then
\begin{equation} \label{Adjoint 1}
\langle {\bf A} \, \GRAD \, s , \vec{n}^\star \rangle_\StarFaces =
- \langle s  , \frac{1}{a} \, \DIVGstar \, \vec{n}^\star \rangle_\Nodes
\end{equation}
If $\vec{t} \in \SpaceE$ and $\vec{t}^\star \in \StarSpaceE$ then
\begin{equation} \label{Adjoint 2}
\langle {\bf B}^{-1} \, \CURL \vec{t} , \vec{t}^\star \rangle_\StarEdges
= - \langle \vec{t} \,, {\bf A}^{-1} \CURLstar \vec{t}^\star \rangle_\Edges
\end{equation}
If $\vec{n} \in \SpaceF$ and $\vec{t}^\star \in \StarSpaceN$ then
\begin{equation} \label{Adjoint 3}
\langle b^{-1} \, \DIVG \vec{n} , s^\star \rangle_\StarNodes =
- \langle \vec{n} \,, {\bf B} \GRADstar \vec{s}^\star \rangle_\Faces
\end{equation}

For an example, consider \eqref{Adjoint 1}:
\begin{align*}
& \langle {\bf A} \, \GRAD \, s , \vec{n}^\star \rangle_\StarFaces = \\
& \sum_{i,j,k} \left( 
 \frac{s_{i+1,j,k}-s_{i,j,k}}{\dx} \, n^\star_{i+\half, j, k} +  
 \frac{s_{i,j+1,k}-s_{i,j,k}}{\dy} \, n^\star_{i, j+\half, k} + 
 \frac{s_{i,j,k+1}-s_{i,j,k}}{\dz} \, n^\star_{i, j, k+\half}
\right) \dx \dy \dz
\end{align*}
\begin{align*}
& \langle s  , \frac{1}{a} \, \DIVGstar \, \vec{n}^\star \rangle_\Nodes = \\
& \sum_{i,j,k} s_{i,j,k} \,
\left( \frac{n^\star_{i+\half,j,k}- n^\star_{i-\half,j,k}}{\dx}+
        \frac{n^\star_{i,j+\half,k}- n^\star_{i,j-\half,k}}{\dy}+
        \frac{n^\star_{i,j,k+\half}- n^\star_{i,j,k-\half}}{\dz} \right)
\dx \dy \dz \,.
\end{align*}
Three summation by parts will prove \eqref{Adjoint 1}. The remaining
formulas \eqref{Adjoint 2} and \eqref{Adjoint 3} can be proved in the same way.

\newpage \clearpage
\setcounter{equation}{0}
\section{Discretizing Wave Equations in 3D
\label{Discretize Wave}}

The previous results will be used to discretize the scalar
wave equation and Maxwell's wave equation in three dimensions.
Here the simulation region is all space, boundary conditions for
bound domain will be discussed later.

\subsection{The Scalar Wave Equation \label{Scalar Wave 3d}}

The second order scalar wave equation \eqref{Second Order 1} can be written
as as first order system as in \eqref{System 1}. But here the notation will
be changed to match that in Section \ref{Mimetic Discretizations}:
\begin{equation*}
\frac{\partial s}{\partial t} = a^{-1} \divg \, \vec{v} \,,\quad
\frac{\partial \vec{v}}{\partial t} = {\bf A} \grad  s  \,,
\end{equation*}
with $s \in H_P$ and $\vec{v} \in H_S$.
This system will be discretized using the operators described in Section
\ref{Mimetic Discretizations} so now let $s \in \SpaceN$ and
$v \in \StarSpaceF$ and then the leapfrog discretization is
\begin{equation}
\frac{s^{n+1} - s^{n}}{\dt} =  a^{-1}\, \DIVGstar v^{n+\half} \,,\quad
\frac{v^{n+\half} - v^{n-\half}}{\dt} = {\bf A} \GRAD s^n \,.
\label{Discrete 3D Wave}
\end{equation}
If $s^0$ and $v^{\half}$ are given then the leapfrog scheme for $n \geq 0$ is
\[
s^{n+1} = s^{n} + \dt \, a^{-1}\, \DIVGstar v^{n+\half} \,,\quad
v^{n+\thalf} = v^{n+\half} + \dt \, {\bf A} \GRAD s^{n+1} \,.
\]
This gives a discretization of a second order scalar wave equation
\eqref{Second Order 1} as
\begin{equation}
\frac{s^{n+1} - 2\, s^{n} + s^{n-1}}{\dt^2}
= a^{-1}\, \DIVGstar {\bf A} \GRAD s^n \,, 
\end{equation}
and as discussed in Section \ref{Wave Equations With Variable Materials} on continuum
wave equations.

This is also a discretization of the vector wave equation
\begin{equation}
\frac{v^{n+\thalf} - 2\, v^{n+\half} + v^{n-\half}}{\dt^2}
= {\bf A} \GRAD a^{-1}\, \DIVGstar v^{n+\half} \,. 
\end{equation}

The results in Section \ref{ODEs} give two conserved quantities
for the discretization. To see this using \eqref{Conserved A half}
and \eqref{Conserved A* full}  set
$A$ to $a^{-1} \, \DIVGstar$,
$A^*$ to $-{\bf A} \GRAD$,
$g^{n+\half}$ to $v^{n+\half}$ and
$f^n$ to $s^n$
to get the conserved conserved quantities
\begin{equation}
C^n =
 \norm{s^{n}}_\Nodes^2
 + \norm{\frac{v^{n+1/2} + v^{n-1/2}}{2}}_\StarFaces^2
 - \frac{\Delta t^2}{4} \norm{ {\bf A} \GRAD \, s^{n} }_\StarFaces^2 \,,
\end{equation}
\begin{equation}
C^{n+1/2} =
\norm{v^{n+1/2}}_\StarFaces^2
+ \norm{\frac{s^{n+1} + s^{n}}{2}}_\Nodes^2
- \frac{\Delta t^2}{4} \norm{ a^{-1} \, \DIVGstar \, v^{n+1/2} }_\Nodes^2 \,.
\end{equation}
These formulas agree with those derived in detail in Appendix
\ref{Appendix Details}.

The programs {\tt ScalarWave.m} and {\tt ScalarWaveStar.m} were used
to show that the energies $C_n$ and $C_{n+\half}$ are constant to
less than 1 part in $10^{15}$ for the discretization described above
and the one that changes $\GRAD$ to $\GRADstar$ and $\DIVGstar$ to
$\DIVG$.

\subsection{Maxwell's Equations \label{E and M section}}

Assume that ${\vec E} \in \SpaceE$ and $\vec{H} \in \StarSpaceE$ so that,
using the notation in the previous sections, the Maxwell system
\ref{Maxwell System} will be discretized as 
\[
\frac{{\vec E}^{n+1} - {\vec E}^n}{\dt} = \epsilon^{-1} \, \CURLstar \vec{H}^{n+\half} \,,\quad
\frac{\vec{H}^{n+\half} - \vec{H}^{n-\half}}{\dt} = - \mu^{-1} \, \CURL {\vec E}^n \,.
\]
where in Exact Sequence diagram \eqref{Discrete-Exact-Sequences}
$A = \epsilon$ and $B = \mu$. Here $\epsilon$ and $\mu$
can be symmetric positive definite matrices.  If
$\vec{E} = (Ex, Ey,Ez)$ and  $\vec{H} = (Hx, Hy,Hz)$
then Table \ref{Dual Primal Grids} shows that $Ex$ and $Hx$ are indexed as
\[
Ex^n_{i+\half,j,k} \,,\quad Hx^{n+\half}_{i,j+\half,k+\half}
\]
just as in Yee's paper \cite{Yee1966}.  If $\vec{E}^0$ and $\vec{H}^\half$
are given then the leapfrog scheme for $n \geq 0$ is
\[
{\vec E}^{n+1} = {\vec E}^n
+ \dt \, \epsilon^{-1} \, \CURLstar \vec{H}^{n+\half} \,,\quad
\vec{H}^{n+3/2} = \vec{H}^{n+\half}
- \dt \, \mu^{-1} \, \CURL {\vec E}^{n+1} \,.
\]

Using a similar argument, it is easy to see that
\[
C_{n+1/2} =
\norm{\frac{\vec{E}^{n+1} + \vec{E}^{n}}{2}}_\Edges^2
+ \norm{ \vec{H}^{n+1/2}}_\StarEdges^2 
- \frac{\Delta t^2}{4}
\norm{ \epsilon^{-1} \, \CURLstar \, \vec{H}^{n+1/2} }_\Edges^2 
\]
is a conserved quantity and that
\[
C_{n+1/2} \geq
\norm{\frac{\vec{E}^{n+1} + \vec{E}^{n}}{2}}_\Edges^2 +
\left( 1 - \frac{\Delta t^2}{4} \norm{\epsilon^{-1} \, \CURLstar}^2 \right)
\norm{ \vec{H}^{n+1/2}}_\StarEdges^2  \,.
\]
So $C_{n+1/2} \geq 0$ for $\Delta t$ sufficiently small provided
$\norm{\epsilon^{-1} \, \CURLstar}$ is finite.

Also
\[
C_n =  \norm{\vec{E}^{n}}_\Edges^2 
 -\frac{\Delta t^2}{4} \norm{ \mu^{-1} \CURL \, \vec{E}^{n} }_\StarFaces^2 
 +\norm{\frac{\vec{H}^{n+1/2} + \vec{H}^{n-1/2}}{2}}_\StarFaces^2 \,.
\]
is a conserved quantity and 
\[
\norm{C_n} \ge	
\left(1 - \frac{\Delta t^2}{4}
\norm{\mu^{-1} \CURL}^2\right)
\norm{\vec{E}^n}_\Edges^2
+ \norm{\frac{\vec{H}^{n+1/2} + \vec{H}^{n-1/2}}{2}}_\StarFaces^2 \,.
\]
So $\norm{C_n}$ is positive for sufficiently small $\Delta t$ if
$\norm{ \mu^{-1} \CURL \, \vec{E}^{n} }$ is finite.
Also, these formulas agree with those derived in detail in Appendix
\ref{Appendix Details}.

The codes {\tt Maxwell.m} and {\tt MaxwellStar.m} confirm that our algorithms
conserve $C_{n+1/2}$ and $C_n$ to two parts in $10^{16}$. Additionally,
the divergence of the curl of the electric and magnetic fields are constant
to one part in $10^{14}$ when there are no sources.

\newpage \clearpage
\setcounter{equation}{0}
\section{Implementation in 2D \label{Implementaion 2D}}

The goal is to illustrate how scalar and vector functions are discretized
in 2D and then show how the gradient and divergence are discretized.
This will provide some intuition on how to discretize functions and
differential operators in 3D.  For 3D the primal and dual grids are shown
in Figure \ref{Figure-Dual-Grid} and described in Table \ref{Dual Primal Grids}.
In some ways the 2D discretization is more complicated than the 3D 
because in 3D cells have edges and faces while in 2D the edge of a
cell can correspond to either an edge or a face in 3D.

The implementation of mimetic finite difference discretizations in a bounded
region can be annoying due to the two staggered grids and indices that contain
a half. Additionally, the components of vector fields are discretized
at different points. There are also problems at the boundary that the
star operators will be used fix.  Initially the star operators will be
trivial, except at the boundary. The simulation region will be the unit square.

For 3D a detailed description of the mimetic discretization is given
in Section \ref{Mimetic Discretizations}. In particular the left most
box in Figure \ref{Discrete-Exact-Sequences} will be used for the 2D
discretization.

To facilitate programming, all indices are integers greater than zero,
for example $i \geq 1$ and $j \geq 1$.

\subsection{The 2D primal and dual grids}

\begin{figure}
\begin{center}
\includegraphics[width=5.00in,trim=100 247 100 250,clip]{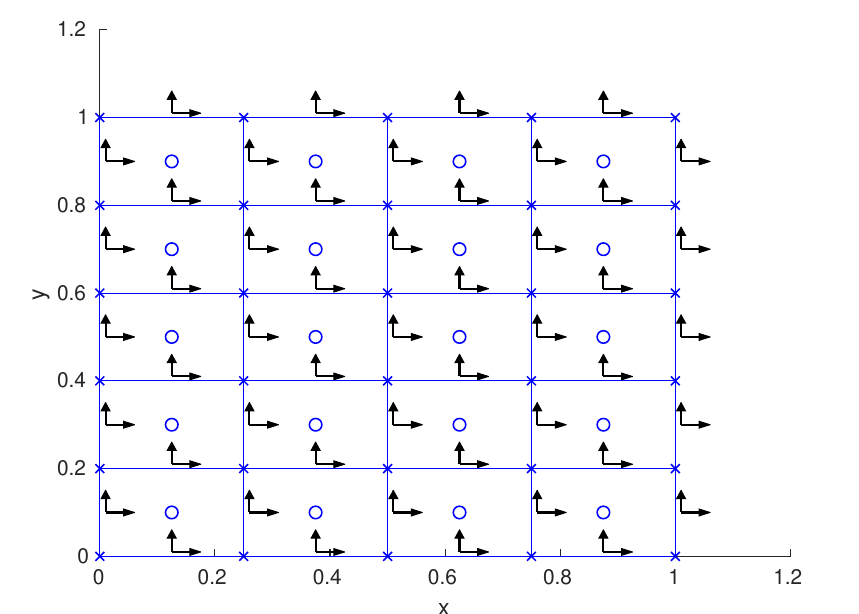}
\end{center}
\caption{2D primal grid tangent and normal vector fields, $Nx = 4$, $Ny = 5$
\label{Tangent Normal Primal}}
\end{figure}
\begin{figure}
\begin{center}
\includegraphics[width=5.00in,trim=100 247 100 250,clip]{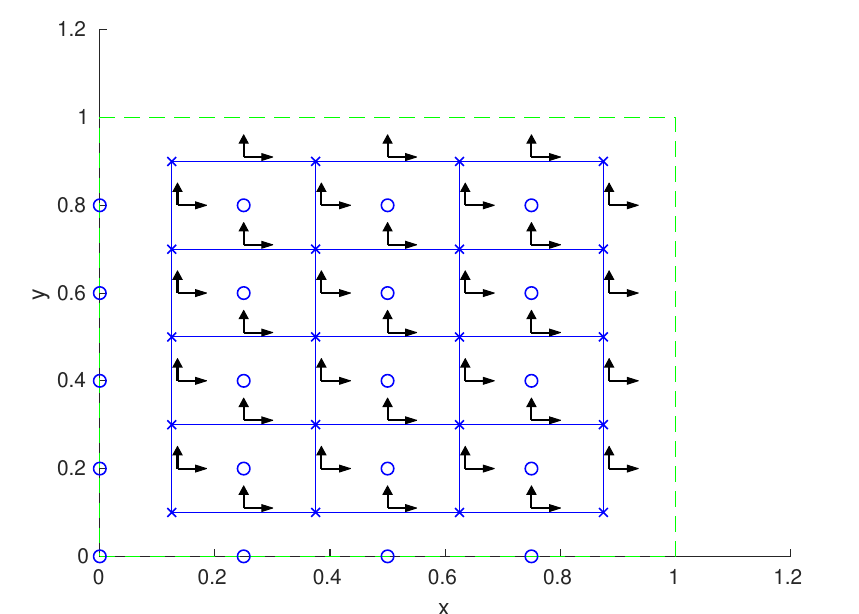}
\end{center}
\caption{2D dual grid tangent and normal vector fields, $Nx = 4$, $Ny = 5$
\label{Tangent Normal Dual}}
\end{figure}

Let $Nx,Ny$ be positive integers and then set $dx = 1/Nx$, $dy = 1/Ny$.
The primal grid nodes (cell corners) are given by
\begin{align*}
(xp(i,j), yp(i,j)) & = ((i-1) \, dx, (j-1) \, dy) \,,\quad
1 \leq i \leq Nx+1 \,,
1 \leq j \leq Ny+1 \,,
\end{align*}
and the dual grid nodes are given by
\begin{align*}
(xd(i,j), yd(i,j)) & = ( (i-1/2) \, dx, (j-1/2) \, dy) \,, \quad
1 \leq i \leq Nx \,,
1 \leq j \leq Ny \,.
\end{align*}
Note that in the primal grid $Nx$ is the number of cells in
$x$ direction and $Ny$ is the number of cells in the $y$ direction.  For
$Nx = 4$ and $Ny = 5$ the positions of scalar and vector function on the
primal and dual grids are illustrated in Figures \ref{Tangent Normal Primal}
and \ref{Tangent Normal Dual}, see {\tt FigurePrimalDual2.m}.
Compare these to Figure \ref{Figure-Dual-Grid} for a 3D grid.

Important points are that a primal grid cell center is given by a dual
grid node and the dual grid centers are given by given by the {\em interior}
primal grid nodes. Additionally the location of dual grid tangent vectors
are the same as the location of the primal grid interior normal vectors
and the position of the dual grid normal vectors are the same as the position
of the interior grid tangent vectors. On the boundary of the primal grid the
positions of points and vectors do not correspond to anything in the dual grid.
As will be seen this is important for representing boundary conditions for
partial differential equations.

\subsection{Discretizing continuum functions on the 2D grids}

\begin{figure}
\begin{center}
\begin{tabular}{c}
\includegraphics[width=6.00in,height=8.00in,trim=100 240 110 215]{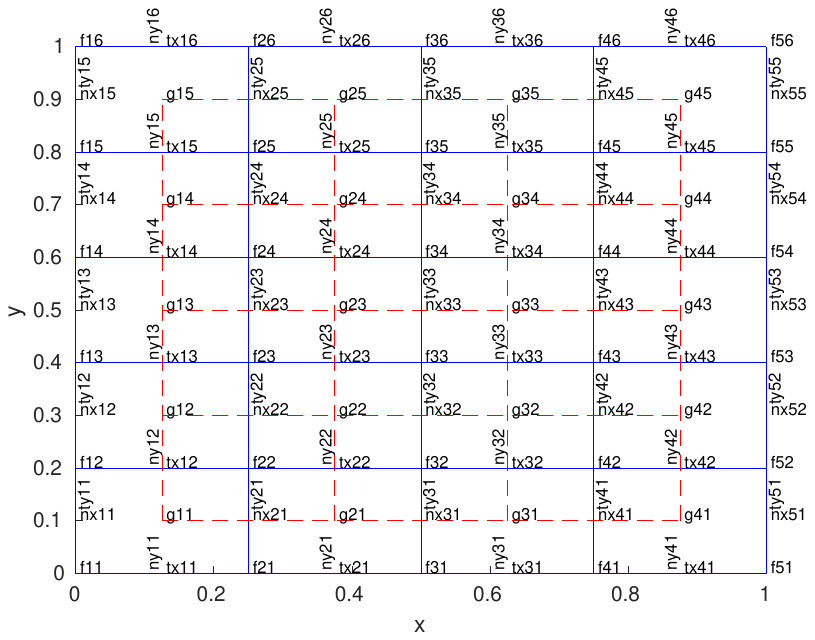}
\end{tabular}
\end{center}
\caption{Scalar and vector fields on the primal grid.
\label{Scalar Vector Primal}}
\end{figure}

\begin{figure}
\begin{center}
\begin{tabular}{c}
\includegraphics[width=6.00in,height=8.00in,trim=100 240 110 215]{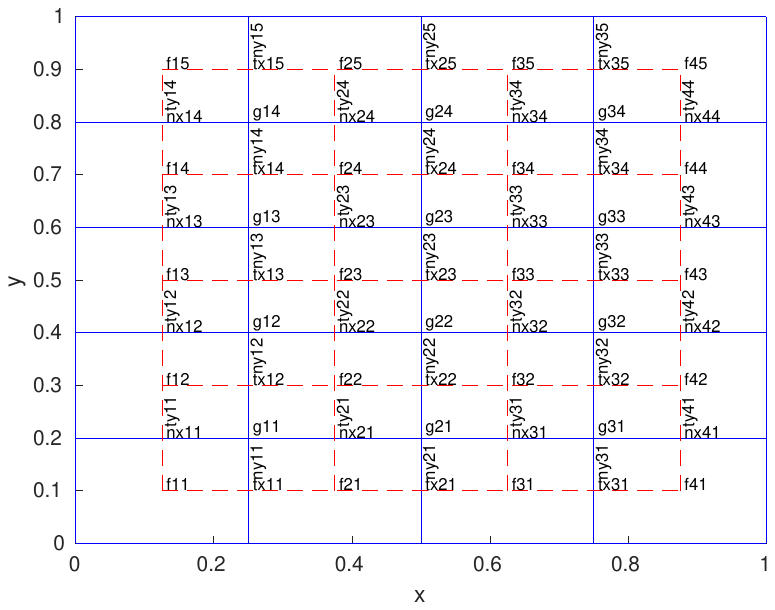}
\end{tabular}
\end{center}
\caption{Scalar and vector fields on the dual grid,
\label{Scalar Vector Dual}}
\end{figure}

Scalar functions will be discretized at either the nodes (corners) or
cell centers of the primal and dual grids while vector fields will be
discretized at the centers of the  of the edges of the cells as illustrated
in Figures \ref{Tangent Normal Primal} and \ref{Tangent Normal Dual} for grids
in a region that is a unit square (see {\tt FigureDetailsPrimalDual2.m}).
There are a total of eight types of discretized functions.
In the text, function names on the primal grid have a $p$ appended as in
$fp$ while functions on the dual grid have a $d$ appended as in $fd$.
In the figures which grid the functions are on is clear.

For the primal grid there are four cases.  A scalar function $f(x,y)$
with spatial weight $1$ is discretized at primal cell nodes:
\begin{equation}
fp(i,j) = f(xp(i),yp(j)) \,,\quad i \leq Nx+1 \,,\quad j \leq Ny+1 \,.
\end{equation}
A scalar function $g(x,y) $ with spatial weights $1/d^3$ is discretized at
primal grid cell centers:
\begin{equation}
gp(i,j) = g(xd(i),yd(j)) \,,\quad i \leq Nx \,,\quad j \leq Ny \,.
\end{equation}
A vector function $(tx(x,y),ty(x,y)$ with spatial weight $1/d$ (tangent on
primal grid) is discretized at cell edge centers:
\begin{align}
txp(i,j) = & tx(xd(i),yp(j)) \,,\quad i \leq Nx   \,,\quad j \leq Ny+1 \,;\\
typ(i,j) = & ty(xp(i),yd(j)) \,,\quad i \leq Nx+1 \,,\quad j \leq Ny \,. 
\end{align}
A vector function $(nx(x,y),ny(x,y)$ with spatial weight $1/d^2$ (normal on
primal grid) is also discretized at cell edge centers:
\begin{align}
nxp(i,j) = & nx(xp(i),yd(i)) \,,\quad i \leq Nx+1 \,,\quad j \leq Ny \,;\\
nyp(i,j) = & ny(xd(i),yp(j)) \,,\quad i \leq Nx   \,,\quad j \leq Ny+1 \,. 
\end{align}

For the dual grid there are also four cases.  A scalar function $f(x,y)$ with
spatial weight $1$ is discretized at the dual cell nodes:
\begin{equation}
fd(i,j) = f(xd(i),yd(j)) \,,\quad i \leq Nx \,,\quad j \leq Ny \,.
\end{equation}
A scalar function $g(x,y) $ with spatial weights $1/d^3$ is discretized at
dual grid cell centers nodes:
\begin{equation}
gd(i,j) = g(xp(i+1),yp(j+1)) \,,\quad
	i \leq Nx-1 \,,\quad j \leq Ny-1 \,.
\end{equation}
A vector function $(tx(x,y),ty(x,y)$ with spatial weight $1/d$ (tangent on
dual grid) is discretized at cell edge centers:
\begin{align}
txd(i,j) = & tx(xp(i+1),yd(j)) \,,\quad i \leq Nx-1 \,,\quad j \leq Ny \,;\\
tyd(i,j) = & ty(xd(i),yp(j+1)) \,,\quad i \leq Nx \,,\quad j \leq Ny \,. 
\end{align}
A vector function $(nx(x,y),ny(x,y)$ with spatial weight $1/d^2$ (normal
on the dual grid) is discretized at cell face centers:
\begin{align}
nxd(i,j) = & nx(xd(i),yp(j+1)) \,,\quad i \leq Nx \,,\quad j \leq Ny-1 \,;\\
nyd(i,j) = & ny(xp(i+1),yd(j)) \,,\quad i \leq Nx-1 \,,\quad j \leq Ny \,. 
\end{align}

Examples of 2D scalar and vector fields can be generated using
{\tt ScalarField2p.m},
{\tt ScalarField2d.m},
{\tt TangentField2p.m},
{\tt TangentField2d.m},
{\tt NormalField2p.m},
{\tt NormalField2d.m},
{\tt DensityField2p.m},
{\tt DensityField2d.m}
and tested using {\tt TestFields2pd.m}

\subsection{The Star Operators}

For simplicity the discussion of star operators will begin for material
properties that are constant. These star operators will map quantities
defined on the primal grid to the dual grid while their inverses do the
opposite.  For constant and isotropic materials the star operators are
multiplication by a constant for quantities that are defined at the same
point in the grids.  For scalar variables this is done by multiplication
by a constant $a>0$ with spatial dimension $1/d^3$. For vectors a constant
diagonal matrix
\[
{\bf A} =
\left[
\begin{matrix}
A11 & 0  \\
0 & A22 \\
\end{matrix}
\right]
\]
with $A11>0$ and $A22>0$ spatial dimension $1/d$ will be used.
If $A11 \neq A22$ then there is a simple anisotropy.  Away from the
boundaries of the region, the two star operators are inverses of each
other.  For boundary value problems, the mismatch between the sizes of the
primal and dual grids will be used to represent the boundary conditions.
The star operators will be given in pairs that are essentially inverse of
each other, first the mapping from the primal grid to the dual grid and then
the inverse.

If $gd = \star fp$ then
\begin{equation}
gd(i,j) = a \,fp(i+1,j+1)
\,,\quad 1 \leq i \leq N_x-1
\,,\quad 1 \leq j \leq N_y -1
\,.
\end{equation}
If $fp = \star gd$ then
\begin{equation}
fp(i,j) = \frac{1}{a} \, gd(i-1,j-1)
\,,\quad 2 \leq i \leq N_x
\,,\quad 2 \leq j \leq N_y
\,.
\end{equation}
In this case $fp$ is not defined on the boundary of the primal grid.

If $fd = \star gp$ then
\begin{equation}
gp(i,j) = a\,fd(i,j)
\,,\quad 1 \leq i \leq N_x
\,,\quad 1 \leq j \leq N_y 
\,.
\end{equation}
If $gp = \star fd$ then
\begin{equation}
gp(i,j) = \frac{1}{a}\,fd(i,j)
\,,\quad 1 \leq i \leq N_x
\,,\quad 1 \leq j \leq N_y 
\,.
\end{equation}

If $(txd,tyd) = \star (nxp,nyp)$ then
\begin{align}
txd(i,j) &= \frac{1}{A11} \, nxp(i+1,j)
\,,\quad 1 \leq i \leq N_x-1
\,,\quad 1 \leq j \leq N_y
\,. \nonumber \\
tyd(i,j) &= \frac{1}{A22} \, nyp(i,j+1)
\,,\quad 1 \leq i \leq N_x
\,,\quad 1 \leq j \leq N_y-1
\,.
\end{align}

If $(nxp,nyp) = \star (txd,tyd)$ then
\begin{align}
nxp(i,j) &= A11 \, txd(i-1,j)
\,,\quad 2 \leq i \leq N_x
\,,\quad 1 \leq j \leq N_y
\,. \nonumber \\
nyp(i,j) &= A22 \, tyd(i,j-1)
\,,\quad 1 \leq i \leq N_x
\,,\quad 2 \leq j \leq N_y
\,.
\end{align}

If $(nxd, nyd) = \star(txp,typ)$ then
\begin{align}
nxd(i,j) &= A11 \, txp(i,j+1)
\,,\quad 1 \leq i \leq N_x
\,,\quad 1 \leq j \leq N_y-1
\,. \nonumber \\
nyd(i,j) &= A22 \, typ(i+1,j)
\,,\quad 1 \leq i \leq N_x-1
\,,\quad 1 \leq j \leq N_y
\,. \label{nxd star txp}
\end{align}

If $(txp,typ) = \star (nxd, nyd)$ then
\begin{align}
txp(i,j) &= \frac{1}{A11} \, nxd(i,j-1)
\,,\quad 1 \leq i \leq N_x
\,,\quad 2 \leq j \leq N_y
\,. \nonumber \\
typ(i,j) &= \frac{1}{A22} \, nyd(i-1,j)
\,,\quad 2 \leq i \leq N_x
\,,\quad 1 \leq j \leq N_y
\,.
\end{align}

\subsection{2D Differential Operators}

The 2D discrete differential operators are the gradient and divergence
which are given by {\tt Grad2p.m} and {\tt Grad2d.m} and {\tt Div2p.m}
and {\tt Div2d.m}. All the differential operators use centered differences.
The operators on the primal and dual grids differ in their indexing.
The programs {\tt TestGradDiv2p.m} and {\tt TestGradDiv2d.m} show that
the gradient and divergence operators are second order accurate on the
primal and dual grids.

\subsection{Scalar Wave Equation}

The discretization of the 2D scalar wave equation is the same as the
discretization of 3D scalar wave equation given in \eqref{Discrete 3D Wave}
and has been implemented in {\tt Wave2D.m}. The wave equation is written
as a first order system of differential equations that are discretized
using staggered space-time grids and {\tt Grad2p} and {\tt Div2d}.
In general the approximate solutions of the discrete wave equation are
second order accurate. For some cases the solutions are forth order accurate
and there is at least one example where the solution is exact
(set $m1 = n1 = 1$ and eliminate the $m2$, $n2$ part of the test
solution to see this).
Both conservation laws are constant to at least 1 part in $10^{15}$
when the star operator is trivial.

What star was used in Wave2D.m? See \eqref{nxd star txp}.

\subsection{Boundary Conditions}

The needs a rewrite.

Figures \ref{Scalar Vector Primal} and \ref{Scalar Vector Dual} illustrates
the position of the scalar field $f$ and the scalar function $g$ which is
the Laplacian of $f$, that is the divergence of the gradient of $f$.
It also illustrates the positions of the boundary conditions which must
specify $g$ on the boundary:
\begin{align*}
g_{1,j}     \,,\, & 1 \leq j \leq Ny+1 \,\;\\
g_{Nx+1,j}  \,,\, & 1 \leq j \leq Ny+1 \,\;\\
g_{i,1}     \,,\, & 1 \leq i \leq Nx+1 \,\;\\
g_{i,Ny+1,} \,,\, & 1 \leq i \leq Nx+1 \,.
\end{align*}
Note that the values of $g$ are defined twice at the corner points
of the region:
$(1,1)$; $(1,Ny+1)$; $(Nx+1,1)$ and $(Nx+1,Ny+1)$.
In fact in standard discrete boundary value problems these values of $g$
are not needed and can be assigned any value. On the other hand creating
a data structure that doesn't have these values creates a programming mess.
The figures were generated using {\tt Figure2DDiv.m}, {\tt Figure2DGrad.m}
and {\tt Figure2DLap.m}.

The typical boundary condition is of mixed or Robin type, that is,
\[
\alpha \, \vec{n} \bdot \vec{v} + \beta f = \gamma \,,
\]
which in the discrete setting becomes
\begin{align*}
\text{at  } y = 0 \quad &
\alpha_{1,j} \, vy_{1,j} +\beta_{1,j} \, g_{1,j} =
\gamma_{1,j} \,; \\
\text{at  } y = 1 \quad &
\alpha_{Nx+1,j} \, vy_{Nx+1,j} +\beta_{Nx+1,j} \, g_{Nx+1,j} =
\gamma_{Nx+1,j} \,; \\
\text{at  } x = 0 \quad &
\alpha_{i,1} \, vy_{i,1} +\beta_{i,1} \, g_{i,1} =
\gamma_{i,1} \,; \\
\text{at  } x = 1 \quad &
\alpha_{i,Ny+1} \, vy_{i,Ny+1} +\beta_{i,Ny+1} \, g_{i,Ny+1} =
\gamma_{i,Ny+1} \,. 
\end{align*}
These equations can be trivially solved for
$ g_{i,1} $, $ g_{i,Ny+1} $, $ g_{1,j} $, $ g_{Nx+1,j} $ 
which could give two different values of $f$ for the corner points.
In the typical explicit time stepping algorithms for wave equations,
these values are never used.

Check this, probably not correct.
For Dirichlet boundary conditions
the program {\tt Wave2D.m} confirms that the solution of the 2D wave
equation is second order accurate while the program {\tt Wave2DExact.m}
illustrates some cases where the solutions are accurate up to some
small multiple of {\tt eps}.


%
\newpage \clearpage
\bibliography{citations}
\begin{appendix}
\newpage \clearpage
\setcounter{equation}{0}
\section{Energy Preserving Discretizations of the Harmonic Oscillator
\label{C-N} }

Here the well-know fact that the Crank-Nicholson discretization conserves
the discrete analog of the energy for the harmonic oscillator is shown.
It is also shown that the methods introduced in \cite{WanBihloNave2015}
produce a discretization that is equivalent to the Crank-Nicholson
discretization.

\subsection{Conserving the Simple Energy}

The Crank-Nicholson discretization does preserve the simple
energy \eqref{Simple Conserved}:
\[
\frac{ u_{n+1}-u_n }{\Delta t} = - \omega \, \frac{v_{n+1}+v_n}{2} \,,\quad
\frac{ v_{n+1}-v_n }{\Delta t} =   \omega \, \frac{u_{n+1}+u_n}{2} \,.
\]
This gives a discretization of the second order differential equation:
\begin{equation}\label{2-CN}
\frac{u_{n+2} - 2 \, u_{n+1} +u_n }{\Delta t^2} + 
\omega^2 \frac{u_{n+2} + 2 \, u_{n+1} +u_n }{4} = 0  \,.
\end{equation}
Then
\begin{align*}
C_{n+1}^2 - C_n^2 & = 
\frac{1}{2} \left(v_{n+1} + v_n\right) \, \left(v_{n+1} - v_n\right) + 
\frac{1}{2} \left(u_{n+1} + u_n\right) \, \left(u_{n+1} - u_n\right) \\
& =
  \frac{\Delta t \, \omega}{4} \,
  \left(v_{n+1} + v_n\right) \,
  \left(u_{n+1} + u_n\right) -
  \frac{\Delta t \, \omega}{4} \,
  \left(u_{n+1} + u_n\right) \,
  \left(v_{n+1} + v_n\right) \equiv 0 \,,
\end{align*}
so $C_n$ is conserved.  Write the system as
\begin{align*}
u_{n+1} + \frac{\Delta t \, \omega}{2} v_{n+1} & =
   u_n  - \frac{\Delta t \, \omega}{2} v_n \\
v_{n+1} - \frac{\Delta t \, \omega}{2} u_{n+1} & =
   v_n  + \frac{\Delta t \, \omega}{2} u_n
\,,
\end{align*}
so that the scheme is implicit, that is it involves the inversion of a
$2\times2$ matrix. The coefficient matrix is always invertible, so there
is no restriction on the size of $\Delta t$, that is, the scheme is
unconditionally stable.

\subsection{The Conservation Law First}

Following the discussion in \cite{WanBihloNave2015} it is easy to show that
the only reasonable discretization that conserves the simple conservation
law \eqref{Simple Conserved} is equivalent to the Crank-Nicholson
discretization. First compute using \eqref{Simple Conserved} that
\[
C_{n+1}^2 - C_n^2 = 
\left(u_{n+1}-u_n\right)
\left(u_n+u_{n+1}\right) +
\left( v_{n+\frac{3}{2}} - v_{n-\frac{1}{2}} \right)
\frac{v_{n-\frac{1}{2}}+2 v_{n+\frac{1}{2}}+v_{n+\frac{3}{2}}}{4} \,.
\]
Choosing
\[
\frac{u_{n+1}-u_n}{\Delta t} =
- \omega \frac{v_{n-\frac{1}{2}}+2 v_{n+\frac{1}{2}}+v_{n+\frac{3}{2}}}{4}
\]
and
\[
\frac{v_{n+\frac{3}{2}} - v_{n-\frac{1}{2}}}{2 \, \Delta t} =
\omega \frac{u_n+u_{n+1}}{2}
\]
will make the $C_n$ constant.
If $\alpha =  \Delta t \, \omega / 2 $ then these equations can be written
\begin{align*}
u_{n+1} + \frac{\alpha}{2} v_{n+3/2} & =
	u_n - \alpha v_{n+1/2} -\frac{\alpha}{2} v_{n-1/2} \,, \\
- 2 \, \alpha \, u_{n+1}  + v_{n+3/2} & = 2 \, \alpha u_{n} +v_{n-1/2} \,.
\end{align*}
So the difference equations are implicit.

It is easy to check that $u_n$ satisfies the second order difference
equation \eqref{2-CN}. Unfortunately, this discretization produces the
same $u_n$ values as the Crank-Nicholson scheme but with a greater
computational cost.  Setting
\[
v_n = \frac{v_{n+1/2}+v_{n-1/2}}{2} \,,
\]
converts this scheme along with it's conserved quantity to the
Crank-Nicholson scheme along with it's conserved quantity.  

\newpage \clearpage
\setcounter{equation}{0}

\section{Details for Discrete Conserved Quantities \label{Appendix Details}}

\subsection{Scalar Wave}

As before a second order discrete equation and a second order average
will be needed
\begin{align*}
\frac{u^{n+1} - 2\, u^{n} + u^{n-1}}{\dt^2}
& = a^{-1}\,\frac{\DIVGstar v^{n+\half}- \DIVGstar v^{n-\half}}{\dt} \\
& = a^{-1}\, \DIVGstar \frac{v^{n+\half}- v^{n-\half}}{\dt} \\
& = a^{-1}\, \DIVGstar {\bf A} \GRAD u^n \\
\frac{u^{n+1} + 2\, u^{n} + u^{n-1}}{4}
& =   u^n + \frac{u^{n+1} - 2\, u^{n} + u^{n-1}}{4} \\
& =   u^n + \frac{\dt^2}{4}\frac{u^{n+1} - 2\, u^{n} + u^{n-1}}{\dt^2} \\
& =   u^n + \frac{\dt^2}{4} a^{-1}\, \DIVGstar {\bf A} \GRAD u^n 
\\
\end{align*}

To find a conserved quantity let
\begin{align*}
C1_{n+1/2} & = \norm{\frac{u^{n+1} + u^{n}}{2}}_\Nodes^2 \,, \\
C2_{n+1/2} & = \norm{ v^{n+1/2}}_\StarFaces^2 \,, \\
C3_{n+1/2} & = \Delta t^2 \norm{ a^{-1} \, \DIVGstar \, v^{n+1/2} }_\Nodes^2 \,.
\end{align*}
As before compute:
\begin{align*}
C1_{n+1/2}- C1_{n-1/2}
& = \langle \frac{u^{n+1} + 2 \, u^{n}  + u^{n-1}}{4} , u^{n+1} - u^{n-1} \rangle_\Nodes \\
& = \langle u^{n} + \frac{\Delta t^2}{4} a^{-1} \, \DIVGstar \, {\bf A} \GRAD u^{n} , u^{n+1} - u^{n-1} \rangle_\Nodes \,;  \\
& =
\langle u^{n} , u^{n+1} - u^{n-1} \rangle_\Nodes 
+
\frac{\Delta t^2}{4}
\langle a^{-1} \, \DIVGstar \, {\bf A} \GRAD u^{n}
, u^{n+1} - u^{n-1} \rangle_\Nodes \,;
\end{align*}
Using the adjoint equation \eqref{Adjoint 1} gives
\begin{align*}
C2_{n+1/2}- C1_{n-1/2}
& = \langle v^{n+1/2} + v^{n-1/2} , v^{n+1/2} - v^{n-1/2} \rangle_\StarFaces \\
& = \langle v^{n+1/2} + v^{n-1/2} , \Delta t \, {\bf A} \GRAD u^n \rangle_\StarFaces \\
& = - \Delta t \, \langle a^{-1} \, \DIVGstar \, v^{n+1/2} + a^{-1} \, \DIVGstar \, v^{n-1/2} , u^n \rangle_\Nodes \\
& = -\Delta t \, \langle \frac{u^{n+1}-u^{n-1}}{\Delta t} , u^n \rangle_\Nodes \\
& = - \langle u^n , u^{n+1}-u^{n-1} \rangle_\Nodes \,;
\end{align*}
Also
\begin{align*}
C3_{n+1/2}- C1_{n-1/2} & = \Delta t^2
\langle a^{-1} \, \DIVGstar \, v^{n+1/2} - a^{-1} \, \DIVGstar \, v^{n-1/2} \,,\, a^{-1} \, \DIVGstar \, v^{n+1/2} + a^{-1} \, \DIVGstar \, v^{n-1/2} \rangle_\Nodes \\
& = \Delta t^ 2 \langle a^{-1} \, \DIVGstar \left( v^{n+1/2} - v^{n+1/2} \right) \,,\, \frac{u^{n+1}-u^{n-1}}{\Delta t} \rangle_\Nodes \\
& = \Delta t^ 2 \langle - \Delta t \, a^{-1} \, \DIVGstar \, -{\bf A} \GRAD u^{n} \,,\,
	\frac{u^{n+1}-u^{n-1}}{\Delta t} \rangle_\Nodes \\
& = \Delta t^ 2 \langle a^{-1} \, \DIVGstar \, {\bf A} \GRAD u^{n} \,,\, u^{n+1}-u^{n-1} \rangle_\Nodes \,.
\end{align*}
Consequently $C = C1 + C2 - C3/4$ is a conserved quantity:
\[
C_{n+1/2} = \norm{\frac{u^{n+1} + u^{n}}{2}}^2 + \norm{ v^{n+1/2}}^2 
 - \frac{\Delta t^2}{4} \norm{ a^{-1} \, \DIVGstar \, v^{n+1/2} }^2 \,.
\]
This implies that 
\[
C_{n+1/2} \geq \norm{\frac{u^{n+1} + u^{n}}{2}}^2 +
\left( 1 - \frac{\Delta t^2}{4} \norm{a^{-1} \, \DIVGstar}^2 \right)
\norm{ v^{n+1/2}}^2  \,.
\]
So $C_{n+1/2} \geq 0$ for $\Delta t$ sufficiently small provided
$\norm{a^{-1} \, \DIVGstar}$ is finite.

Next look for an analog $C_{n}$ of the scaler conserved quantity
\begin{align*}
C1_{n} & = \norm{\frac{v^{n+1/2} + v^{n-1/2}}{2}}_\StarFaces^2 \,, \\
C2_{n} & = \norm{u^{n}}_\Nodes^2 \,, \\
C3_{n} & = \Delta t^2 \norm{{\bf A} \GRAD \, u^{n} }_\StarFaces^2 \,. 
\end{align*}
First compute
\begin{align*}
C1_{n+1}- C1_{n}
& = \langle
\frac{v^{n+3/2} + 2 \, v^{n 1/2} + v^{n-1/2}}{4} \,,\,
   v^{n+3/2} -v^{n-1/2} \rangle_\StarFaces \\
& =
\langle v^{n+1/2} , v^{n+3/2} -v^{n-1/2} \rangle_\StarFaces
+ \frac{\Delta t^2}{4}
\langle
{\bf A} \GRAD \, a^{-1} \, \DIVGstar \, v^{n+1/2} \,,\,
v^{n+3/2} -v^{n-1/2} \rangle\StarFaces
\,.
\end{align*}
Using the adjoint equation \eqref{Adjoint 1} gives
\begin{align*}
C2_{n+1} - C2_{n} 
& = \langle u^{n+1}-u^{n} \,,\, u^{n+1}+u^{n} \rangle_\Nodes \\
& = \langle \Delta t \,  a^{-1} \, \DIVGstar \, v^{n+1/2} \,,\, u^{n+1}+u^{n} \rangle_\Nodes \\
& = \Delta t \langle v^{n+1/2} \,,\, -{\bf A} \GRAD u^{n+1}-{\bf A} \GRAD u^{n} \rangle \\
& = \Delta t \langle v^{n+1/2} \,,\,
	-\frac{v^{n+3/2}- v^{n-1/2}}{\Delta t} \rangle \\
& = - \langle v^{n+1/2} \,,\,
	v^{n+3/2}- v^{n-1/2} \rangle  \,.
\end{align*}
Also
\begin{align*}
C3_{n+1}- C3_{n} & = \Delta t^2 
\langle
   {\bf A} \GRAD u^{n+1}-{\bf A} \GRAD u^n  \,,\,  {\bf A} \GRAD u^{n+1}+{\bf A} \GRAD u^n \rangle_\StarFaces \\
& = \Delta t^2 \langle  {\bf A} \GRAD u^{n+1} - {\bf A} \GRAD u^n  \,,\, \frac{v^{n+3/2}-v^{n-1/2}}{\Delta t}
\rangle_\StarFaces \\
& = \Delta t^2 \langle
   \Delta t \, {\bf A} \GRAD a^{-1} \, \DIVGstar \, v^{n+1/2} \,,\, \frac{v^{n+3/2}-v^{n-1/2}}{\Delta t}
\rangle \\
& = \Delta t^2 \langle
   {\bf A} \GRAD a^{-1} \, \DIVGstar \, v^{n+1/2} \,,\, v^{n+3/2}-v^{n-1/2} \rangle \,.
\end{align*}
Consequently $C_n = C1_n + C2_n - C3_n/4$ is a conserved quantity: 
\[
C_n =  \norm{u^{n}}^2 
 -\frac{\Delta t^2}{4} \norm{ {\bf A} \GRAD \, u^{n} }^2 
 +\norm{\frac{v^{n+1/2} + v^{n-1/2}}{2}}^2 \,.
\]
This implies that
\[
\norm{C_n} \ge	 \left(1 - \frac{\Delta t^2}{4} \norm{{\bf A} \GRAD}^2\right) \norm{u^n}^2
	+ \norm{\frac{v^{n+1/2} + v^{n-1/2}}{2}}^2 \,,
\]
so $\norm{C_n}$ is positive for sufficiently small $\Delta t$ if
$\norm{ {\bf A} \GRAD \, u^{n} }$ is finite.

\subsection{Maxwell}

To study conserved quantities for Maxwell's equations the second order
disctete difference and average will be needed:
\begin{align*}
\frac{\vec{E}^{n+1} - 2\, \vec{E}^{n} + \vec{E}^{n-1}}{\dt^2}
& = - \epsilon^{-1}\, \CURLstar \mu^{-1} \CURL \vec{E}^n \\
\frac{\vec{E}^{n+1} + 2\, \vec{E}^{n} + \vec{E}^{n-1}}{4}
& =   \vec{E}^n -
\frac{\dt^2}{4} \epsilon^{-1}\, \CURLstar \mu^{-1} \CURL \vec{E}^n 
\end{align*}

To find a conserved quantity $C_{n+1/2}$ let
\begin{align*}
C1_{n+1/2} & = \norm{\frac{\vec{E}^{n+1} + \vec{E}^{n}}{2}}_\Edges^2 \,, \\
C2_{n+1/2} & = \norm{ \vec{H}^{n+1/2}}_\StarEdges^2 \,, \\
C3_{n+1/2} & = \Delta t^2 \norm{\epsilon^{-1} \, \CURLstar \, \vec{H}^{n+1/2} }_\Edges^2  \,.
\end{align*}

As before compute:
\begin{align*}
C1_{n+1/2}- C1_{n-1/2}
& = \langle \frac{\vec{E}^{n+1} + 2 \, \vec{E}^{n}  + \vec{E}^{n-1}}{4} , \vec{E}^{n+1} - \vec{E}^{n-1} \rangle_\Edges \\
& = \langle \vec{E}^{n} - \frac{\Delta t^2}{4} \epsilon^{-1} \, \CURLstar \, \mu^{-1}
\CURL \vec{E}^{n} , \vec{E}^{n+1} - \vec{E}^{n-1} \rangle_\Edges \,;  \\
& =
\langle \vec{E}^{n} , \vec{E}^{n+1} - \vec{E}^{n-1} \rangle_\Edges 
- \frac{\Delta t^2}{4}
\langle \epsilon^{-1} \, \CURLstar \, \mu^{-1} \CURL \vec{E}^{n}
, \vec{E}^{n+1} - \vec{E}^{n-1} \rangle_\Edges \,;
\end{align*}

Using the adjoint equation \eqref{Adjoint 1} gives
\begin{align*}
C2_{n+1/2}- C1_{n-1/2}
& = \langle \vec{H}^{n+1/2} + \vec{H}^{n-1/2} \,, \vec{H}^{n+1/2} - \vec{H}^{n-1/2} \rangle_\StarFaces \\
& = \langle \vec{H}^{n+1/2} + \vec{H}^{n-1/2} \,, - \Delta t \, \mu^{-1} \CURL \vec{E}^n \rangle_\StarFaces \\
& = - \Delta t \, \langle \epsilon^{-1} \, \CURLstar \, \vec{H}^{n+1/2} + \epsilon^{-1} \, \CURLstar \, \vec{H}^{n-1/2} , \vec{E}^n \rangle_\Edges \\
& = -\Delta t \, \langle \frac{\vec{E}^{n+1}-\vec{E}^{n-1}}{\Delta t} , \vec{E}^n \rangle_\Edges \\
& = - \langle \vec{E}^n , \vec{E}^{n+1}-\vec{E}^{n-1} \rangle_\Edges \,;
\end{align*}
Also
\begin{align*}
C3_{n+1/2}- C1_{n-1/2} & = \Delta t^2
\langle \epsilon^{-1} \, \CURLstar \, \vec{H}^{n+1/2} - \epsilon^{-1} \, \CURLstar \, \vec{H}^{n-1/2} \,,\, \epsilon^{-1} \, \CURLstar \, \vec{H}^{n+1/2} + \epsilon^{-1} \, \CURLstar \, \vec{H}^{n-1/2} \rangle_\Edges \\
& = \Delta t^ 2 \langle \epsilon^{-1} \, \CURLstar \left( \vec{H}^{n+1/2} - \vec{H}^{n+1/2} \right) \,,\, \frac{\vec{E}^{n+1}-\vec{E}^{n-1}}{\Delta t} \rangle_\Edges \\
& =\Delta t^ 2 \langle - \Delta t \, \epsilon^{-1} \, \CURLstar \, \mu^{-1} \CURL \vec{E}^{n} \,,\,
	\frac{\vec{E}^{n+1}-\vec{E}^{n-1}}{\Delta t} \rangle_\Edges \\
& = - \Delta t^ 2 \langle \epsilon^{-1} \, \CURLstar \, \mu^{-1} \CURL \vec{E}^{n} \,,\, \vec{E}^{n+1}-\vec{E}^{n-1} \rangle_\Edges \,.
\end{align*}
Consequently $C = C1 + C2 - C3/4$ is a conserved quantity:
\[
C_{n+1/2} =
\norm{\frac{\vec{E}^{n+1} + \vec{E}^{n}}{2}}_\Edges^2
+ \norm{ \vec{H}^{n+1/2}}_\StarEdges^2 
- \frac{\Delta t^2}{4}
\norm{ \epsilon^{-1} \, \CURLstar \, \vec{H}^{n+1/2} }_\Edges^2 \,.
\]
This implies that 
\[
C_{n+1/2} \geq
\norm{\frac{\vec{E}^{n+1} + \vec{E}^{n}}{2}}_\Edges^2 +
\left( 1 - \frac{\Delta t^2}{4} \norm{\epsilon^{-1} \, \CURLstar}^2 \right)
\norm{ \vec{H}^{n+1/2}}_\StarEdges^2  \,.
\]
So $C_{n+1/2} \geq 0$ for $\Delta t$ sufficiently small provided
$\norm{\epsilon^{-1} \, \CURLstar}$ is finite.

Next look for a conserved quantity $C_{n}$:
\begin{align*}
C1_{n} & = \norm{\frac{\vec{H}^{n+1/2} + \vec{H}^{n-1/2}}{2}}_\StarFaces^2 \,, \\
C2_{n} & = \norm{\vec{E}^{n}}_\Edges^2 \,, \\
C3_{n} & = \Delta t^2 \norm{\mu^{-1} \CURL \, \vec{E}^{n} }_\StarFaces^2 \,. 
\end{align*}
First compute
\begin{align*}
C1_{n+1}- C1_{n}
& = \langle
\frac{\vec{H}^{n+3/2} + 2 \, \vec{H}^{n 1/2} + \vec{H}^{n-1/2}}{4} \,,\,
   \vec{H}^{n+3/2} -\vec{H}^{n-1/2} \rangle_\StarFaces \\
& =
\langle \vec{H}^{n+1/2} , \vec{H}^{n+3/2} -\vec{H}^{n-1/2} \rangle_\StarFaces
+ \frac{\Delta t^2}{4}
\langle
\mu^{-1} \CURL \, \epsilon^{-1} \, \CURLstar \, \vec{H}^{n+1/2} \,,\,
\vec{H}^{n+3/2} -\vec{H}^{n-1/2} \rangle_\StarFaces
\,.
\end{align*}
Using the adjoint equation \eqref{Adjoint 1} gives
\begin{align*}
C2_{n+1} - C2_{n} 
& = \langle \vec{E}^{n+1}-\vec{E}^{n} \,,\, \vec{E}^{n+1}+\vec{E}^{n} \rangle_\Edges \\
& = \langle \Delta t \,  \epsilon^{-1} \, \CURLstar \, \vec{H}^{n+1/2} \,,\, \vec{E}^{n+1}+\vec{E}^{n} \rangle_\Edges \\
& = \Delta t \langle \vec{H}^{n+1/2} \,,\, -\mu^{-1} \CURL \vec{E}^{n+1}-\mu^{-1} \CURL \vec{E}^{n} \rangle_\StarFaces \quad\text{(adjoint)} \\
& = \Delta t \langle \vec{H}^{n+1/2} \,,\,
	-\frac{\vec{H}^{n+3/2}- \vec{H}^{n-1/2}}{\Delta t} \rangle_\StarFaces \\
& = - \langle \vec{H}^{n+1/2} \,,\,
	\vec{H}^{n+3/2}- \vec{H}^{n-1/2} \rangle_\StarFaces  \,.
\end{align*}
Also
\begin{align*}
C3_{n+1}- C3_{n} & = \Delta t^2 
\langle
   \mu^{-1} \CURL \vec{E}^{n+1}-\mu^{-1} \CURL \vec{E}^n  \,,\,  \mu^{-1} \CURL \vec{E}^{n+1}+\mu^{-1} \CURL \vec{E}^n \rangle_\StarFaces \\
& = \Delta t^2 \langle  \mu^{-1} \CURL \vec{E}^{n+1} - \mu^{-1} \CURL \vec{E}^n  \,,\, \frac{\vec{H}^{n+3/2}-\vec{H}^{n-1/2}}{\Delta t}
\rangle_\StarFaces \\
& = \Delta t^2 \langle
   \Delta t \, \mu^{-1} \CURL \epsilon^{-1} \, \CURLstar \, \vec{H}^{n+1/2} \,,\, \frac{\vec{H}^{n+3/2}-\vec{H}^{n-1/2}}{\Delta t}
\rangle_\StarFaces \\
& = \Delta t^2 \langle
   \mu^{-1} \CURL \epsilon^{-1} \, \CURLstar \, \vec{H}^{n+1/2} \,,\, \vec{H}^{n+3/2}-\vec{H}^{n-1/2} \rangle_\StarFaces \,.
\end{align*}
Consequently $C_n = C1_n + C2_n - C3_n/4$ is a conserved quantity: 
\[
C_n =  \norm{\vec{E}^{n}}_\Edges^2 
 -\frac{\Delta t^2}{4} \norm{ \mu^{-1} \CURL \, \vec{E}^{n} }_\StarFaces^2 
 +\norm{\frac{\vec{H}^{n+1/2} + \vec{H}^{n-1/2}}{2}}_\StarFaces^2 \,.
\]
This implies that
\[
\norm{C_n} \ge	
\left(1 - \frac{\Delta t^2}{4}
\norm{\mu^{-1} \CURL}^2\right)
\norm{\vec{E}^n}_\Edges^2
+ \norm{\frac{\vec{H}^{n+1/2} + \vec{H}^{n-1/2}}{2}}_\StarFaces^2 \,,
\]
so $\norm{C_n}$ is positive for sufficiently small $\Delta t$ if
$\norm{ \mu^{-1} \CURL \, \vec{E}^{n} }$ is finite.

The codes {\tt Maxwell.m} and {\tt MaxwellStar.m} confirm that our algorithms
conserve $C_{n+1/2}$ and $C_n$ to two parts in $10^{16}$. Additinally,
the divergence of the curl of the electric and magnetic fields are constant
to one part in $10^{14}$ when there are no sources.

\newpage \clearpage
\setcounter{equation}{0}
\section{Conservation Laws and Positive Solutions}

Conservation laws that say the total amount of some positive substance
is conserved play an important role in modeling using partial differential
equations, for example the Navier-Stokes equations
\cite{PeyretTaylor83}(equations 1.5, 1.6 and 1.7) can be put into this form.
To provided some insight into discretizing such conservation laws,
two important but simple cases will be considered.  For a similar discussion
see Chapeter 11 in \cite{Oliver16}.

\subsection{Transport}

The transport equation in one dimension is given by
\[
\frac{\partial \rho}{\partial t} +
\frac{\partial \, v \, \rho}{\partial x} = 0 \,,
\]
where $\rho = \rho(x,t)$ is a density and $v = v(x)$ is the velocity of
transport. An important assumption is that $\rho \geq 0$ as it typically
represents the density of some substance.  The general solution of this
equation is
\[
\rho (x,t) = w(x-vt) \,,
\]
where $w(x) = \rho(x,0)$ is the initial data. This solution is
a right translation of $w(x)$. This equation also has an important
conservation law:
\[
\int_{-\infty}^{\infty} \rho (x,t) \, d x =
\int_{-\infty}^{\infty} w(x) \, d x  \,.
\]
The conserved quantity is the total amount of material being transported.
Also note that if $w(x) \geq 0$ then $\rho(x,t) \geq 0$ for all $t$. These
two properties are central to this discussion.
Our interest is in finite difference discretizations of equations that
have a similar conservation law and maintain the positivity of the solution.

We assume that $\Delta x > 0 $ and use two grids: a primal grid
$x_i = i \, \Delta x$ that has cells $[x_i,x_{i+1}]$ and a grid of
cell centers $x_{i+\half} = (i+\half) \, \Delta x$ where $-\infty < i < \infty$.
Note that if $\rho$ is a density then it has spatial dimension $1/d^k$ in
a space of dimension $k$ suggesting that $\rho$ should be in a cells. If
a primal grid is chosen then the discretization of $\rho$ is
\[
\rho^{n+\half}_{i+\half} \,.
\]
We will use the conservation of material
\[
\Delta x \, \rho^{n+\half}_{i+\half} 
\]
in a cell to discretize this equation as
\[
\Delta x \, \rho^{n+3/2}_{i+\half} = 
	\Delta x \, \rho^{n+\half}_{i+\half} 
	+ \Delta t \, v_i \, \rho^{n+\half}_{i-\half}
	- \Delta t \, v_{i+1} \, \rho^{n+\half}_{i+\half} \,.
\]
Rewrite this as
\[
\frac{\rho^{n+3/2}_{i+\half} - \rho^{n+\half}_{i+\half}}{\Delta t}
+ \frac{v_{i+1} \, \rho^{n+\half}_{i+\half} -v_i \, \rho^{n+\half}_{i-\half}}
{\Delta x} =  0 \,,
\]
to see that the discretization is a first order approximation of the
differential equation. As an update of the density the equation becomes
\[
\rho^{n+3/2}_{i+\half} =
\rho^{n+\half}_{i+\half} 
	+ \frac{\Delta t}{\Delta x} \, v_i \, \rho^{n+\half}_{i-\half}
	- \frac{\Delta t}{\Delta x} \, v_{i+1} \, \rho^{n+\half}_{i+\half} \,.
\]
Now if 
\[
\frac{\Delta t}{\Delta x} \, v_i \geq 0  \,,\quad
1 - \frac{\Delta t}{\Delta x} \, v_{i+1} \geq 0 \,,
\]
that is if 
\[
v_i \geq 0 \,,\quad
\frac{\Delta t}{\Delta x} \, v_{i+1} \leq 1 \,,
\]
then the discretization preserves the positivity of the discrete solution
and is the well known upwind scheme. This scheme is not useful if 
the velocity $v = v(x)$ has both negative and positive values.
To fix this consider $v$ rather than $\rho$.

\begin{figure}
\begin{center}
\begin{tabular}{cc}
   \includegraphics[width=2.75in,trim = 0 0 0 350]{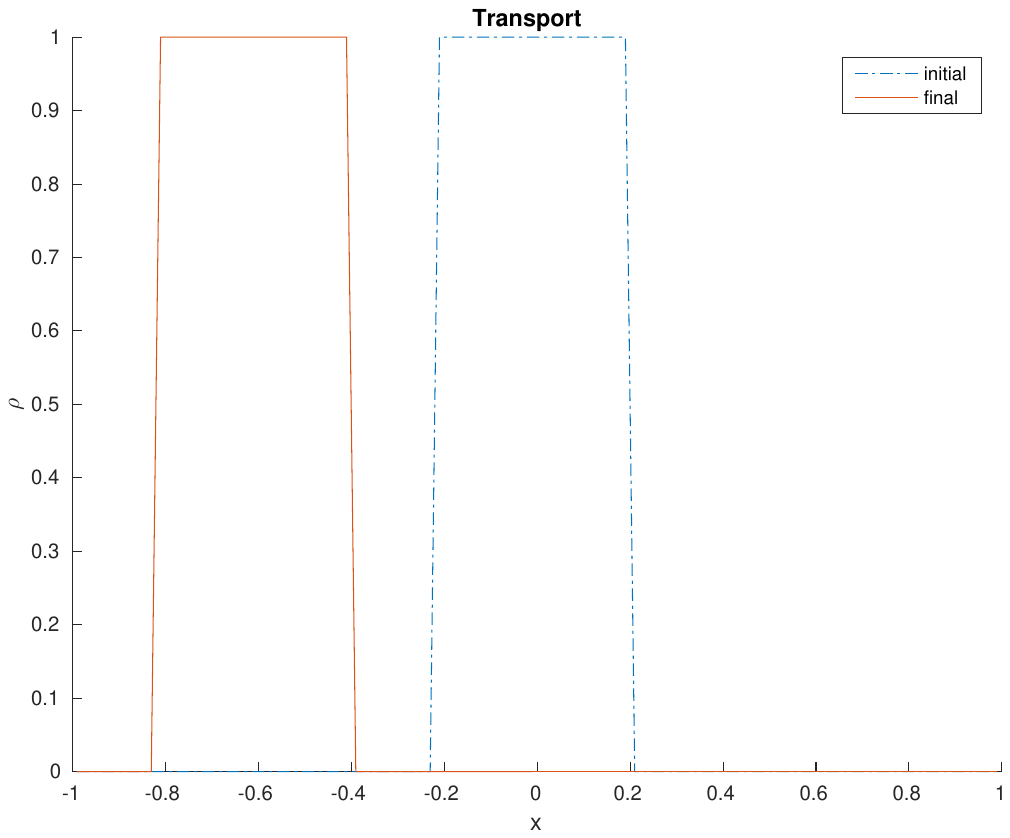} &
   \includegraphics[width=2.75in,trim = 0 0 0 350]{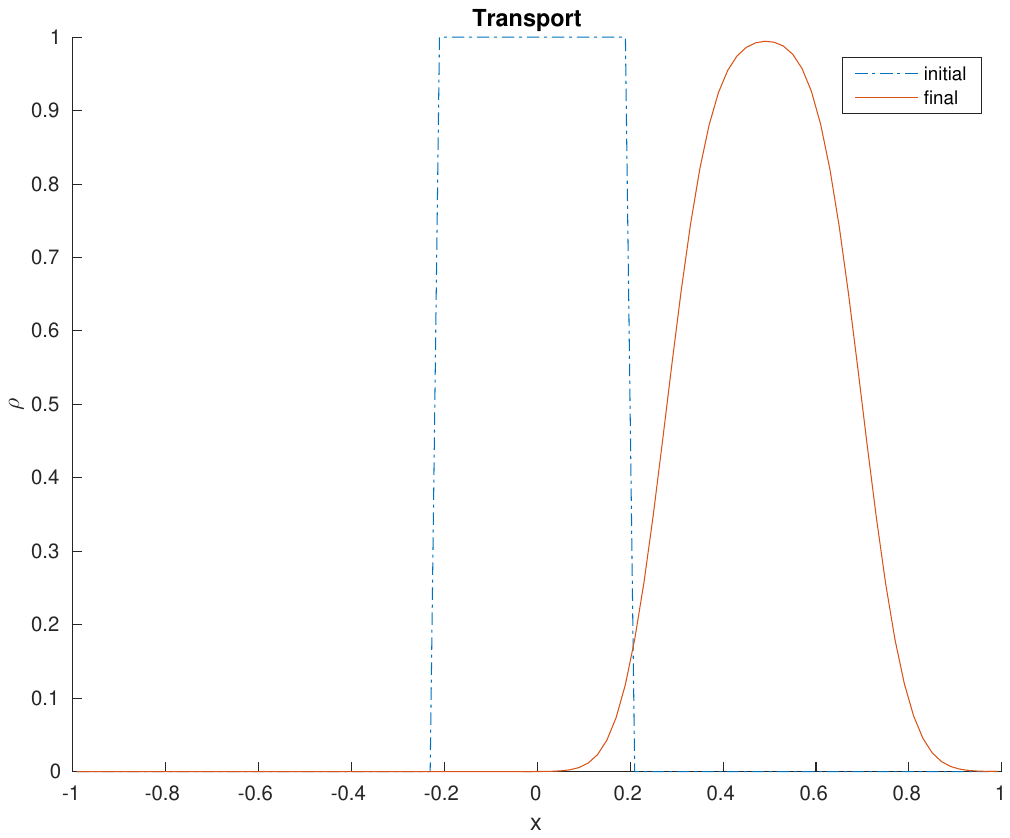} \\
A & B \\
\end{tabular}
\caption{A: Left transport of a square wave $v \, \Delta t / \Delta x = -1$.
B: Right transport of a square wave with $v = 0.4167$. (See {\tt Transport.m})}
\label{Transport Figure}
\end{center}
\end{figure}

So consider the edges of the cells and compute the amount of material
being transferred between the neighboring cells, that is for each time
step $n$, for all $i$ compute the discrete solution as follows:
\begin{align*}
\text{if } v_i \geq 0 \text{ then} \quad & 
\rho^{n+3/2}_{i-\half} = \rho^{n+3/2}_{i-\half}
- v_i \frac{\Delta t}{\Delta x} \, \rho^{n+\half}_{i-\half} \,;\\
& \rho^{n+3/2}_{i+\half} = \rho^{n+3/2}_{i+\half}
+ v_i \frac{\Delta t}{\Delta x} \, \rho^{n+\half}_{i-\half} \,; \\
\text{if } v_i \leq 0 \text{ then}\quad &
\rho^{n+3/2}_{i-\half} = \rho^{n+3/2}_{i-\half}
-  v_i \frac{\Delta t}{\Delta x} \, \rho^{n+\half}_{i+\half} \,;\\
& \rho^{n+3/2}_{i+\half} = \rho^{n+3/2}_{i+\half}
+ v_i \frac{\Delta t}{\Delta x} \, \rho^{n+\half}_{i+\half} \,.
\end{align*}
If $v_i$ is positive then this removes some material from cell
$i-\half$ and put it into cell $i+\half$ and conversely if $v_i$
is negative. If $V = \max(|v_i|)$ then the most material that can
be removed from cell $i-\half$ is 
\[
V \,  \frac{\Delta t}{\Delta x} \rho^{n+\half}_{i-\half} \,,
\]
so to keep $\rho \geq 0$ it must be that
\[
V \,  \frac{\Delta t}{\Delta x}  \leq 1 \,.
\]
An interesting feature of this algorithm is that for 
$v_i \,  \Delta t / \Delta x = \pm 1$ it gives an exact solution
of solution as shown in Figure \ref{Transport Figure}.
This is an upwind scheme for velocities that change direction that
keeps that preserves $\rho \geq 0$ and conserves the amount material
being transported. As done in {\tt Transport.m} this scheme can be
implemented with out the conditional in the update loop.

\begin{figure}
\begin{center}
\begin{tabular}{cc}
   \includegraphics[width=2.75in,trim = 0 0 0 350]{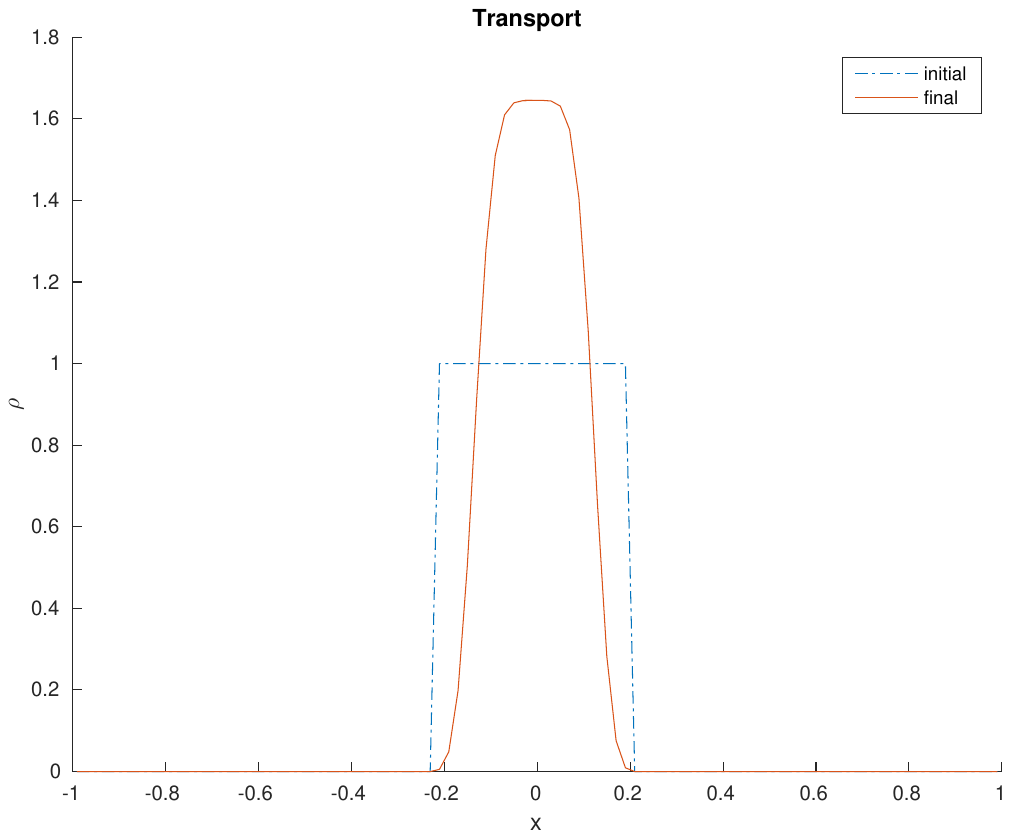} &
   \includegraphics[width=2.75in,trim = 0 0 0 350]{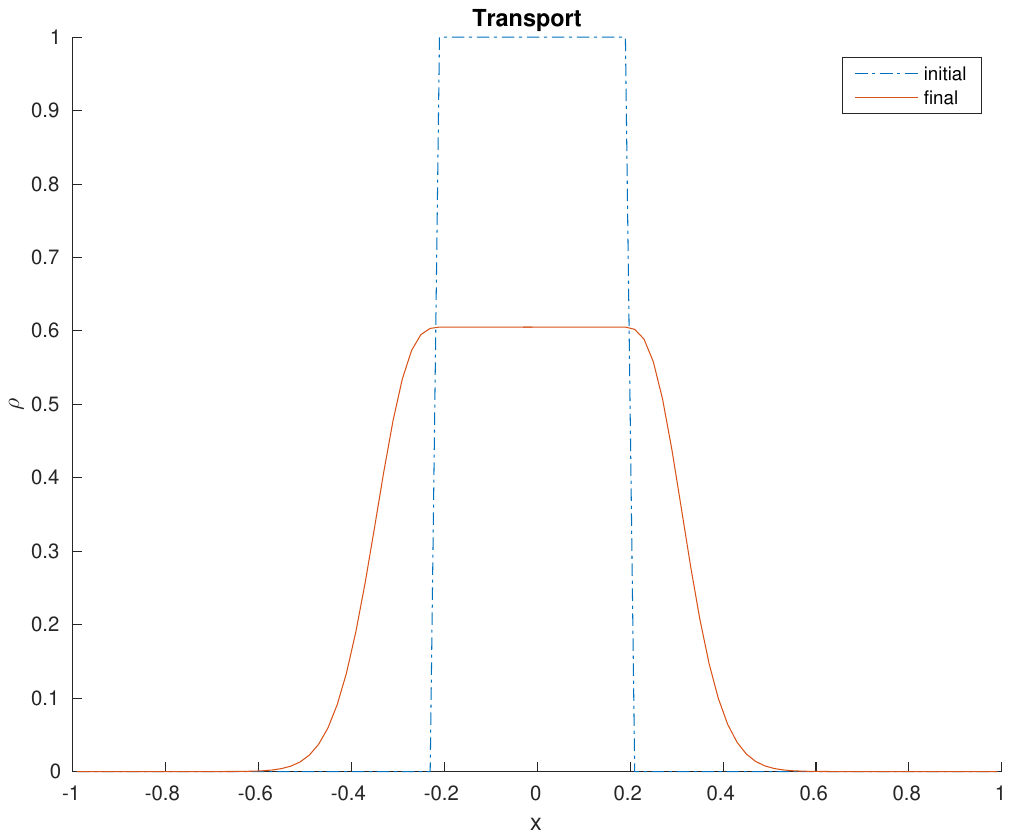} \\
A & B \\
\end{tabular}
\caption{A: Collapse with $v = -x$. B: Expand with $v = x$.
(See {\tt Transport.m})}
\label{Colapse Expand Figure}
\end{center}
\end{figure}

Not all discretizations preserve positive solutions, for example
the Lax-Wendroff, Richtmyer, and MacCormac schemes do not for linear
equations (see {\tt Lax-Wendroff-Positive.nb}).
This can also be seen by by choosing initial data $f_i$ that are
all zero except for one $i$ where $f_i = 1$. For linear equations
the Richtmyer and MacCormac schemes produce the same solution as
the Lax-Wendroff scheme.

\subsection{Diffusion}

The diffusion equation in one dimension is given by
\[
\frac{\partial \rho}{\partial t} =
\frac{\partial}{\partial x}  D \, \frac{\partial \rho}{\partial x} \,,
\]
where $\rho = \rho(x,t)$ is the heat density $D = D(x) \geq 0$ is the
diffusion coefficient. For this discussion $t \geq 0$ and $\rho$ is smooth
and zero for large values of $|x|$. Then integrating the differential
equation gives
\[
\int_{-\infty}^{\infty} \rho(x,t) \, dx = 0 \,.
\]
If $\rho(x,0) \geq 0$ then the solution of the equation is given by
convolution with a Gaussian so then $\rho(x,t)\geq 0 $ for $t \geq 0$.

The standard forward time center space finite difference discretization
of this equation is given by
\[
\frac{\rho_{i+\half}^{n+\half} - \rho_{i+\half}^{n-\half}}{\Delta t} =
\frac{1}{\Delta x} \left(
D_{i+1} \frac{\rho_{i+\thalf}^{n-\half} - \rho_{i+\half}^{n-\half}}{\Delta x}-
D_{i} \frac{\rho_{i+\half}^{n-\half} - \rho_{i-\half}^{n-\half}}{\Delta x}
\right)
\]
or in computational form
\[
\rho_{i+\half}^{n+\half} =
\rho_{i+\half}^{n-\half} +
\frac{\Delta t}{\Delta x^2}
\left(
 D_{i+1} \, \rho_{i+\thalf}^{n-\half}
-\left( D_{i+1} + D_{i} \right)  \, \rho_{i+\half}^{n-\half}
+D_{i}   \, \rho_{i-\half}^{n-\half}
\right).
\]
This algorithm will preserve positive solutions for
\[
\left( D_{i+1} + D_{i} \right) \, \frac{\Delta t}{\Delta x^2} \leq 1 \,,
\]
which is the standard stability constraint for this discretization.

\end{appendix}

\newpage \clearpage
\end{document}